\documentclass[a4paper,10pt]{amsart}
\usepackage{amsfonts}
\usepackage[colorlinks=true,linkcolor=blue,draft=false]{hyperref}
\usepackage{etoolbox}
\usepackage{amsmath}
\usepackage{amssymb}
\usepackage{amsthm}
\usepackage{dsfont}
\usepackage{pifont}
\usepackage{mathrsfs}
\usepackage{graphicx}
\usepackage{enumerate}
\usepackage[all]{xy}
\usepackage{geometry}
\usepackage{amsfonts}
\usepackage{fancyhdr}
\usepackage{tikz}
\usepackage{tikz-cd}
\usetikzlibrary{arrows}
\usepackage[section]{algorithm}
\usepackage{algorithmic}

\usepackage{dsfont}

\usepackage{amsaddr}

\usepackage[backend=biber,sorting=nty,maxnames=5]{biblatex}
\newcommand{\Aa}{\mathcal{A}}
\newcommand{\Bb}{\mathcal{B}}
\newcommand{\Cc}{\mathcal{C}}

\newcommand{\Gg}{\mathcal{G}}
\newcommand{\Hh}{\mathcal{H}}

\newcommand{\Ll}{\mathcal{L}}

\newcommand{\Nn}{\mathcal{N}}

\newcommand{\Rr}{\mathcal{R}}
\newcommand{\Ss}{\mathcal{S}}

\newcommand{\Uu}{\mathcal{U}}

\newcommand{\C}{\mathbb{C}}

\newcommand{\R}{\mathbb{R}}
\renewcommand{\S}{\mathfrak{S}}

\newcommand{\Z}{\mathbb{Z}}
\newcommand{\D}{\mathbb{D}}

\DeclareMathOperator\Ker{Ker}

\DeclareMathOperator\GL{GL}

\DeclareMathOperator\Div{Div}

\DeclareMathOperator\Ob{Ob}

\DeclareMathOperator\LL{LL}
\DeclareMathOperator\clbl{clbl}
\DeclareMathOperator\cc{cc}

\DeclareMathOperator\PC{PC}
\DeclareMathOperator\SPC{SPC}

\DeclareMathOperator\SCPC{SCPC}
\DeclareMathOperator\sw{sw}

\DeclareMathOperator\olbl{olbl}

\DeclareMathOperator\loc{loc}

\newcommand{\LLL}{\bbar{\LL}}

\newcommand{\xdownarrow}[1]{%
  {\left\downarrow\vbox to #1{}\right.\kern-\nulldelimiterspace}
}

\newcommand{\infl}{\ar@{{(}->}}
\newcommand{\defl}{\ar@{->>}}

\newcommand{\intv}[1]{[\![#1]\!]}
\newcommand{\muls}[1]{ \left\{\kern-0.6em \left\{ #1\right\}\kern-0.6em \right\} }

\newcommand{\epsi}{\varepsilon}

\newcommand{\pphi}{\varphi}
\newcommand{\sous}{\backslash}
\newcommand{\RE}{\mathsf{R}}

\newcommand{\bbar}[1]{\overline{#1}}

\newcommand{\ttilde}[1]{\widetilde{#1}}
\newcommand{\nit}[1]{\underline{#1}}

\theoremstyle{plain}
\newtheorem{prop}{Proposition}[section]
\newtheorem{prop-def}[prop]{Proposition-Definition}
\newtheorem{lem}[prop]{Lemma}
\newtheorem{theo}[prop]{Theorem}
\newtheorem{cor}[prop]{Corollary}

\newtheorem*{theo*}{Theorem}
\newtheorem*{cor*}{Corollary}
\newtheorem*{prop*}{Proposition}
\newtheorem*{fact*}{Fact}

\theoremstyle{definition}
\newtheorem{definition}[prop]{Definition}

\theoremstyle{remark}
\newtheorem{rem}[prop]{Remark}

 \normalbaselines
\title[Parabolic subgroups of complex braid groups: the remaining case]{Parabolic subgroups of complex braid groups: the remaining case}
\author{Owen Garnier}
\address{LAMFA, Université de Picardie Jules Verne, CNRS UMR 7352,\\ 33, rue Saint-Leu, 80000, Amiens, France.}
\email{o.garnier@u-picardie.fr}
\date{\today}

\subjclass[2020]{Primary 20F36, Secondary 20F55}
\keywords{Complex braid groups, Parabolic subgroups, Garside category, Curve complex.}
\begin{document}

\begin{abstract}
Recently, Marin and Gonz\'alez-Meneses introduced a class of ``parabolic'' subgroups for generalized braid groups associated to arbitrary complex reflection groups. Using notably Garside group structures on these generalized braid groups, they proved general results on parabolic subgroups for all cases but one. This last case is that of the complex braid group $B(G_{31})$, which has no Garside group structure known so far, but instead a Garside groupoid structure. Using this Garside groupoid structure, we complete the results of Marin and Gonz\'alez-Meneses by proving their main theorems for the parabolic subgroups of the complex braid group $B(G_{31})$.
\end{abstract}

\maketitle
\tableofcontents
\addtocontents{toc}{\protect\setcounter{tocdepth}{1}}
\section{Introduction}
Let $W$ be a complex reflection group, and let $B(W)$ be the generalized braid group associated to it, as in \cite{bmr}. In \cite{paratresses}, Marin and Gonz\'alez-Meneses provide a purely topological definition of a parabolic subgroup of $B(W)$, which is (in particular) a subgroup isomorphic to the braid group of some parabolic subgroup of $W$. 

When $B(W)$ is endowed with a ``natural'' Garside group structure (e.g. dual braid monoid, Artin-Tits monoid), they show that the collection of parabolic subgroups of $B(W)$ coincides with the collection of parabolic subgroups as defined by Godelle in \cite{parabolicgarside}. This allows the authors to use Garside-theoretic arguments to prove structural results on parabolic subgroups, most notably that they form a lattice under inclusion.

This approach leaves out the irreducible groups $W$ for which $B(W)$ do not admit (to our knowledge) a convenient Garside group structure. In \cite{paratresses}, the authors are able to use topological methods to extend their results to almost every such groups. In this article, we settle the only remaining case, that is $B(G_{31})$. More precisely, we show the two following theorems.

\begin{theo*}(Theorem \ref{theo:parabolic_closure}) For every element $x\in B(G_{31})$, there exists a unique minimal parabolic subgroup $\PC(x)$ of $B(G_{31})$ containing it, and we have $\PC(x^m)=\PC(x)$ for every $m\neq 0$.
\end{theo*}

This is the analogue of of \cite[Theorem 1.1]{paratresses} in the case of $B(G_{31})$. The group $\PC(x)$ is called the \nit{parabolic closure} of $x$.

\begin{theo*}(Theorem \ref{theo:intersection_parabolic_subgroups}) If $B_1,B_2$ are two parabolic subgroups of $B(G_{31})$, then $B_1\cap B_2$ is a parabolic subgroup of $B(G_{31})$. More generally, every intersection of a family of parabolic subgroups is a parabolic subgroup.
\end{theo*}

This is the analogue of \cite[Theorem 1.2]{paratresses} in the case of $B(G_{31})$.
 
For any complex reflection group $W$, the irreducible parabolic subgroups of $B(W)$ are the vertices of a graph $\Gamma$, called the \nit{curve graph}, where two parabolic subgroups $B_1,B_2$ are adjacent if and only if $B_1\neq B_2$ and either $B_1\subset B_2$, or $B_2\subset B_1$, or $B_1\cap B_2=[B_1,B_2]=1$. The clique complex of $\Gamma$, that is, the flag complex made of all the simplices whose edges belong to the curve graph, is endowed with a faithful action of $B(W)/Z(B(W))$ by \cite[Proposition 2.10]{paratresses}. This construction generalizes the usual curve complex associated to the usual braid group.

Another result proven in \cite{paratresses} in all cases except $B(G_{31})$ is a convenient characterization of adjacency in the curve graph. Again we prove that this result also holds for $B(G_{31})$. Recall from \cite{bmr} and \cite{beskpi1} that the center of an irreducible complex braid group is infinite cyclic. We can associate to an irreducible parabolic subgroup $B_0$ of $B(G_{31})$ the positive generator $z_{B_0}$ of its center (here, positive means that the image of $z_{B_0}$ under the natural map $B(G_{31})\to \Z$ is positive). We have

\begin{theo*}(Theorem \ref{theo:characterization_of_adjacency}) If $B_0\subset B(G_{31})$ is an irreducible parabolic subgroup, then $B_0=\PC(z_{B_0})$. Moreover, if $B_1,B_2\subset B(G_{31})$ are two irreducible parabolic subgroups, then an element $g\in B(G_{31})$ satisfies $(B_1)^g=B_2$ if and only if $(z_{B_1})^g=z_{B_2}$, and $B_1,B_2$ are adjacent in $\Gamma$ if and only if $z_{B_1}$ and $z_{B_2}$ are distinct and commute.
\end{theo*}

The proof of these three results goes as follows. As we said, the group $B(G_{31})$ has no known Garside group structure. However, Bessis defined in \cite[Section 11]{beskpi1} a Garside groupoid $\Bb_{31}$ which is equivalent to $B(G_{31})$. This groupoid is a particular case of a broader construction: Consider $W<\GL(V)$ an irreducible complex reflection group. Assume that $W$ is well-generated, that is, generated by $\dim V$ reflections. Consider $g$ a regular element of $W$ in the sense of Springer (possibly for the eigenvalue $1$). The centralizer $W_g$ of $g$ in $W$ is known to also be a complex reflection group. The \nit{Springer groupoid} \cite[Definition 11.23]{beskpi1} associated to $W_g$ and $W$ is a Garside groupoid $\Bb$, which is equivalent to $B(W_g)$ considered as a groupoid with one object.

Our method consists in using the Springer groupoid in the same way that Marin and Gonz\'alez-Meneses use the dual braid monoid associated to a well-generated group. In order to do this, we first define a Garside-theoretic notion of standard parabolic subgroupoid of a Springer groupoid. Our definition is more restrictive than the one given in the reference book \cite{ddgkm} on Garside categories. However, the following result, which is the goal of Section 1, proves that our definition of standard parabolic subcategories is a good analogue of the notion of standard parabolic subgroup, in that parabolic subgroups are exactly the conjugates (in the groupoid sense) of the standard parabolic subgroups. Recall that the equivalence between $B(W_g)$ and $\Bb$ allows us to identify $B(W_g)$ with the group $\Bb(u,u)$ of endomorphisms of some fixed object $u$ in $\Bb$.
\begin{theo*}(Theorem \ref{theo:parabolics_up_to_conj})
Let $B_0\subset B(W_g)$ be a subgroup. The following assertions are equivalent
\begin{enumerate}[(1)]
\item The subgroup $B_0$ is a parabolic subgroup of $B(W_g)$.
\item There is a standard parabolic subgroupoid $\Bb_0$ of $\Bb$ and a morphism $f\in \Bb(u,v)$, where $v\in \Ob(\Bb_0)$, such that conjugation by $f$ in $\Bb$ identifies the image of $B_0$ in $\Bb(u,u)$ with $\Bb_0(v,v)$.
\end{enumerate}
\end{theo*}

This general theorem specializes to the case where $W_g=G_{31}$ and $W=G_{37}$, giving us a complete description of the parabolic subgroups of $B(G_{31})$ in Garside-theoretic terms, using the Springer groupoid $\Bb_{31}$. The second section is dedicated to this specific case. This restriction allows us to prove several results by direct computations on the (finite set) of standard parabolic subgroupoids of $\Bb_{31}$. In particular, we prove a complete description of the parabolic subgroups of $B(G_{31})$ up to conjugacy in terms of braid diagrams. By \cite[Theorem 4.2 and Section 4.3.1]{springercat}, the group $B(G_{31})$ admits the following presentation
\[B(G_{31})=\left\langle s,t,u,v,w\left| \begin{array}{l} st=ts,~vt=tv,~wv=vw,\\ suw=uws=wsu,\\svs=vsv,~vuv=uvu,~utu=tut,~twt=wtw\end{array}\right. \right\rangle.\]

\begin{theo*}(Theorem \ref{theo:lattice_of_parabolics_of_B31}) The parabolic subgroups of $B(G_{31})$ are given up to conjugacy by the admissible subdiagrams of its braid diagram. That is, parabolic subgroups of $B(G_{31})$ are exactly the conjugates in $B(G_{31})$ of the following subgroups
\[\begin{array}{rclrclrcl} \langle \varnothing\rangle&=&\{1\},& \langle s \rangle&\simeq& B(A_1),& \langle t,v\rangle&\simeq& B(A_1\times A_1),\\
\langle s,v \rangle&\simeq& B(A_2),&\langle s,u,w \rangle&\simeq& B(G(4,2,2)),&\langle s,t,v\rangle&\simeq &B(A_1\times A_2),\\
\langle s,u,v,w\rangle&\simeq& B(G(4,2,3)),&\langle t,u,v\rangle&\simeq& B(A_3),&\langle s,t,u,v,w \rangle&= &B(G_{31}),\end{array}\]
where $A_n$ denotes the complex reflection group $G(1,1,n+1)$. For each such parabolic subgroup $\langle S\rangle$, a presentation is given by taking the relations in the above presentation of $B(G_{31})$ involving only elements of $S$ (in particular, we obtain the given isomorphism type for each group).
\end{theo*}

\subsection*{Acknowledgments.} Work on this article began during a research stay in Sevilla in autumn 2022, where I worked with Juan Gonz\'alez-Meneses. I thank him for his hospitality and for his very useful insights, especially on Lemma \ref{lem:sundial}. I also thank my PhD advisor Ivan Marin, for guiding this work with his numerous remarks and profound questions. This work is part of my PhD thesis.
\addtocontents{toc}{\protect\setcounter{tocdepth}{2}}

\section{Parabolic subgroupoids of Springer groupoids}
The goal of this section is to introduce both a topological and combinatorial concept of parabolic subgroupoid of the Springer groupoid associated to a regular centralizer in a well-generated complex reflection group. These two notions should coincide in the sense that a topological parabolic subgroupoid should be conjugate to a combinatorial (standard) parabolic subgroupoid. To construct such notions, we need several reminders about the topological and combinatorial construction of Springer groupoids.

\subsection{Complex braid groups, regular braids}

Let $V$ be a finite dimensional complex vector space. A subgroup $W<\GL(V)$ is a \nit{complex reflection group} if it is finite and generated by (pseudo-)reflections, i.e. elements in $\GL(V)$ whose fixed subspace is a hyperplane. 

Let $W<\GL(V)$ be a complex reflection group, and let $v\in V$. By a classical theorem of Steinberg \cite[Theorem 9.44]{lehrertaylor}, the stabilizer of $v$ in $W$ is generated by the reflections of $W$ which fix $v$. Such a stabilizer is called a \nit{parabolic subgroup} of $W$.

Let $\Aa$ be the set of hyperplanes associated to reflections in $W$ (reflecting hyperplanes). We consider the intersection lattice $\Ll(\Aa):=\{\bigcap_{H\in A}H~|~ A\subset \Aa\}$. Elements of $\Ll(\Aa)$ are called \nit{flats}. For $L\in \Ll(\Aa)$, we denote by $L^0$ the complement in $L$ of the flats strictly induced in $L$. The sets $L^0$ for $L\in \Ll(\Aa)$ form a stratification of the space $V$. Since $W$ acts on flats, this stratification induces a stratification of the quotient space $V/W$, that we call the \nit{discriminant stratification}.

Let $X$ be the complement in $V$ of the union of the reflecting hyperplanes of $W$. We have that $V=\bigcap_{H\in \varnothing}H$ is a flat, and $X=V^0$ with the above notation. The \nit{braid group} $B(W)$ attached to $W$ is then defined as the fundamental group of $X/W$. 

A complex reflection group $W<\GL(V)$ is \nit{irreducible} if there are no subspaces of $V$ which are globally $W$-invariant, apart from $\{0\}$ and $V$. Irreducible complex reflection groups were classified by Shephard and Todd in \cite{shetod}. From now on we restrict our attention to irreducible complex reflection groups.

An element $g\in W$ is called \nit{regular} for the eigenvalue $\zeta$ if $g$ admits a $\zeta$-eigenvector in $X$. An element $g\in W$ which is regular for the eigenvalue $\zeta_d:=\exp(\frac{2i\pi}{d})$ will be called $d$-regular. If such an element exist, we say that $d$ is a \nit{regular number} for $W$. Note that, as the identity element in $W$ is $1$-regular by definition, $d=1$ is always a regular number for $W$. The following Theorem is due to Springer \cite{sprireg} (statement (1)) and to Denef-Loeser \cite{denefloeser} (statement (2)).

\begin{theo}\label{theo:regular_elements}\cite[Theorem 1.9]{beskpi1} Let $W<\GL(V)$ be an irreducible complex reflection group, and let $g\in W$ be a $d$-regular element. 
\begin{enumerate}
\item The group $W_g:=C_W(g)$ acts as a complex reflection group on $V_g:=\Ker(g-\zeta)$.
\item We denote by $X_g$ the complement in $V_g$ of the reflecting hyperplanes of $W_g$. The natural inclusion $V_g\to V$ induces a homeomorphism of topological pairs $(X_g/W_g,V_g/W_g)\simeq ((X/W)^{\mu_d},(V/W)^{\mu_d})$. In particular, $B(W_g)$ is isomorphic to the fundamental group of $(X/W)^{\mu_d}$.
\end{enumerate}
\end{theo}

By the Chevalley-Shephard-Todd Theorem \cite[Theorem V.5.4]{boulie}, the subalgebra $S(V^*)^W$ of the symmetric algebra $S(V^*)$ is a polynomial algebra. More precisely, one can consider a system of homogeneous polynomials $(f_1,\ldots,f_n)$ in $S(V^*)$ (where $n=\dim_\C V$) such that $S(V^*)^W=\C[f_1,\ldots,f_n]$. We call $(f_1,\ldots,f_n)$ a \nit{system of basic invariants} for $W$. The sequence of degrees of a system of basic invariants (in nondecreasing order) is independent of the choice of the system. We call this sequence the \nit{degrees} of the complex reflection group $W$.

Let $f=(f_1,\ldots,f_n)$ be a system of basic invariants for $W$. Geometrically, the isomorphism between $S(V^*)^W$ and a polynomial algebra induced by $f$ allows us to identify the orbit space $V/W$ with the affine space $\C^n$. Under this identification, the complement $\Hh$ of $X/W$ in $V/W$ becomes an algebraic hypersurface, called the \nit{discriminant hypersurface}. We denote by $\Delta_f$ the reduced equation of $\Hh$ in $\C[X_1,\ldots,X_n]$.

Let $W< \GL(V)$ be irreducible and \nit{well-generated} (i.e. generated by a family of $\dim_\C V$ reflections). By \cite[Theorem 2.4]{beskpi1}, the highest degree $h$ of $W$ is unique (i.e. it appears exactly once in the sequence of degrees) and it is a regular number for $W$.
\subsection{Dual braid monoid}
Our main references are \cite{beskpi1} and \cite{besgar}. We fix $W$ an irreducible well-generated complex reflection group. We denote by $\Rr$ the set of all its reflections, and by $h$ its highest degree. 

The \nit{reflection length} of an element $w\in W$ is defined as the minimal $k\geqslant 0$ such that $w$ is written $w=r_1\cdots r_k$ with $r_i\in \Rr$, and is denoted by $\ell_\Rr(w)$. We take the convention that $\ell_\Rr(1)=0$. For $w\in W$ and $m\geqslant 0$, we denote by $D_m(w)$ the set of \nit{length-additive factorizations} of $w$ as a product of $m$ elements
\[D_m(w):=\left\{(u_1,\ldots,u_m)~|~ u_1\cdots u_m=w \text{ and }\ell_\Rr(u_1)+\cdots+ \ell_\Rr(u_m)=\ell_\Rr(w)\right\}.\]
These sets are important for defining and handling Springer groupoids. We also define the set $D_\bullet(w):=\bigsqcup_{m\geqslant 1} D_m(w)$ of length-additive factorizations of $w$.

The reflection length is used to define two partial orderings on $W$, setting $a\preccurlyeq b$ if $\ell_\Rr(a)+\ell_\Rr(a^{-1}b)=\ell_\Rr(b)$, and $a\preccurlyeq^{\Lsh} b$ if $\ell_\Rr(ba^{-1})+\ell_\Rr(a)=\ell_\Rr(b)$. By definition, for $w\in W$, the set $[1,w]=\{a\in W~|~1\preccurlyeq a\preccurlyeq w\}$ is in one-to-one correspondence with $D_2(w)$ by $x\mapsto (x,x^{-1}w)$. Note that, as the set $\Rr$ is invariant under conjugacy, we have $x\preccurlyeq y$ if and only if $x\preccurlyeq^\Lsh y$. That is, every element of $W$ is \nit{balanced}.

We fix a \nit{Coxeter element} $c\in W$ (a $h$-regular element). The \nit{dual braid monoid} $M(c)$ attached to the data $(W,c)$ is defined by generators and relations, taking a formal copy of $[1,c]$ as generators, and endowing it with relations $s\cdot t=u$ whenever $(s,t)\in D_2(u)$. 

Moreover, the set $D_m(c)$ is endowed with an map $\tau$, defined by $\tau(x_1,\ldots,x_m):=(x_2,\ldots,x_m,x_1^c)$. We can then define $D_m^n(c):=\{x\in D_m(c)~|~\tau^n(x)=x\}$.

\subsection{Lyashko-Looijenga maps, Springer groupoids}
Let $W<\GL(V)$ be an irreducible well-generated complex reflection group with highest degree $h$, and let $d$ be a regular number for $W$. We set
\[p:=\frac{d}{d\wedge h}\text{~and~}q:=\frac{h}{d\wedge h}\]
where $d\wedge h$ is the gcd of $d$ and $h$.

\subsubsection{Topological definition of Springer groupoids}Let $n$ denote the dimension of $V$. We fix a system of basic invariants $f=(f_1,\ldots,f_n)$ for $W$.
\begin{definition}\cite[Definition 7.25]{beskpi1} Let $x=W.v\in V/W$. The multiset $\LLL(x)$ is defined as the solutions of the univariate polynomial
\[P(T):=\Delta_f(f_1(v),\ldots,f_{n-1}(v),f_n(v)+T)\in \C[T].\]
By construction, $\LLL(x)$ lies in $E_n:=\C^n/\S_n$. The map $\LLL:V/W\to E_n$ is called the (extended) \nit{Lyashko-Looijenga morphism}.
\end{definition}
Recall that $(f_1,\ldots,f_n):V\to \C^n$ induces a homeomorphism $V/W\simeq \C^n$. Furthermore, since $f_1,\ldots,f_n$ are polynomials, $P(T)$ is also a polynomial. Recall also that the elementary symmetric polynomials provide an identification $E_n\simeq \C^n$ (they form a system of basic invariant for the complex reflection group $\S_n<\GL_n(\C)$). 

Since $W$ is well-generated by assumption, the map $\LLL$ is a ramified covering map by \cite[Theorem 5.3]{beskpi1}. By construction, we have $0\in \LLL(x)$ if and only if $x\in \Hh$. If we denote by $E_n^\circ\subset E_n$ the subset of multisets not containing $0$, then for $x\in V/W$, we have $x\in X/W$ if and only if $\LLL(x)\in E_n^\circ$. 

The following proposition is very useful for constructing paths in $V/W$.

\begin{prop}\label{prop:lift_lll}\cite[Remark 7.21]{beskpi1} Let $\gamma:[0,1]\to E_n$ be a path such that points are not unmerged as $t$ increases (they can be merged). If $x_0\in V/W$ is such that $\LLL(x_0)=\gamma(0)$, then there exists a unique lift $\ttilde{\gamma}:[0,1]\to V/W$ of $\gamma$ such that $\ttilde{\gamma}(0)=x_0$.
\end{prop}

The open subset $\Uu:=\{x\in V/W~|~ \LLL(x)\cap i\R_{\geqslant 0}=\varnothing\}$ of $X/W$ is contractible \cite[Lemma 6.3]{beskpi1}, thus $\pi_1(X/W,\Uu)$ is well-defined as the set of classes of paths from some point of $\Uu$ to some other one, up to a homotopy leaving the endpoints in $\Uu$. Following \cite[Definition 6.4]{beskpi1}, we set $B(W)=\pi_1(X/W,\Uu)$ from now on.

\begin{definition}\cite[Definition 11.24]{beskpi1} A \nit{circular tunnel} is a couple $T=(x,L)\in \Uu\times [0,2\pi/h]$ such that $e^{iL}x\in \Uu$. It is identified with the path $\gamma_T:[0,1]\to X/W$ sending $t$ to $e^{itL}x$.
\end{definition}
By definition of $\pi_1(X/W,\Uu)$, the path $\gamma_T$ associated to a circular tunnel $T$ induces a well-defined element of $B(W):=\pi_1(X/W,\Uu)$. Of course, several distinct circular tunnels may represent the same element of $B(W)$. In particular, we have the following useful lemma, which is a variation of \cite[Lemma 6.15]{beskpi1}.

\begin{lem}(Hurwitz rule)\label{lem:hurwitzrule}
Let $\lambda:[0,1]\to [0,\frac{2\pi}{h}]$ and let $\gamma:[0,1]\to \Uu$ be continuous paths.\\ If $T_t:=(\gamma(t),\lambda(t))$ is a circular tunnel for all $t\in [0,1]$, then all the $T_t$ for $t\in [0,1]$ represent the same element of $B(W)$.
\end{lem}
\begin{proof}
Let $t\in [0,1]$. Defining $H(r,s):=e^{is\lambda(rt)}\gamma(rt)$ for $r,s\in [0,1]$ yields a homotopy between the paths $H(0,s):s\mapsto e^{i\lambda(0)}\gamma(0)$ and $H(1,s):s\mapsto e^{is\lambda(t)}\gamma(t)$ in $X/W$. By assumptions, the endpoints of this homotopy, namely $H(r,0)=\gamma(rt)$ and $H(r,1)=e^{i\lambda(rt)}\gamma(rt)$, lie in $\Uu$. Thus $T_0$ and $T_t$ represent the same element of $B(W)$ by definition.
\end{proof}

In particular, this lemma allows us to see that for any $x\in \Uu$, the couple $(x,\frac{2\pi}{h})$ is a circular tunnel, which represents the same element $\Delta$ of $B(W)$. Indeed, $(x,\frac{2\pi}{h})$ is always a circular tunnel since $\LLL(e^{\frac{2i\pi}{h}}x)=\LLL(x)$ for $x\in \Uu$. Then, for $x,y\in \Uu$, we can consider a path $\gamma:[0,1]\to \Uu$ from $x$ to $y$ since $\Uu$ is path connected. We then have that $(x,\frac{2\pi}{h})$ and $(y,\frac{2\pi}{h})$ represent the same element of $B(W)$ by Lemma \ref{lem:hurwitzrule}. Moreover, we have the following theorem.

\begin{theo}\cite[Proposition 8.5 and Lemma 11.10]{beskpi1} The projection map $B(W)\to W$ sends $\Delta$ to a Coxeter element $c$ of $W$. It also induces a one-to-one correspondence between homotopy classes of circular tunnels and the set $[1,c]$. This one-to-one correspondence in turn induces an isomorphism between $B(W)$ and the enveloping group of the dual braid monoid $M(c)$.
\end{theo}

From now on, we fix $c$ to be the image of $\Delta$ in $B(W)$. Using this theorem, we can define the image of a circular tunnel in $W$ as the image of its homotopy class in $B(W)$ under the projection map $B(W)\to W$. The image of a circular tunnel in $W$ belongs to $[1,c]$.

We aim to describe the fibers of the Lyashko-Loooijenga morphism on $V/W$, in the spirit of \cite[Proposition 11.13]{beskpi1} (which covers the case of $X/W$). In \cite[Proposition 11.13]{beskpi1} this is done by defining a notion of cyclic label of an element $x$ of $X/W$ which, coupled with $\LLL(x)$, completely characterizes $x$. We will extend the definition of the cyclic label to points belonging to $\Hh$. In \cite{beskpi1}, the notion of cyclic label of an element of $X/W$ is defined using an earlier notion of label \cite[Definition 6.9]{beskpi1}. We provide a direct definition here.

Let $\kappa\in E_n^\circ$. The arguments of the points composing $\kappa$ are seen in $[\frac{-3\pi}{2},\frac{\pi}{2}[$. Following \cite[Definition 11.8]{beskpi1}, the points composing $\kappa$ are ordered by decreasing argument. Points with identical arguments are ordered by increasing modulus. Since we consider arguments in $[\frac{-3\pi}{2},\frac{\pi}{2}[$, points on the vertical half-line $i\R_+$ always come last, as in the following example

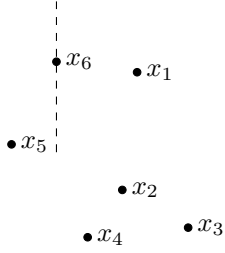
\begin{center}
\begin{tikzpicture}
\draw [dashed] (0,0)--(0,2);

\node[draw,circle,inner sep=1pt,fill] (P1) at (1.5*0.7071,1.5*0.7071) {};
\node[right] at (1.5*0.7071,1.5*0.7071) {$x_1$};

\node[draw,circle,inner sep=1pt,fill] (P2) at (-30:1) {};
\node[right] at (P2) {$x_2$};

\node[draw,circle,inner sep=1pt,fill] (P2) at (-30:2) {};
\node[right] at (P2) {$x_3$};

\node[draw,circle,inner sep=1pt,fill] (P2) at (-70:1.2) {};
\node[right] at (P2) {$x_4$};

\node[draw,circle,inner sep=1pt,fill] (P2) at (-190:0.6) {};
\node[right] at (P2) {$x_5$};

\node[draw,circle,inner sep=1pt,fill] (P2) at (90:1.2) {};
\node[right] at (P2) {$x_6$};
\end{tikzpicture}
\end{center}

The \nit{cyclic support} of $\kappa$ is then defined as the sequence $(\kappa_1,\ldots,\kappa_m)$ of distinct points in $\kappa$, ordered as above \cite[Definition 11.8]{beskpi1}. If $\kappa$ contains $0$, then we extend this definition by defining $(\kappa_1,\ldots,\kappa_m,0)$ as the cyclic support of $\kappa$, where $(\kappa_1,\ldots,\kappa_m)$ is the cyclic support of $\kappa\setminus \{0\}$. These definitions are extended to $V/W$ by defining the cyclic support of $x\in V/W$ as that of $\LLL(x)$. 

Let $x\in X/W$, with cyclic support $(\kappa_1,\ldots,\kappa_m)$. For $j\in \intv{1,m}$, we set $\theta_j:=\frac{\pi}{2}-\arg(\kappa_j)\in ]0,2\pi]$ (again, we see the argument of $\kappa_j$ in $[\frac{-3\pi}{2},\frac{\pi}{2}[$), so that $e^{i\theta_j}\kappa_j\in i\R_+$. The sequence $(\theta_1,\ldots,\theta_m)$ is the \nit{cyclic argument} of $x$ \cite[Definition 11.8]{beskpi1}.

Assume at first that $0<\theta_1<\theta_2<\ldots<\theta_m<2\pi$, and let $\theta\in ]0,\frac{2\pi}{h}]$. By definition, the couple $(x,\theta)$ is a circular tunnel if and only if $\theta$ is different from all the $\frac{\theta_j}{h}$ for $j\in \intv{1,m}$. For $j\in \intv{0,m}$, the Hurwitz rule (Lemma \ref{lem:hurwitzrule}) gives that the class of the circular tunnel $(x,\theta)$ does not depend on $\theta$ such that $\frac{\theta_j}{h}<\theta<\frac{\theta_{j+1}}{h}$ (with the convention that $\theta_0=0$ and $\theta_{m+1}=2\pi$). We denote by $p_j$ the image of this circular tunnel in $W$. The Hurwitz rule again gives that $p_0=1$ and $p_m=c$. The \nit{cyclic label} $\clbl(x)$ of $x$ is then defined as the sequence $(c_1,\ldots,c_m)$ where $c_1=p_1$ and $c_i=(c_1\cdots c_{i-1})^{-1}p_i$ for $i\in \intv{2,m}$ \cite[Definition 11.9]{beskpi1}. By \cite[Lemma 11.10]{beskpi1}, $\clbl(x)$ is always a length-additive decomposition of $c$. 

Now, if the above condition on the cyclic argument is not satisfied, or if $x\in \Hh$, we will define the cyclic label of $x$ as that of a desingularization of $x$. In order to show that such a definition is valid, we need to show that it doesn't depend on the desingularization we choose. 

Let $x\in V/W$, with cyclic support $(\kappa_1,\ldots,\kappa_m)$. For a family of paths $\gamma_1,\ldots,\gamma_m:[0,1]\to \C$, we define a path $\gamma:[0,1]\to E_n$ by setting $\gamma(t)$ as the multiset $\muls{\gamma_1(t),\ldots,\gamma_m(t)}$, where the multiplicity of $\gamma_j(t)$ is the same as that of $\kappa_j$ in $\LLL(x)$. We consider the following conditions on the family $(\gamma_1,\ldots,\gamma_m)$.

\begin{enumerate}[\quad $\bullet$]
\item For all $j\in \intv{1,m}$, we have $\gamma_j(t)=\kappa_j$ (in particular, $\gamma(1)=\LLL(x)$).
\item For all $t\in [0,1]$, we have $\gamma_i(t)\neq \gamma_j(t)$ for $i\neq j$.
\item For $t<1$, we have $\gamma(t)\in E_n^\circ$.
\item For $t\in [0,1]$, the cyclic support of $\gamma(t)$ is $(\gamma_1(t),\ldots,\gamma_m(t))$ (that is, the ordering of the cyclic support is preserved along $\gamma$), and the cyclic argument $(\theta_1(t),\ldots,\theta_m(t))$ of $\gamma(t)$ respects $0<\theta_1(t)<\ldots<\theta_m(t)<2\pi$.
\end{enumerate}
The first two conditions ensure in particular that the path $\gamma$ has a unique lift in $V/W$ terminating at $x$ (by Proposition \ref{prop:lift_lll}). If all these conditions are met, the lift of $\gamma$ in $V/W$ which terminates at $x$ is called a \nit{desingularization path} for $x$. One easily sees that every point of $V/W$ admits a desingularization path. For instance, if $\theta_m=2\pi$, then a path of the form $t\mapsto e^{i\epsi(1-t)}x$ for $\epsi>0$ small enough can give a desingularization path. If $\theta_j=\theta_{j+1}$, for some $j\in \intv{1,m}$ then we can replace $\kappa_i$ by $\gamma_j(t)=e^{i\epsi(1-t)}\kappa_j$ for $\epsi>0$ small enough. And, if $x\in \Hh$ (i.e. if $\kappa_m=0$), then we can replace $\kappa_m$ by a path of the form $\gamma_m(t)=e^{i\epsi(1-t)}i\epsi(1-t)$ for $\epsi>0$ small enough. Note that, if $x\in X/W$ respects the previous condition on the cyclic argument, then the constant path equal to $x$ is a desingularization path for $x$. 

Let $\gamma_1$ and $\gamma_2$ be two desingularization paths for $x$. We need to show that $\clbl(\gamma_1(t))$ and $\clbl(\gamma_2(t))$ are equal for $t$ big enough.

\begin{lem}\label{lem:desing_path_ca_marche}
Let $x\in V/W$, and let $\gamma:[0,1]\to V/W$ be a desingularization path for $x$. For all $t,t'<1$, we have $\clbl(\gamma(t))=\clbl(\gamma(t'))$. Furthermore, if $\gamma_1,\gamma_2$ are two desingularization paths for $x$, then we have $\clbl(\gamma_1(t))=\clbl(\gamma_2(t))$ for all $t<1$.
\end{lem}
\begin{proof}
Let $\gamma:[0,1]\to V/W$ be a desingularization path for $x$, and let $(\theta_1(t),\ldots,\theta_m(t))$ be the cyclic argument of $\gamma(t)$. By definition of a desingularization path, $\theta_j:[0,1]\to ]0,2\pi]$ is a continuous path for $j\in \intv{1,m}$, and we have $0<\theta_1(t)<\ldots<\theta_m(t)<2\pi$. For $j\in \intv{1,m}$ and $t<1$, the circular tunnel $(\gamma(t),\frac{\theta_j(t)+\theta_{j+1}(t)}{2})$ always represent the same element of $B(W)$ by the Hurwitz rule (with the convention that $\theta_0=0$ and $\theta_{m+1}=2\pi$). Since this element is the product of the first $j$ terms of $\clbl(\gamma(t))$, we have the first result.

Now, let $\gamma_1,\gamma_2$ be two desingularization paths for $x$.  For $t<1$, let us write $(\theta_1^i(t),\ldots,\theta_m^i(t))$ the cyclic argument of $\gamma_i(t)$ for $i=1,2$. By assumption, we have $0<\theta_1^i(t)<\ldots<\theta_m^i(t)<2\pi$ for $i=1,2$. For $j\in \intv{1,m}$, and $s\in [0,1]$, we set $\omega_j(s):=s(\theta_j^1(t))+(1-s)(\theta_j^2(t))$. By construction, we have $0<\omega_1(s)<\ldots<\omega_m(s)<2\pi$ for all $s\in [0,1]$. We obtain a path from $\gamma_1(t)$ to $\gamma_2(t)$. Applying the Hurwitz rule to this path gives that $\gamma_1(t)$ and $\gamma_2(t)$ share the same cyclic label. 
\end{proof}

Thanks to this lemma, the following definition is valid and doesn't depend on the choice of a desingularization path.

\begin{definition}\label{def:clbl}
Let $x\in V/W$, and let $\gamma:[0,1]\to V/W$ be a desingularization path for $x$. The \nit{cyclic label} $\clbl(x)$ of $x$ is defined as $\clbl(\gamma(t))$ for $t<1$ big enough. If $x\in \Hh$, then the \nit{outer label} $\olbl(x)$ is defined as $\clbl(x)$ minus its last entry. If $x\in X/W$, then $\olbl(x)$ is defined to be equal to $\clbl(x)$.
\end{definition}

Since the cyclic label of an element of $X/W$ is a length-additive decomposition of $c$, it is easy to recover $\olbl$ from $\clbl$ and vice versa. Note that $\olbl(x)$ is empty if and only if $\LL(x)=\muls{0}$, that is if and only if $x=0$ by \cite[Lemma 5.6]{beskpi1}. We now show that, as we expected, the cyclic label completely describes the fiber of the Lyashko-Looijenga morphism.

\begin{prop}\label{prop:trivialization_lll}
Let $x,x'\in V/W$. If $\LLL(x)=\LLL(x')$ and $\clbl(x)=\clbl(x')$, then $x=x'$. Conversely, let $\kappa\in E_n$ with cyclic support $(\kappa_1,\ldots,\kappa_m)$ and let $s:=(s_1,\ldots,s_m)\in D_m(c)$. If the multiplicity of $\kappa_i$ in $\kappa$ is equal to $\ell_\Rr(s_i)$ for all $i\in \intv{1,m}$, then there is a unique $x\in V/W$ such that $\LLL(x)=\kappa$ and $\clbl(x)=s$.
\end{prop}
\begin{proof} Assume that $x,x'\in V/W$ are such that $\LLL(x)=\LLL(x')$ and $\clbl(x)=\clbl(x')$. Since $x\in \Hh$ if and only if $0\in \LLL(x)$, we have that if either $x$ or $x'$ lie in $X/W$, then so does the other. In this case, we have $x=x'$ by \cite[Proposition 11.13]{beskpi1}. We can thus restrict our attention to the case where $x,x'\in \Hh$. Let $\gamma:[0,1]\to V/W$ be a desingularization path for $x$. Since $\LLL(x)=\LLL(x')$, and since the condition of being a desingularization path only depends on $\LLL\circ \gamma$, the unique lift $\gamma'$ of $\LLL\circ \gamma$ in $V/W$ which terminates at $x'$ is a desingularization path for $x'$. By definition, we have $\clbl(\gamma(t))=\clbl(x)=\clbl(x')=\clbl(\gamma'(t))$ for $t$ big enough. Since $\LLL(\gamma(t))=\LLL(\gamma'(t))$ by assumption, the first part of the proof gives that $\gamma(t)=\gamma'(t)$ for $t<1$ big enough. The respective endpoints of $\gamma$ and $\gamma'$, that is, $x$ and $x'$, are then equal.

Let now $(\kappa,s)\in E_n\times D_\bullet(c)$ be such that the cyclic support of $\kappa$ has the same length as $s$, and that the multiplicity of $\kappa_i$ in $\kappa$ is equal to $\ell_\Rr(s_i)$ for $i\in \intv{1,m}$. If the cyclic argument of $\kappa$ is strictly increasing and in $]0,2\pi[$, then in particular, $\kappa\in E_n^\circ$ and \cite[Proposition 11.13]{beskpi1} gives that there is a unique $x\in X/W$ such that $\LLL(x)=\kappa$ and $\clbl(x)=s$. Otherwise, let $\gamma_1,\ldots,\gamma_m:[0,1]\to \C$ be a family of paths which respects the conditions of a desingularization path for $(\kappa_1,\ldots,\kappa_m)$. By \cite[Proposition 11.13]{beskpi1}, there is a unique $x\in X/W$ such that $\LLL(x)=\gamma(0)$ and $\clbl(x)=s$. By construction, the unique lift $\ttilde{\gamma}$ of $\gamma$ inside $V/W$ starting at $x$ is a desingularization path for $\gamma(1)$. By Lemma \ref{lem:desing_path_ca_marche}, we have $\clbl(\gamma(t))=s$ for all $t<1$. By definition, $\gamma(1)$ is then a point of $V/W$ such that $\LLL(\gamma(1))=\kappa$ and $\clbl(\gamma(1))=s$.
\end{proof}

We can reformulate the first statement in the above proposition using outer labels instead of cyclic labels as follows.

\begin{cor}\label{cor:lll_olbl_caractérisent_égalité}
Let $x,x'\in V/W$. If $\LLL(x)=\LLL(x')$ and $\olbl(x)=\olbl(x')$, then $x=x'$.
\end{cor}
\begin{proof}
If $x\in X/W$, then $\LLL(x)=\LLL(x')$ implies that $x'\in X/W$. In this case, we have $\clbl(x)=\olbl(x)=\olbl(x')=\clbl(x')$ and $x=x'$ by Proposition \ref{prop:trivialization_lll}. If $x'\in \Hh$, then $\LLL(x)=\LLL(x')$ implies that $x'\in \Hh$. Let $\olbl(x)=(c_1,\ldots,c_m)$. Since the cyclic label must be a length-additive decomposition of $c$, we have
\[\clbl(x)=(c_1,\ldots,c_m,(c_1\cdots c_m)^{-1}c)=\clbl(x'),\] thus $x=x'$ by Proposition \ref{prop:trivialization_lll}.
\end{proof}

We now investigate how some (lifts of) paths in $E_n$ impact the cyclic label of an element of $V/W$. Let $x\in X/W$, with cyclic support $(\kappa_1,\ldots,\kappa_m)$. Suppose that we swap two consecutive points $\kappa_i$ and $\kappa_{i+1}$ of $\LLL(x)$, there are two natural ways to do so, one going ``farther'' than the other:

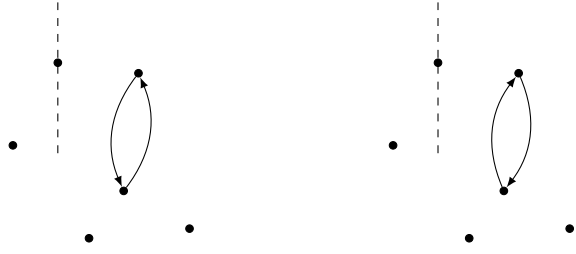
\begin{center}
\begin{tikzpicture}
\draw [dashed] (0,0)--(0,2);

\node[draw,circle,inner sep=1pt,fill] (P1) at (1.5*0.7071,1.5*0.7071) {};
\node[draw,circle,inner sep=1pt,fill] (P12) at (-30:1) {};
\node[draw,circle,inner sep=1pt,fill] (P2) at (-30:2) {};
\node[draw,circle,inner sep=1pt,fill] (P2) at (-70:1.2) {};
\node[draw,circle,inner sep=1pt,fill] (P2) at (-190:0.6) {};
\node[draw,circle,inner sep=1pt,fill] (P2) at (90:1.2) {};

\draw[->,>=latex] (P1) to[bend right] (P12);
\draw[->,>=latex] (P12) to[bend right] (P1);

\begin{scope}[shift=(0:5)]
\draw [dashed] (0,0)--(0,2);

\node[draw,circle,inner sep=1pt,fill] (P1) at (1.5*0.7071,1.5*0.7071) {};
\node[draw,circle,inner sep=1pt,fill] (P12) at (-30:1) {};
\node[draw,circle,inner sep=1pt,fill] (P2) at (-30:2) {};
\node[draw,circle,inner sep=1pt,fill] (P2) at (-70:1.2) {};
\node[draw,circle,inner sep=1pt,fill] (P2) at (-190:0.6) {};
\node[draw,circle,inner sep=1pt,fill] (P2) at (90:1.2) {};

\draw[->,>=latex] (P1) to[bend left] (P12);
\draw[->,>=latex] (P12) to[bend left] (P1);
\end{scope}
\end{tikzpicture}
\end{center}

By Proposition \ref{prop:lift_lll}, both of these paths in $E_n$ have a unique lift in $V/W$, that we denote by  $\gamma_i^{+}$ and $\gamma_i^-$, respectively.

\begin{lem}\label{lem:cyclic_hurwitz_moves}
Let $x\in X/W$, with cyclic support $(\kappa_1,\ldots,\kappa_m)$ and cyclic label $(c_1,\ldots,c_m)$. The cyclic labels of $\gamma_i^+(1)$ and $\gamma_i^{-}(1)$ are given by
\[\clbl(\gamma_i^+(1))=(c_1,\ldots,c_{i-1},c_{i+1},c_{i}^{c_{i+1}},c_{i+2},\ldots,c_m),\]
\[\clbl(\gamma_i^-(1))=(c_1,\ldots,c_{i-1},{}^{c_i}c_{i+1},c_{i},c_{i+2},\ldots,c_m).\]
\end{lem}
\begin{proof}
The concatenation of $\gamma_i^-$ and $\gamma_i^+$ gives a homotopically trivial path from $x$ to itself. In particular the assertion about $\clbl(\gamma_i^-(1))$ follows from that on $\clbl(\gamma_i^{+}(1))$. By construction of the cyclic label, we can assume that the cyclic argument $(\theta_1,\ldots,\theta_m)$ of $x$ satisfies $0<\theta_1<\ldots<\theta_m<2\pi$. Let then $\theta$ be such that $\frac{\theta_{i-1}}{h}<\theta<\frac{\theta_{i}}{h}$. By \cite[Lemma 11.1]{beskpi1}, one can replace $x$ by $e^{i\theta}x$ and consider $x':=\gamma_1^+(x)$.

Now if $\theta'$ is such that $\frac{\theta_2}{h}<\theta'<\frac{\theta_3}{h}$, then the circular tunnels $(x,\theta')$ and $(x',\theta')$ represent the same element in $B(W)$ by the Hurwitz rule. That is the product of the first two terms of $\clbl(x)$ and $\clbl(x')$ are equal.  The Hurwitz rule also shows that the terms of $\clbl(x)$ and $\clbl(x')$ are equal for $i>2$.

Lastly, the assertion that the first two terms of $\clbl(x')$ are $(c_2,c_1^{c_2})$ follows from \cite[Lemma 11.11]{beskpi1} and \cite[Corollary 6.20]{beskpi1}.
\end{proof}

By induction, we see that in general, moving a point of the cyclic support of some $x\in X/W$ may only change the terms of the cyclic label corresponding to points of the cyclic support with lower modulus.

Another result that will be useful to us is the behavior of the cyclic label when ones bring points of the cyclic support affinely towards the center.

\begin{lem}\label{lem:shrinkinglemma}
Let $x\in X/W$ with cyclic support $(\kappa_1,\ldots,\kappa_m)$ and cyclic label $(c_1,\ldots,c_m)$. We assume that all the $\kappa_i$ have distinct arguments. Let $I\sqcup J=\intv{1,m}$ be a partition of $\intv{1,m}$, with $J=\{j_1<\ldots<j_p\}$.
Let $\gamma_I:[0,1]\to E_n$ be the path starting at $\LLL(x)$ and such that the cyclic support of $\gamma(t)$ is made of the points $\kappa_j$ for $j\in J$ and the points $(1-t)\kappa_i$ for $i\in I$. There is a unique lift $\ttilde{\gamma}_I$ of $\gamma_I$ in $V/W$, and we have
\[\clbl(\ttilde{\gamma}_I(1))=(c_{j_1},\ldots,c_{j_p},(c_{j_1}\cdots c_{j_p})^{-1}c).\]
\end{lem}
\begin{proof}
First, up to a rotation, we can assume that the last term of the cyclic argument of $x$ is strictly smaller than $2\pi$. The existence and uniqueness of $\ttilde{\gamma}_I$ comes from Proposition \ref{prop:lift_lll}. Let $\gamma_1$ be the path starting from $\LLL(x)$ consisting in sliding the points $x_i$ with $i\in I$ together clockwise (and going closer than $0$ than all the $x_j$ with $j\in J$) and next to the vertical half-line $i\R_{\geqslant 0}$ (closer than $\kappa_{j_{\max}}$, with $j_{\max}:=\max J$). Let then $\gamma_2$ be the path starting from $\gamma_1(1)$ consisting in sliding this last point down towards $0$.

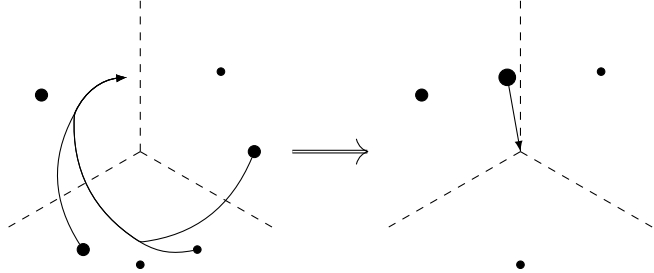
\begin{center}
\begin{tikzpicture}[scale=1]
\draw [dashed] (0,0)--(0,2);
\draw [dashed] (0,0)--(+2*0.866,-1);
\draw [dashed] (0,0)--(-2*0.866,-1);

\node[draw,circle,inner sep=1pt,fill] (P1) at (1.5*0.7071,1.5*0.7071) {};
\node[draw,circle,inner sep=1.6pt,fill] (P2) at (1.5,0) {};

\node[draw,circle,inner sep=1pt,fill] (P3) at (1.5*0.5,-1.5*0.8660) {};
\node[draw,circle,inner sep=1pt,fill] (P4) at (0,-1.5) {};
\node[draw,circle,inner sep=1.6pt,fill] (P5) at (-1.5*0.5,-1.5*0.8660) {};

\node[draw,circle,inner sep=1.6pt,fill] (P6) at (-1.5*0.8660,1.5*0.5) {};

\draw[->,>=latex,rounded corners=1pt] (P2) to[bend left] (0,-1.2) to[bend left] (-0.8660,0.5) to[bend left] (100:1);
\draw[->,>=latex,rounded corners=1pt] (P3) to[bend left=20] (0,-1.2) to[bend left] (-0.8660,0.5) to[bend left] (100:1);
\draw[->,>=latex,rounded corners=1pt] (P5) to[bend left] (-0.8660,0.5) to[bend left] (100:1);

\draw [double,->] (2,0) to (3,0);

\begin{scope}[shift=(0:5)]
\draw [dashed] (0,0)--(0,2);
\draw [dashed] (0,0)--(2*0.866,-1);
\draw [dashed] (0,0)--(-2*0.866,-1);

\node[draw,circle,inner sep=1pt,fill] (P1) at (1.5*0.7071,1.5*0.7071) {};

\node[draw,circle,inner sep=1pt,fill] (P4) at (0,-1.5) {};

\node[draw,circle,inner sep=1.6pt,fill] (P6) at (-1.5*0.8660,1.5*0.5) {};
\node[draw,circle,inner sep=2.2pt,fill] (P7) at (100:1) {};

\draw[->,>=latex] (P7) to (0,0);
\end{scope}
\end{tikzpicture}
\end{center}

The path $\ttilde{\gamma}_I$ is homotopic to the unique lift in $V/W$ starting at $x$ of the concatenation $\gamma_1*\gamma_2$. Since moving points in the cyclic support of some $x\in V/W$ may only affect the terms of the cyclic label corresponding to points of the cyclic support with lower modulus by Lemma \ref{lem:cyclic_hurwitz_moves}, and since the cyclic label is a length-additive decomposition of $c$, we have 
\[\clbl(\gamma_1(1))=(c_{j_1},\ldots,c_{j_p},(c_{j_1}\cdots c_{j_p})^{-1}c).\]
The path $\gamma_2$ is a desingularization path for $\gamma_2(1)$. For all $t<1$, we have by the Hurwitz rule that $\clbl(\gamma_2(t))=\clbl(\gamma_2(0))=\clbl(\gamma_1(1))=(c_{j_1},\ldots,c_{j_p},(c_{j_1}\cdots c_{j_p})^{-1}c).$
By definition of $\clbl(\gamma_2(1))$, we deduce that $\clbl(\gamma_2(1))=(c_{j_1},\ldots,c_{j_p},(c_{j_1}\cdots c_{j_p})^{-1}c)$ as claimed.
\end{proof}

Now, Proposition \ref{prop:trivialization_lll} gives a complete description of the topological pair $(X/W,V/W)$ via the map $(\LLL,\clbl)$. We want to restrict this map to $((X/W)^{\mu_d},(V/W)^{\mu_d})$ in order to get a convenient description of this pair, in the spirit of \cite[Lemma 11.14]{beskpi1}. First, as $E_n$ is made of multisets inside $\C$, it comes equipped with an action of $\C^*$. In particular we can define $(E_n)^{\mu_k}$ as the fixed points of $E_n$ under the action of the $k$-th roots of unity. Since the map $\LLL$ is homogeneous of degree $h$ by \cite[Lemma 11.1]{beskpi1}, we have $\LLL(\zeta_d x)=(\zeta_d)^h\LLL(x)$. As $(\zeta_d)^h$ is a primitive $p$-th root of unity, we obtain that $\LLL(x)\in (E_n)^{\mu_p}$ if $x\in (V/W)^{\mu_d}$. Thus we are interested in describing points $x\in V/W$ such that $\LLL(x)\in (E_n)^{\mu_p}$.

\begin{lem}\label{lem:bes11.14}
Let $x\in V/W$ be nonzero. Assume that $\LLL(x)\in (E_n)^{\mu_p}$. Then $\olbl(x)$ contains $pk$ terms for some positive integer $k$, and $\olbl(\zeta_{ph}x)=\tau^k(\olbl(x))$.
\end{lem}
\begin{proof}
If $x\in X/W$, then this is simply \cite[Lemma 11.14]{beskpi1}. We restrict our attention to the case where $x\in \Hh$. Since $\zeta_p\LLL(x)=\LLL(x)$ by assumption, the cyclic support of $\LLL(x)$ is a reunion of $\mu_p$-orbits in $\C$. As each nonzero $\mu_p$-orbit in $\C$ contains exactly $p$ points, we obtain that $\olbl(x)$ contains $pk$ terms, where $k$ is the number of $\mu_p$-orbits in the cyclic support of $\LLL(x)\setminus \{0\}$ (since $x\neq 0$, we have $k>0$ by \cite[Lemma 5.6]{beskpi1}). We write $\olbl(x)=(c_1,\ldots,c_{pk})$ and $c_{pk+1}:=(c_1\cdots c_{pk})^{-1}c$ so that $\clbl(x)=(c_1,\ldots,c_{pk},c_{pk+1})$. 

Let $\gamma:[0,1]\to V/W$ be a desingularization path for $x$. We denote by $\gamma_j:[0,1]\to \C$ the $j$-th component of the cyclic support of $\LLL\circ \gamma$, and let $(\theta_1(t),\ldots,\theta_{pk}(t),\theta_{pk+1}(t))$ be the cyclic argument of $\gamma(t)$. By \cite[Lemma 11.11]{beskpi1}, the cyclic label of $\zeta_{ph}\gamma(t)$ is given for $t$ big enough by $(c_{k+1},c_{k+2},\ldots,c_{pk},c_0,c_1^c,\ldots,c_k^c)=\tau^k(\clbl(\zeta_{ph}\gamma(t)))$
The cyclic argument of $\zeta_{ph}\gamma(t)$ is given by
\[\left(\theta_{k+1}(t)-\frac{2\pi}{p},\ldots,\theta_{pk+1}(t)-\frac{2\pi}{p},\theta_1(t)-\frac{2\pi}{p}+2\pi,\ldots,\theta_{k}(t)-\frac{2\pi}{p}+2\pi\right).\]
In particular, $\zeta_{ph}\gamma$ is not a desingularization path for $\zeta_{ph}x$, since it doesn't respect the condition on (the path $\zeta_{ph}\gamma_{pk+1}$ which terminates at $0\in \LLL(x)$ does not give the last term of the cyclic support of $\zeta_{ph}\gamma$). We fix $r>0$ such that, for $t<1$ big enough, we have $|\gamma_{pk+1}(t)|<r<|\gamma_{j}(t)|$ for all $j\neq pk+1$. We can then rotate $\zeta_{ph}\gamma_{pk+1}(t)$ by an angle of $K(t)$ so that $\theta_k(t)-\frac{2\pi}{p}+2\pi<\theta_{pk+1}(t)-\frac{2\pi}{p}+K(t)<2\pi$. We obtain a new path $\gamma_{pk+1}'(t)$. The family $(\zeta_{ph}\gamma_1(t),\ldots,\zeta_{ph}\gamma_{pk}(t),\gamma_{pk+1}'(t))$ induces a desingularization path $\gamma'$ for $\zeta_{ph}x$, which is homotopic to $\zeta_{ph}\gamma$.

However, by Lemma \ref{lem:cyclic_hurwitz_moves}, the cyclic label of $\gamma'(t)$ is given (by definition), for $t$ big enough, by 
\[\clbl(\gamma'(t))=(c_{k+1},c_{k+2},\ldots,c_{pk},c_1^c,\ldots,c_k^c,(c_0)^{c^{-1}c_1\cdots c_k c}).\]
Thus, $\olbl(x)=(c_{k+1},c_{k+2},\ldots,c_{pk},c_1^c,\ldots,c_k^c)$ as claimed.
\end{proof}

Using this Lemma, we get a complete description of $(V/W)^{\mu_d}$ using the map $(\LLL,\clbl)$.

\begin{lem}\label{lem:caractérisation_v/wmu_d} For all nonzero $x\in V/W$, the following assertions are equivalent
\begin{enumerate}
\item $x\in (V/W)^{\mu_d}$.
\item $\LLL(x) \in (E_n)^{\mu_p}$ and $\tau^{kq}(\olbl(x))=\olbl(x)$, where $pk$ is the cardinality of the cyclic support of $\LLL(x)\setminus \{0\}$.
\end{enumerate}
\end{lem}
\begin{proof}
Assume $(i)$. We have $\LLL(x)=(\zeta_d)^h\LLL(x)$, thus $\LLL(x)\in (E_n)^{\mu_p}$ since $(\zeta_d)^h$ is a primitive $p$-th root of unity. By Lemma \ref{lem:bes11.14}, we have $\olbl(\zeta_d x)=\olbl(\zeta_{d'h}^qx)=\tau^{kq}(\olbl(x))$, where $k$ is the number of points of the cyclic support of $\LLL(x)$ with argument $\theta\in ]0,2\pi/p]$.

Conversely, assuming $(ii)$, we conclude that $x$ and $\zeta_dx$ satisfy $\LLL(x)=\LLL(\zeta_d x)$ and $\olbl(x)=\tau^{kq}(\olbl(x))=\olbl(\zeta_d x)$ by Lemma \ref{lem:bes11.14}. We then have $x=\zeta_dx$ by Corollary \ref{cor:lll_olbl_caractérisent_égalité}.
\end{proof}

We finish this section by giving the topological definition of the Springer groupoid. We consider the following union of half-lines
\[D:=\bigcup_{\zeta\in \mu_p} \zeta i\R_{\geqslant 0}\]
and $\Uu^{\mu_d}=(X/W)^{\mu_d}\cap \Uu$. The complement of $D$ in $\C$ consists of $p$ sectors $P_1,\ldots,P_p$, which are labeled so that $P_k$ contains all nonzero points of arguments $\theta\in [\frac{\pi}{2}-k\frac{2\pi}{p},\frac{\pi}{2}-(k-1)\frac{2\pi}{p}[$, as in the example below (where $p=5$)
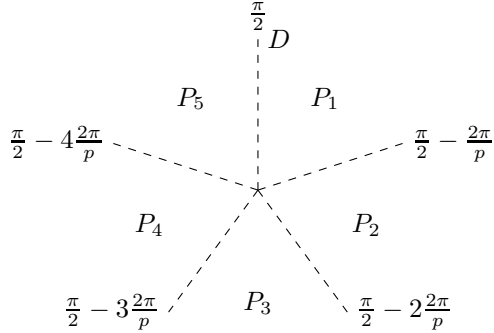
\begin{center}\begin{tikzpicture}[scale=1]

\draw [dashed] (0,0)--(90:2);
\node[above] at (90:2) {$\frac{\pi}{2}$};
\draw [dashed] (0,0)--(90+72:2);
\node[left] at (90+72:2) {$\frac{\pi}{2}-4\frac{2\pi}{p}$};
\draw [dashed] (0,0)--(90+2*72:2);
\node[left] at (90+2*72:2) {$\frac{\pi}{2}-3\frac{2\pi}{p}$};
\draw [dashed] (0,0)--(90+3*72:2);
\node[right] at (90+3*72:2) {$\frac{\pi}{2}-2\frac{2\pi}{p}$};
\draw [dashed] (0,0)--(90+4*72:2);
\node[right] at (90+4*72:2) {$\frac{\pi}{2}-\frac{2\pi}{p}$};

\node[right] at (0,2) {$D$};

\node at (90+36+0*72:1.5) {$P_5$};
\node at (90+36+1*72:1.5) {$P_4$};
\node at (90+36+2*72:1.5) {$P_3$};
\node at (90+36+3*72:1.5) {$P_2$};
\node at (90+36+4*72:1.5) {$P_1$};
\end{tikzpicture}\end{center}

Since $\LLL(x)\in (E_n)^{\mu_p}$ when $x\in (V/W)^{\mu_d}$, we have
\[\Uu^{\mu_d}=\{x\in (X/W)^{\mu_d}~|~\LLL(x)\cap D=\varnothing\}.\]
Let $x\in (V/W)^{\mu_d}$. The points of the cyclic support of $\LLL(x)\setminus \{0\}$ can be partitioned into $p$ groups, according to which sector $P_i$ they lie in for $i\in \intv{1,p}$. 

\begin{definition}
Let $x\in V/W$. The \nit{cyclic support} $\cc(x)$ of $x$ is the sequence $(c_1,\ldots,c_p)$, where $c_i$ is the product (in the ordering given by the cyclic support of $x$) of the terms of $\clbl(x)$ associated to points in the sector $P_i$, for $i\in \intv{1,p}$.
\end{definition}
This is a generalization of \cite[Definition 11.18]{beskpi1}, which covers the case where $x\in (X/W)^{\mu_d}$. Note that $\cc(x)=(1,\ldots,1)$ if and only if $x=0$, again by \cite[Lemma 5.6]{beskpi1}. A first consequence of Lemma \ref{lem:caractérisation_v/wmu_d} is that $\tau^q\cc(x)=\cc(x)$ if $x\in (V/W)^{\mu_d}$. In this case, since $\cc(x)$ has length $p$, \cite[Lemma 1.28]{springercat} gives that $\cc(x)$ has the form $(c_1,c_1^{c^{\eta}},\ldots,c_1^{c^{(p-1)\eta}})$ and is entirely determined by its first term. Moreover, if $x\in \Uu^{\mu_d}$, then $\cc(x)$ is a length-additive decomposition of $c$, and $\cc(x)\in D_p^q(c)$. 

By \cite[Lemma 11.22]{beskpi1}, the cyclic content induces a bijection between the connected components of $\Uu^{\mu_d}$ and $D_p^q(c)$. The connected components of $\Uu^{\mu_d}$ are furthermore open and contractible \cite[Lemma 11.22]{beskpi1}. We can define the \nit{Springer groupoid} (associated to $W$ and $d$) as the groupoid $\Bb:=\pi_1((X/W)^{\mu_d},\Uu^{\mu_d})$ made of classes of paths from some point of $\Uu^{\mu_d}$ to some other one, up to a homotopy leaving the endpoints inside $\Uu^{\mu_d}$. Let $u\in \pi_0(\Uu^{\mu_d})\simeq D_p^q(c)$. By definition of $\Bb$, the group $\Bb(u,u)$ is canonically identified with every group $\pi_1((X/W)^{\mu_d},x)$ with $x\in u$.

\subsection{Parabolic subgroups and local braid groupoids}
Let $W<\GL(V)$ be an irreducible complex reflection group (not necessarily well-generated), and let $x_0\in V/W$ lie on the discriminant hypersurface $\Hh$. Let also $*\in X/W$ be a basepoint so that $B(W)=\pi_1(X/W,*)$. We start by recalling the main definitions of \cite[Section 2]{paratresses} about parabolic subgroups of complex braid groups.

First consider a \nit{normal ray} \cite[Section 2.1]{paratresses} for the topological pair $(X/W,V/W)$, based at $*$ and terminating at $x_0$. By definition of a normal ray, the set $\eta(]\alpha,1[)\subset X/W$ is always simply connected when $\alpha$ is big enough. The fundamental groups $\pi_1(X/W,\eta(t))$ for $t>\alpha$ are then canonically identified with one another. We denote $\pi_1(X/W,\eta):=\pi_1(X/W,\eta(t))$ for $t$ big enough. In particular, we have $\pi_1(X/W,\eta)\simeq \pi_1(X/W,\eta(0))=B(W)$. By \cite[Proposition 2.1]{paratresses}, the image of $\pi_1(X/W\cap U,\eta)$ in $\pi_1(X/W,\eta)$ for $U\subset V/W$ a neighborhood of $x_0$ does not depend on $U$ for $U$ small enough. This image is denoted by $\pi_1^{\loc}(X/W,\eta)$ and is called a \nit{local fundamental group} (of 
$X/W$ at $\eta$). The image of $B(W)$ of such a local fundamental group is, by definition, a \nit{parabolic subgroup} of $B(W)$.

Parabolic subgroups of $B(W)$ are stable under conjugacy, and a change of basepoint $B(W)\simeq \pi_1(X/W,x)$ induces a bijection between the associated collections of parabolic subgroups \cite[Proposition 2.5]{paratresses}.

If $B_0\subset B(W)$ is a parabolic subgroup, then its image $W_0$ in $W$ is a parabolic subgroup of $W$ by \cite[Proposition 2.5]{paratresses}, and we have $B_0\simeq B(W_0)$. Furthermore, by \cite[Proposition 2.6]{paratresses}, two parabolic subgroups $B_1,B_2\subset B(W)$ are conjugate if and only if their images in $W$ are conjugate. In particular, if $\eta$ is a normal ray based at $*$ and terminating at some $x_0$, then the induced parabolic subgroup $B_0$ depends -up to conjugacy- only on the stratum of the discriminant stratification to which $x_0$ belongs. 

We assume again from now on that $W$ is well-generated, and we use the notation from the last Section. Let $g\in W$ be a $\zeta_d$-regular element. By Theorem \ref{theo:regular_elements}, we have a homeomorphism of topological pairs $(X_g/W_g,V_g/W_g)\simeq ((X/W)^{\mu_d},(V/W)^{\mu_d})$. Since the definitions of normal ray, local fundamental group and parabolic subgroup are purely topological, this homeomorphism induces a bijection between the parabolic subgroups of $B(W_g)$ and the parabolic subgroups defined for the topological pair $((X/W)^{\mu_d},(V/W)^{\mu_d})$. 

In \cite{paratresses}, the parabolic subgroups of $W$ are understood in terms of a Garside group structure on $B(W)$. In the case of $(X/W)^{\mu_d}$, the Garside structure will only exist at the level of the Springer groupoid. Thus we wish to adapt the notion of local fundamental group to a notion of local fundamental groupoid, more suitable for understanding the relation with the Springer groupoid. We begin by defining a particular neighborhood basis of elements of $(V/W)^{\mu_d}$.

Let $\kappa:=\muls{z_1,\ldots,z_n}\in E_n$. By definition of the topology of $E_n$, a neighborhood of $\kappa$ in $E_n$ is the image of a neighborhood of $(z_1,\ldots,z_n)$ in $\C^n$ under the projection map $\C^n \to E_n$. In particular, choosing a neighborhood $U_{z_i}$ in $\C$ of each of the $z_i$ induces a neighborhood $U$ of $\kappa$ in $E_n$ (the image of $\prod_{i=1}^n U_{z_i}$ under the projection map $\C^n\to E_n$). For $z\in \C$ and $r>0$, we denote by $\D(z,r)$ the open disk with center $z$ and radius $r$. 

\begin{definition} Let $\kappa=\muls{z_1,\ldots,z_n}\in E_n$ be such that $\kappa\cap D\subset \{0\}$. A \nit{confining neighborhood} of $\kappa$ in $E_n$ is the image in $E_n$ of a set of the form $\prod_{i=1}^n \D(z_i,r_i)\subset \C^n$ such that
\begin{enumerate}
\item $r_i=r_j$ if $z_i=z_j$.
\item $\D(z_i,r_i)\cap \D(z_j,r_j)=\varnothing$ if $z_i\neq z_j$.
\item For $z_i\in \kappa\setminus \{0\}$, we have $\D(z_i,r_i)\cap D=\varnothing$. 
\item There is some $r>0$ such that, for all $i\in \intv{1,n}$, we either have
\begin{itemize}
\item $z_i=0$, in which case $r_i<r$ and $\D(z_i,r_i)\subset \D(0,r)$,
\item $z_i\neq 0$, in which case $r+r_i<|z_i|$ and $\D(z_i,r_i)\cap \D(0,r)=\varnothing$.
\end{itemize} 
\end{enumerate}
If $U$ is a confining neighborhood of $\kappa$, the sets $\D(z_i,r_i)$ with $z_i\neq 0$ are called \nit{outer disks}, and the set $\D(z_i,r_i)$ with $z_i=0$ (which is unique by condition (1)) is called the \nit{central disk}. \\Let $z\in (V/W)^{\mu_d}$ be a point such that $\LLL(z)\cap D\subset\{0\}$. By extension, we say that a neighborhood of $z$ in $(V/W)^{\mu_d}$ is \nit{confining} if it is path connected and if its image under $\LLL$ is a confining neighborhood of $\LLL(z)$.
\end{definition}

Note that, by condition $(2)$, the notion of outer and central disk of a confining neighborhood depends only on the confining neighborhood. Condition $(2)$ is what justifies the name confining, the different disks isolates the points of the cyclic support of $\kappa$ from one another, so that no singularization can happen inside of a confining neighborhood.  Condition $(4)$ is a strengthening of condition $(2)$ in the case where one of the points $z_i$ is $0$. This condition will allow us to apply Lemma \ref{lem:cyclic_hurwitz_moves}. Indeed, since points in the central disk cannot move farther than points in outer disks, they cannot change their cyclic label by Lemma \ref{lem:cyclic_hurwitz_moves}. This will be important in the proof of Lemma \ref{lem:terms_of_outer_points_constant_in_confining_neighborhood}.

Let $\kappa\in E_n$ be such that $\kappa\cap D\subset \{0\}$. It is easily seen that confining neighborhoods of $\kappa$ form a neighborhood basis of $\kappa$ in $E_n$. Let $z\in (V/W)^{\mu_d}$ be such that $\LLL(z)\cap D\subset \{0\}$, and let $U$ be a confining neighborhood of $\LLL(z)$. By Proposition \ref{prop:lift_lll}, the path connected component $U'$ of $z$ in $\LLL^{-1}(U)$ is a confining neighborhood of $z$ in $(V/W)^{\mu_d}$ such that $\LLL(U')=U$. From this we see that confining neighborhoods of $z$ form a neighborhood basis of $z$ in $(V/W)^{\mu_d}$. Furthermore, if $U$ is a confining neighborhood of $z$, then $U\cap \Uu^{\mu_d}$ is nonempty since $\Uu^{\mu_d}$ is dense in $(V/W)^{\mu_d}$.

Let $z\in (V/W)^{\mu_d}$ be such that $\LLL(z)\cap D\subset \{0\}$, and let $U$ be a confining neighborhood of $z$ in $(V/W)^{\mu_d}$. For $x\in U$, we distinguish between the \nit{outer points} of $\LLL(x)$, which lie in an outer disk of $U$, and the \nit{central points} of $\LLL(x)$, which lie in the central disk of $U$.

We see that the definition of confining neighborhood does not make sense for arbitrary points of $(V/W)^{\mu_d}$. However, in our study of parabolic subgroups up to conjugacy, we can restrict our attention to a family of points for which confining neighborhoods are defined.

\begin{prop}\label{prop:path_preserves_multiplicity_lll}
Let $x\in V/W$, and let $\gamma:[0,1]\to V/W$ be a path starting at $x$. If the multiplicity of $\Hh$ at $\gamma(t)$ stays the same for all $t\in [0,1]$, then $\gamma(1)$ and $x$ lie on the same stratum of the discriminant stratification. 
\end{prop}
\begin{proof}
By \cite[Corollary 5.9]{beskpi1}, the multiplicity of the discriminant hypersurface $\Hh$ at $x$ is equal to the multiplicity of $0$ in $\LLL(x)$.

Let $p:V\to V/W$ be the canonical projection map, and let $\Aa$ be the set of reflecting hyperplanes of $W$. The discriminant hypersurface is the reunion of the $p(H)$ for $H\in \Aa$. For $x=p(v)\in V/W$, the stratum of the discriminant stratification to which $x$ belong is entirely determined by the sets $p(H)$ for $H\in \Aa$ to which $x$ belong. 

Now, for all strata $L$, the set $I_L=\{t\in [0,1]~|~\gamma(t)\in \bbar{L}\}$ is closed. If $L_0$ is the stratum containing $x$, we claim that $I_{L_0}=[0,1]$. If this is not the case, then we can consider $[0,t_0]$ the connected component of $0$ in $I_{L_0}$, with $t_0<1$.

Let $v\in V$ be such that $p(v)=\gamma(t_0)$. We can consider a neighborhood $\ttilde{U}$ of $v$ in $V$ which intersects only the reflecting hyperplanes which contain $v$. Since $p$ is a branched covering map \cite[Proposition 6.106]{orlikterao}, $U:=p(\ttilde{U})$ is a neighborhood of $\gamma(t_0)$ in $V/W$. By construction, $U$ intersects some $p(H)$ if and only if $\gamma(t_0)\in p(H)$. 

For $\epsi>0$ small enough, we have $\gamma(t_0+\epsi)\in U$. If $H$ is such that $\gamma(t_0+\epsi)\in p(H)$, then $p(H)\cap U\neq \varnothing$ and $\gamma(t_0)\in p(H)$. Thus we have
\[\{p(H)~|~\gamma(t_0+\epsi)\in p(H)\}\subsetneq \{p(H)~|~\gamma(t_0)\in p(H)\}.\] 
The multiplicity of $0\in \LLL(\gamma(t_0+\epsi))$ is then smaller than that of $0\in \LLL(\gamma(t_0))$, which contradicts the assumption.

We have that $\gamma(t)\in \bbar{L_0}$ for all $t\in [0,1]$. Since the only stratum contained in $\bbar{L_0}$ on which the multiplicity of the discriminant is equal to the multiplicity of $0$ in $\LLL(x)$ is $L_0$, we obtain the result.
\end{proof}

In particular, in the situation of the above proposition, two normal rays terminating at $x$ and $\gamma(1)$ define conjugate parabolic subgroups.

A strong consequence of the above proposition is that two points sharing the same cyclic content belong to the same stratum of the discriminant stratification. As we won't use this result for now, we postpone the proof until Proposition \ref{prop:replace_by_z_s}, which is a stronger statement. For now, we show that, in our study of parabolic subgroups up to conjugacy, we can restrict our attention to points $x$ such that $\LLL(x)$ admits one point in each sector $P_i$ for $i\in \intv{1,p}$.

\begin{cor}\label{cor:replace_by_standard_point}
Let $x\in (V/W)^{\mu_d}$ be nonzero. There is a point $x'\in (X/W)^{\mu_d}$ lying on the same stratum of the discriminant stratification as $x$, such that $\LLL(x')$ admits exactly one point in each sector $P_i$ for $i\in \intv{1,p}$, and such that $\LLL(x')\cap D\subset \{0\}$.
\end{cor}
\begin{proof}
Consider the path $\gamma$ in $E_n$, starting from $\LLL(x)$ and consisting in sliding affinely together nonzero points in each sector $P_i$ for $i\in \intv{1,p}$ (such points exist since $x$ is nonzero). By Proposition \ref{prop:lift_lll}, the path $\gamma$ has a unique lift $\ttilde{\gamma}$ inside $V/W$ starting at $x$. By Lemma \ref{lem:caractérisation_v/wmu_d}, the path $\ttilde{\gamma}$ lies inside $(V/W)^{\mu_d}$. By the Hurwitz rule, we have $\cc(\ttilde{\gamma}(t))=\cc(x)$ for all $t\in [0,1]$. By construction of $\gamma$, the multiplicity of $\Hh$ at $\gamma(t)$ is the same for all $t\in [0,1]$, Proposition \ref{prop:path_preserves_multiplicity_lll} then gives that $x$ and $\gamma(1)$ lie on the same stratum of the discriminant stratification, and $\gamma(1)$ contains one point in each sector $P_i$ for $i\in \intv{1,p}$ by construction.

Now, up to replacing $\gamma(1)$ by $e^{i\theta}\gamma(1)$ for $\theta>0$ small enough, we can assume that $\LLL(\gamma(1))\cap D\subset \{0\}$. Note that $\gamma(1)$ and $e^{i\theta}\gamma(1)$ for $\theta>0$ small enough share the same cyclic label and belong to the same stratum of the discriminant stratification.
\end{proof}
We say that a point $x\in (V/W)^{\mu_d}$ is \nit{standard} if either $x=0$, or if $x\neq 0$ satisfies the condition of the above lemma. If $x\in (V/W)^{\mu_d}$ is such that $\LL(x)\cap D\subset \{0\}$, and such that $\LLL(x)$ admits at most one point in each sector $P_i$ for $i\in \intv{1,p}$, then by Lemma \ref{lem:caractérisation_v/wmu_d}, we either have that $\LLL(x)$ admits exactly one point in each sector $P_i$ for $i\in \intv{1,p}$, in which case $\olbl(x)=\cc(x)$, or that $\LLL(x)=\{0\}$, in which case we have $x=0$ by \cite[Lemma 5.6]{beskpi1}. In both cases, we see that $x$ is standard. We will later restrict our attention to a particular set of standard points in Section \ref{sec:correspondence_parabolic_localbraid}.

\begin{lem}\label{lem:terms_of_outer_points_constant_in_confining_neighborhood}
Let $z\in (V/W)^{\mu_d}$ be a standard point, and let $U$ be a confining neighborhood of $z$ in $(V/W)^{\mu_d}$. For $x\in U$ and $i\in \intv{1,p}$, the product of the terms of $\clbl(x)$ corresponding to the outer points in the sector $P_i$ (in the ordering given by the cyclic support of $x$) is equal to the $i$-th term of $\cc(z)=\olbl(z)$.
\end{lem}
\begin{proof}
By definition, there is a path $\gamma$ from $x$ to $z$ in $U$. For $t\in [0,1]$, let $\beta(t)$ be the subsequence of $\clbl(\gamma(t))$ corresponding to outer points in the sector $P_1$. By definition of a confining neighborhood, the outer points of $\LLL(\gamma(t))$ always have strictly greater modulus than central points. By Lemma \ref{lem:cyclic_hurwitz_moves}, moving points in the cyclic support of some point in $(V/W)^{\mu_d}$ may only affect the terms of the cyclic label corresponding to points of the cyclic support with lower modulus. Thus, $\beta(t)$ can be obtained from $\beta(0)$ by a sequence of elementary moves including
\begin{enumerate}[\quad $\bullet$]
\item Replacing a term $\beta_i$ by a subsequence $\beta_{i,1},\ldots,\beta_{i,m}$ such that $\beta_{i,1}\cdots \beta_{i,m}=\beta_i$.
\item Replacing a subsequence $\beta_{i},\ldots,\beta_{i+m}$ by its product $\beta_{i}\cdots \beta_{i+m}$.
\item Replacing a length 2 subsequence $\beta_i,\beta_{i+1}$ by either $\beta_{i+1}, \beta_i^{\beta_{i+1}}$ or ${}^{\beta_i} \beta_{i+1}, \beta_{i}$.
\end{enumerate}
Since all these elementary moves preserve the product of the sequence, we obtain that the product of $\beta(t)$ is equal to the product of $\beta(0)$. Since this is true for all $t\in [0,1]$, we obtain that the product of $\beta(0)$ is equal to the product of $\beta(1)$, i.e. to the first term of $\cc(z)=\olbl(z)$. We apply the same reasoning in all the sectors $P_i$ to obtain the result.
\end{proof}

Let $U\subset (V/W)^{\mu_d}$ be a confining neighborhood of some standard point $z\in (V/W)^{\mu_d}$. We wish to consider $U\cap \Uu^{\mu_d}$ as a basepoint with several path connected components for $U\cap (X/W)^{\mu_d}$. To do this, we need to show that the connected components of $U\cap \Uu^{\mu_d}$ are simply connected. We actually show an analogue of \cite[Lemma 11.22]{beskpi1} in this context. That is, path connected components of $U\cap \Uu^{\mu_d}$ are in bijection with the possible values of the cyclic content, and they are contractible.

\begin{prop}\label{prop:conn_compo=cyc_cont}
Let $U$ be a confining neighborhood of a standard point in $(V/W)^{\mu_d}$. Two points of $U\cap \Uu^{\mu_d}$ are in the same path-connected component if and only if they share the same cyclic content. Moreover, the path connected components of $U\cap \Uu^{\mu_d}$ are contractible.
\end{prop}
\begin{proof}
If $x,y$ lie in the same path connected component of $U\cap \Uu^{\mu_d}$, then, in particular, $x$ and $y$ lie in the same connected component of $\Uu^{\mu_d}$, thus $\cc(x)=\cc(y)$ by \cite[Lemma 11.22]{beskpi1}.

Now, let $x,y\in U\cap \Uu^{\mu_d}$ share the same cyclic content. First, we consider the paths in $E_n$ which consist in 
\begin{enumerate}
\item In each sector $P_i$ for $i\in \intv{1,p}$, slide affinely together the points in the central disk, and rotate the resulting point counterclockwise next to the associated half-line. 
\item Slide affinely together the points in each outer disks towards its center.
\end{enumerate} 
As in the example below (with $p=3$)

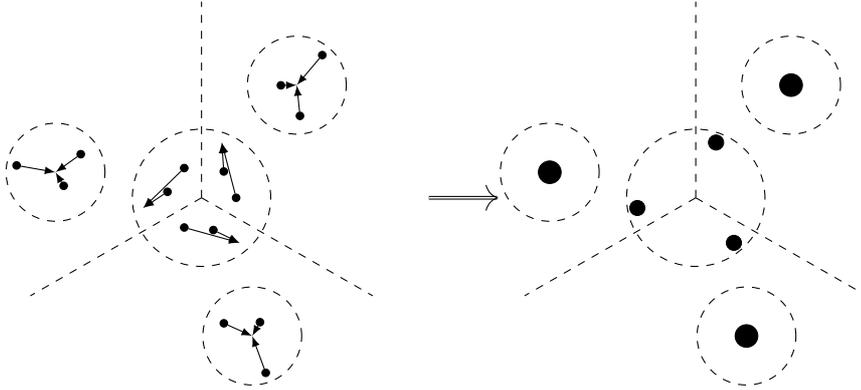
\begin{center}\begin{tikzpicture}[scale=1.3]

\draw [dashed] (0,0)--(90:2);
\draw [dashed] (0,0)--(90+120:2);
\draw [dashed] (0,0)--(90+240:2);

\node[draw,circle,inner sep=1pt,fill] (P11) at (90-40:1.9) {};
\node[draw,circle,inner sep=1pt,fill] (P12) at (90-35:1.4) {};
\node[draw,circle,inner sep=1pt,fill] (P13) at (90-50:1.3) {};
\draw[->,>=latex] (P11) to (90-40:1.5);
\draw[->,>=latex] (P12) to (90-40:1.5);
\draw[->,>=latex] (P13) to (90-40:1.5);
\draw [dashed] (90-40:1.5) circle (0.5);

\node[draw,circle,inner sep=1pt,fill] (P11) at (120+90-40:1.9) {};
\node[draw,circle,inner sep=1pt,fill] (P12) at (120+90-35:1.4) {};
\node[draw,circle,inner sep=1pt,fill] (P13) at (120+90-50:1.3) {};
\draw[->,>=latex] (P11) to (120+90-40:1.5);
\draw[->,>=latex] (P12) to (120+90-40:1.5);
\draw[->,>=latex] (P13) to (120+90-40:1.5);
\draw [dashed] (120+90-40:1.5) circle (0.5);

\node[draw,circle,inner sep=1pt,fill] (P11) at (240+90-40:1.9) {};
\node[draw,circle,inner sep=1pt,fill] (P12) at (240+90-35:1.4) {};
\node[draw,circle,inner sep=1pt,fill] (P13) at (240+90-50:1.3) {};
\draw[->,>=latex] (P11) to (240+90-40:1.5);
\draw[->,>=latex] (P12) to (240+90-40:1.5);
\draw[->,>=latex] (P13) to (240+90-40:1.5);
\draw [dashed] (240+90-40:1.5) circle (0.5);

\node[draw,circle,inner sep=1pt,fill] (P01) at (90-40:0.35) {};
\node[draw,circle,inner sep=1pt,fill] (P02) at (90-90:0.35) {};
\draw[->,>=latex] (P01) to (90-20:0.6);
\draw[->,>=latex] (P02) to (90-20:0.6);

\node[draw,circle,inner sep=1pt,fill] (P01) at (120+90-40:0.35) {};
\node[draw,circle,inner sep=1pt,fill] (P02) at (120+90-90:0.35) {};
\draw[->,>=latex] (P01) to (120+90-20:0.6);
\draw[->,>=latex] (P02) to (120+90-20:0.6);

\node[draw,circle,inner sep=1pt,fill] (P01) at (240+90-40:0.35) {};
\node[draw,circle,inner sep=1pt,fill] (P02) at (240+90-90:0.35) {};
\draw[->,>=latex] (P01) to (240+90-20:0.6);
\draw[->,>=latex] (P02) to (240+90-20:0.6);

\draw [dashed] (0:0) circle (0.7);

\draw [double,->] (2.3,0) to (3,0);

\begin{scope}[shift=(0:5)]

\draw [dashed] (0,0)--(90:2);
\draw [dashed] (0,0)--(90+120:2);
\draw [dashed] (0,0)--(90+240:2);

\node[draw,circle,inner sep=3pt,fill] (P11) at (90-40:1.5) {};
\draw [dashed] (90-40:1.5) circle (0.5);
\node[draw,circle,inner sep=3pt,fill] (P11) at (120+90-40:1.5) {};
\draw [dashed] (120+90-40:1.5) circle (0.5);
\node[draw,circle,inner sep=3pt,fill] (P11) at (240+90-40:1.5) {};
\draw [dashed] (240+90-40:1.5) circle (0.5);

\node[draw,circle,inner sep=2pt,fill] (P01) at (90-20:0.6) {};
\node[draw,circle,inner sep=2pt,fill] (P01) at (120+90-20:0.6) {};
\node[draw,circle,inner sep=2pt,fill] (P01) at (240+90-20:0.6) {};

\draw [dashed] (0:0) circle (0.7);
\end{scope}

\end{tikzpicture}\end{center}

These paths admit unique lifts $\gamma_x$ and $\gamma_y$ in $U$ starting at $x$ and $y$, respectively. The paths $\gamma_x$ and $\gamma_y$ both lie in $\Uu^{\mu_d}$ by Lemma \ref{lem:caractérisation_v/wmu_d}, thus they are homotopically trivial and they only impact the cyclic label. We claim that $\gamma_x(1)=\gamma_y(1)$. Since the cyclic supports of $\gamma_x(1)$ and $\gamma_y(1)$ are equal by construction, it only remains to show that $\clbl(\gamma_x(1))=\clbl(\gamma_y(1))$. 

If the standard point of which $U$ is a confining neighborhood is $0$, then there are no outer disks in $U$, and we have $\clbl(\gamma_x(1))=\cc(\gamma_x(1))=\cc(x)=\cc(y)=\cc(\gamma_y(1))=\clbl(\gamma_y(1))$. If the standard point of which $U$ is a confining neighborhood is nonzero, then there is exactly one outer disk in $U$ per sector $P_i$ for $i\in \intv{1,p}$. By construction, we then have
\[\clbl(\gamma_x(1))=(\alpha_0,\beta_0,\ldots,\alpha_{p-1},\beta_{p-1}) \text{~and~}\clbl(\gamma_y(1))=(\alpha_0',\beta'_0,\ldots,\alpha_{p-1}',\beta'_{p-1})\]
where $\alpha_i^{(')}$ corresponds to the point of the cyclic support of $\gamma_{x,y}(1)$ in the central disk in the sector $P_{i+1}$, and $\beta_i^{(')}$ corresponds to the point of the cyclic support of $\gamma_{x,y}(1)$ in the outer disk in the sector $P_{i+1}$. By Lemma \ref{lem:terms_of_outer_points_constant_in_confining_neighborhood}, we have $\beta_i=\beta'_i$ for $i\in \intv{0,p-1}$. For $i\in \intv{0,p-1}$, the $i$-th term of $\cc(x)$ (resp. $\cc(y))$ is given by $\alpha_i\beta_i$ (resp. $\alpha'_i\beta'_i)$, We obtain that $\alpha_i=\alpha_i'$ for $i\in \intv{0,p-1}$. By Proposition \ref{prop:trivialization_lll}, we have $\gamma_x(1)=\gamma_y(1)$, thus $x$ and $y$ lie in the same path connected component of $U\cap \Uu^{\mu_d}$.

Furthermore, the motions defined above define a deformation-retraction of any path-connected component of $U\cap \Uu^{\mu_d}$ onto a single point, which proves the contractibility.
\end{proof}

We are also able to determine the image of the cyclic content when restricted to $U\cap \Uu^{\mu_d}$.

\begin{cor}\label{cor:bijection_connected_components}
Let $U$ be a confining neighborhood of a standard point $z\in (V/W)^{\mu_d}$. Let also $b$ be the first term of $\olbl(z)$. The cyclic content restricts to a bijection $\pi_0(U\cap \Uu^{\mu_d})\simeq \{u\in D_p^q(c)~|~b\preccurlyeq u_1\}$.
\end{cor}
\begin{proof}
We know by Proposition \ref{prop:conn_compo=cyc_cont} that the cyclic content induces a bijection between $\pi_0(U\cap \Uu^{\mu_d})$ and the  set
\[\{u\in D_p^q(x)~|~ \exists x\in U\cap \Uu^{\mu_d},\cc(x)=u\},\]
we just have to show that this latter set is equal to $\{u\in D_p^q(c)~|~b\preccurlyeq u_1\}$. First, let $x\in U\cap \Uu^{\mu_d}$. We have $b\preccurlyeq \cc(x)_1$ by Lemma \ref{lem:terms_of_outer_points_constant_in_confining_neighborhood}. Conversely, let $u\in D_p^q(c)$ with $b\preccurlyeq u_1$. We write $u_1=\alpha_1 b$. By \cite[Proposition 1.30]{springercat}, there is a unique $\sigma\in D_{2p}^{2q}(c)$ such that $\sigma_1=\alpha$ and $\sigma_2=b$. We consider the standard image $x_\sigma$ of $\sigma$ in $(X/W)^{\mu_d}$ \cite[Definition 11.20]{beskpi1}. Let $\gamma:[0,1]\to E_n$ be the path starting at $\LLL(x_\sigma)$ and consisting in shrinking the first point in each sector towards $0$, as in Lemma \ref{lem:shrinkinglemma}. Let $\ttilde{\gamma}$ be the unique lift of $\gamma$ in $V/W$ obtained by Proposition \ref{prop:lift_lll}. Again, Lemma \ref{lem:caractérisation_v/wmu_d} proves that the path $\ttilde{\gamma}$ actually lies in $(V/W)^{\mu_d}$. By Lemma \ref{lem:shrinkinglemma}, we have $\clbl(\ttilde{\gamma}(1))=\clbl(z)$ and $\ttilde{\gamma}(1)=z$ by Proposition \ref{prop:trivialization_lll}. By definition, there is some point $\ttilde{\gamma}(t)$ which lies in $U$. Since $\cc(\gamma_\sigma(t))=\cc(x_\sigma)=u$, we have the other inclusion.
\end{proof}

Let $U$ be a confining neighborhood of a standard point in $(V/W)^{\mu_d}$. By Proposition \ref{prop:conn_compo=cyc_cont}, the set $U\cap \Uu^{\mu_d}$ can be used as a basepoint (with several path connected components) for the space $U\cap (X/W)^{\mu_d}$.

\begin{definition}
Let $U$ be a confining neighborhood of a standard point $(V/W)^{\mu_d}$. We define the \nit{relative local braid groupoid} associated to $U$ by 
\[\Bb_U:=\pi_1(U\cap (X/W)^{\mu_d},U\cap \Uu^{\mu_d}).\]
\end{definition}

By Corollary \ref{cor:bijection_connected_components}, we can see $\Ob(\Bb_U)$ as a subset of $D_p^q(c)\simeq \pi_0(\Uu^{\mu_d})= \Ob(\Bb)$, where $\Bb$ is the Springer groupoid. The inclusion $U\cap (X/W)^{\mu_d}\hookrightarrow (X/W)^{\mu_d}$ then give a natural functor $\Bb_U\to \Bb$, which is injective on objects. As we will see later, it is also injective on morphisms, but at this stage we are not supposed to know this.

At this point, it is unclear whether or not the relative local braid groupoid $\Bb_U$ depends on the choice of $U$, and whether or not the functor $\Bb_U\to \Bb$ sends $\Bb_U(u,u)$ to a parabolic subgroup of $\Bb(u,u)$. The next Section is devoted to the study of these problems.

\subsection{Standard parabolic subgroupoids}
Our goal is now to obtain a description of parabolic subgroups up to conjugacy. To do this we need a Garside-theoretic version of the definition of relative local braid groupoid. In this section we detail the combinatorial counterparts of the topological definitions given in the last Section.

\subsubsection{Garside definition of Springer groupoids} First, we give some quick preliminaries on Garside categories. Following \cite[Definition 2.4]{besgar}, we define a \nit{homogeneous Garside category} as a triple $(\Cc,\phi,\Delta)$, where $\Cc$ is a homogeneous category, endowed with an automorphism $\phi$ and a natural transformation $\Delta:1_\Cc\Rightarrow \phi$, which satisfies suitable properties. In particular, the set
\[\Ss:=\bigsqcup_{u\in \Ob(\Cc)} \{f\in \Cc(u,-)~|~f\text{ is a prefix of } \Delta(u)\}\]
is finite, its elements are called the \nit{simple morphisms}. For $u,v\in \Ob(\Cc)$, we will denote by $\Ss(u,-)$ (resp. $\Ss(-,v)$, $\Ss(u,v)$) the elements of $\Ss$ with source $u$ (resp. with target $v$, with source $u$ and target $v$). The enveloping groupoid $\Gg$ of $\Cc$ is called a \nit{Garside groupoid}. Like in \cite[Definition V.1.38]{ddgkm}, for $m\geqslant 0$, we denote by $\Delta(u)^m$ the composition
\[\xymatrixcolsep{3pc}\xymatrix{u\ar[r]^-{\Delta(u)} & \phi(u) \ar[r]^-{\Delta(\phi(u))} &\cdots \ar[r]^-{\Delta(\phi^{m-1}(u))}& \phi^{m}(u)}.\]
For $m<0$, we denote by $\Delta(u)^m$ the inverse of $\Delta(\phi^{m}(u))^{-m}$ in $\Gg(u,\phi^m(u))$. In the conventions of \cite{ddgkm}, what we call a homogeneous Garside category $(\Cc,\phi,\Delta)$ is a homogeneous cancellative category $\Cc$, endowed with a Garside map $\Delta$ such that the set $\Ss:=\Div(\Delta)$ is finite, and $\phi$ is the functor $\phi_\Delta$ of \cite[Proposition V.1.28]{ddgkm}. By definition, we have that, for every $u\in \Ob(\Cc)$, the sets $\Ss(u,-)$ and $\Ss(-,u)$ are lattices when endowed with the relations ``being a prefix of'' and ``being a suffix of'', respectively. Furthermore, these relations can be extended to $\Gg(u,-)$ and $\Gg(-,u)$ by setting $f\preccurlyeq g$ if $f^{-1}g\in \Cc$ and $f\preccurlyeq^\Lsh g$ if $gf^{-1}\in \Cc$. The sets $\Gg(u,-)$ and $\Gg(-,u)$ endowed with these relations are also lattices. We will denote by $a\vee b$ and $a\wedge b$ the join and meet of $a,b\in \Gg(u,-)$, respectively, with respect to $\preccurlyeq$, and we will denote by $a\vee^\Lsh b$ and $a\wedge^\Lsh b$ the join and meed of $a,b\in \Gg(-,u)$, respectively, with respect to $\preccurlyeq^{\Lsh}$.

By definition, a morphism $a\in \Gg$ is \nit{positive} if $a\in \Cc$, it is \nit{negative} if $a^{-1}\in \Cc$. The nontrivial elements of $\Ss$ which do not admit proper prefixes are called \nit{atoms}. Every morphisms $f\in \Gg$ admits a unique \nit{left-weighted factorization} $f=\Delta(u)^ks_1\cdots s_r$, where $s_1,\ldots,s_r$ are proper simple morphisms such that $s_is_{i+1}\wedge \Delta(u_i)=s_i$ for $i\in \intv{1,r-1}$ (where $u_i$ is the source of $s_i$). The integer $\inf(f):=k$ and $\sup(f):=k+r$ are called the \nit{infimum} and the \nit{supremum} of $f$, respectively. Moreover, every morphism $f\in \Gg$ can be written uniquely as a fraction $f=a^{-1}b$, where $a\wedge b$ is trivial. This is the \nit{reduced left-fraction decomposition} of $f$. We will say that $a$ and $b$ are the \nit{left-denominator} and the \nit{left-numerator} of $f$, respectively $a=D_L(f)$ and $b=N_L(f)$.

Now, back to Springer groupoids. Following \cite[Section B.5]{beskpi1} and \cite[Section 1.3]{springercat}, we define
\begin{align*}
O&:=\left\{u\in [1,c]^q~\left|~\begin{cases} p\ell_\Rr(u)=\ell_\Rr(c)=n, \\ u(u^{c^{\eta}})\cdots (u^{c^{(p-1)\eta}})=c.\end{cases}\right.\right\},\\
\Ss &:=\left\{(a,b)\in ([1,c]^q)^2~|~ ab \in O\right \},\\
\mathrm{Rel}&:=\left\{(x,y,z)\in ([1,c]^q)^3~|~ xyz\in O\right \}.
\end{align*}
where $[1,c]^q:=\{x\in [1,c]~|~x^{c^q}=x\}$, and $\eta$ is the integer associated to $(p,q)$ introduced in \cite[Lemma 1.28]{springercat} and \cite[Remark 1.29]{springercat} (it depends only on $p$ and $q$). We will use the sets $O,\Ss$ and $\mathrm{Rel}$ to construct a categorical presentation, which will define a Garside category.

The sets $O,\Ss,\mathrm{Rel}$ are in natural one-to-one correspondence with $D_p^q(c)$, $D_{2p}^{2q}(c)$, $D_{3p}^{3q}(c)$, respectively. More precisely, we have
\begin{prop}\cite[Proposition 1.30]{springercat}

\begin{enumerate}
\item The map $u\mapsto (u,u^{c^\eta},\ldots,u^{c^{(p-1)\eta}})$ induces a bijection between $O$ and $D_p^q(c)$.
\item The map $(a,b)\mapsto (a,b,\ldots,b^{c^{(p-1)\eta}})$ induces a bijection between $\Ss$ and $D_{2p}^{2q}(c)$.
\item The map $(x,y,z)\mapsto (x,\ldots,z^{c^{(p-1)\eta}})$ induces a bijection between $\mathrm{Rel}$ and $D_{3p}^{3q}(c)$.
\end{enumerate}
\end{prop}
We will often implicitly use these identifications from now on. We can endow $\Ss$ with the structure of an oriented graph (with object set $O$) by defining the source (resp. the target) of $(a,b)\in \Ss$ as $ab\in O$ (resp. $ba^{c^\eta}\in O$). 

An element $(x,y,z)\in \mathrm{Rel}$ induces a relation on the oriented graph $\Ss$, stating that the following triangle is commutative
\[\xymatrix{& yzx^{c^\eta} \ar[rd]^-{(y,zx^{c^\eta})}& \\ xyz\ar[ru]^-{(x,yz)} \ar[rr]_-{(xy,z)} & &zx^{c^\eta}y^{c^{\eta}}. }\]
One directly checks that the three couples involved do lie in $\Ss$. Note that, following the usual convention for Garside categories, we compose elements as if they were paths rather than applications, thus the relation induced by $(x,y,z)\in \mathrm{Rel}$ is written $(x,yz)(y,zx^{c^\eta})=(xy,z)$. We identify an element of $\mathrm{Rel}$ with the relation it induces on $\Ss$ from now on.

\begin{definition}\cite[Definition 11.5]{beskpi1} and \cite[Definition 2.26]{springercat} The \nit{Springer category} $\Cc$ is defined by the categorical presentation $\Cc:=\langle \Ss~|~ \mathrm{Rel}\rangle^+$. The enveloping groupoid of $\Cc$, obtained by formally inverting all morphisms in $\Cc$, will be denoted by $\Gg$.
\end{definition}
By construction, we have $\Ob(\Gg)=\Ob(\Cc)=O\simeq D_p^q(c)\simeq \Ob(\Bb)$. For each object $u$ of $\Cc$, we can consider $\Delta(u)=(u,1)\in \Ss(u,-)$. For each $s=(a,b)\in \Ss(u,-)$, the triple $(a,b,1)\in \mathrm{Rel}$ induces the relation $(a,b)(b,a^{c^\eta})=\Delta(u)$, which proves that $s$ is a prefix of $\Delta(u)$. In fact, for $(a,b),(d,e)\in \Ss(u,-)$, \cite[Lemma 1.38]{springercat} gives that $(a,b)$ is a prefix of $(d,e)$ if and only if $a\preccurlyeq d$ in $[1,c]$.

For $s:=(a,b)\in \Ss(u,v)$, we have $(a,b)\Delta(v)=\Delta(u)(a^{c^\eta},b^{c^\eta})$. If we write $\phi(s):=(a^{c^{\eta}},b^{c^\eta})$, we obtain an automorphism $\phi$ of $\Cc$, and $\Delta$ induces a natural transformation from the identity functor of $\Cc$ to $\phi$. 

\begin{theo}\cite[Theorem 4.2 and Theorem 9.5]{besgar} The triple $(\Cc,\phi,\Delta)$ is a homogeneous Garside category, whose simple morphisms are exactly the elements of $\Ss$.
\end{theo}

The length functor on $\Cc$ is given on $\Ss$ by $\ell(a,b):=\ell_\Rr(a)$. Let $s:=(a,b)\in \Ss$, we can consider the standard image $x_s$ of $s$ \cite[Definition 11.20]{beskpi1}. We have $\cc(x_s)=ab$ by definition, and we can consider the circular tunnel $(x_s,\frac{\pi}{ph})$. We denote by $b_s$ the homotopy class of this circular tunnel in the Springer groupoid $\Bb$. By \cite[Lemma 11.26 and Theorem 11.28]{beskpi1}, the map $s\mapsto b_s$ induces a groupoid isomorphism between $\Gg$ and $\Bb$. From now on, we will identify $\Gg$ and $\Bb$ using this isomorphism.

\subsubsection{Garside standard parabolic subgroupoids}
Let $\beta\in [1,c]$ be such that there exists a simple of the form $(x,\beta)\in \Ss$. We define a partial map $\delta_\beta$ from $\Ob(\Cc)$ to $\Cc$ as follows
\begin{enumerate}[\quad $\bullet$]
\item If $\beta\not\preccurlyeq^\Lsh u$, then $\delta_\beta(u)$ is not defined.
\item If $u=\alpha\beta$ for some $\alpha\in [1,c]$, then $\delta_\beta(u):=(\alpha,\beta)\in \Ss(u,\beta \alpha^{c^\eta})$.
\end{enumerate}
The fact that there is a simple morphism of the form $(x,\beta)$ in $\Ss$ ensures that $\delta_\beta(u)$ is defined for at least some object $u$ of $\Cc$ (for instance $u=x\beta$).

\begin{lem}
Let $\beta\in [1,c]$ be such that there is a simple of the form $(x,\beta)\in \Ss$. The partial map $\delta_\beta$ is target injective. That is, if $\delta_\beta(u)$ and $\delta_\beta(v)$ are defined for some $u,v\in \Ob(\Cc)$, then the targets of $\delta_\beta(u)$ and $\delta_\beta(v)$ are distinct if $u$ and $v$ are distinct.
\end{lem}
\begin{proof}
Let $u,v$ be such that $\delta_\beta(u)$ and $\delta_\beta(v)$ are defined. We have $\beta\preccurlyeq^\Lsh u,v$ and we write $u=\alpha\beta,v=\alpha'\beta$. The targets of $\delta_\beta(u)$ and $\delta_\beta(v)$ are $\beta \alpha^{c^\eta}$ and $\beta \alpha'^{c^\eta}$, respectively. If those targets are equal, we get $\alpha^{c^\eta}=\alpha'^{c^\eta}$ by cancellativity, thus $\alpha=\alpha'$ and $u=v$.
\end{proof}

In this context, we define
\[\Div(\delta_\beta):=\bigsqcup_{\substack{u\in \Ob(\Cc)\\ \delta_\beta(u)\text{ defined}}}\{s\in \Ss~|~s\preccurlyeq \delta_\beta(u)\}\text{~and~}\Div^{\Lsh}(\delta_\beta):=\bigsqcup_{\substack{u\in \Ob(\Cc)\\ \delta_\beta(u)\text{ defined}}}\{s\in \Ss~|~s\preccurlyeq^{\Lsh} \delta_\beta(u)\}.\]

\begin{definition}Let $\beta\in [1,c]$ be such that there is a simple of the form $(x,\beta)$ in $\Ss$. The groupoid $\Gg_\beta$ generated by $\Div(\delta_\beta)$ is called a \nit{standard parabolic subgroupoid} of $\Gg$. The category $\Cc_\beta:=\Gg_\beta\cap \Cc$ is called a \nit{standard parabolic subcategory} of $\Cc$.
\end{definition}

Note that the definition of $\delta_\beta$ also makes sense if $\beta$ divides no object of $\Gg$. However, we always get $\Gg_\beta=\varnothing$ in this case. On the other hand, if $\Gg_\beta=\Gg_{\beta'}$ are nonempty, we easily get $\beta=\beta'$. We choose the convention that $\varnothing$ is not a standard parabolic subgroupoid.

We aim to show that $\Cc_\beta$ is naturally endowed with a Garside category structure, with $\delta_\beta$ as its Garside map. To do this we will use general results of \cite{ddgkm}. The following proposition is required in order for these general results to apply.

\begin{prop}\label{prop:standard_parabolic_subgroupoids_in_springer}
Let $\beta\in [1,c]$ be such that there exists a simple of the form $(x,\beta)\in \Ss$. We have
\[\Ss_\beta:=\Div(\delta_\beta)=\Div^\Lsh(\delta_\beta)=\left\{(a,b)\in \Ss~|~ \beta\preccurlyeq b\right\}=\Cc_\beta\cap \Ss.\]
\end{prop}
\begin{proof}
We write $\Ss_\beta:=\{(a,b)\in \Ss~|~\beta\preccurlyeq b\}$. Let $s:=(a,b)\in \Ss(u,v)$. 
\begin{enumerate}[(1)]
\item If $s\in \Div(\delta_\beta)$, then $\delta_\beta(u)$ is defined, and $s\preccurlyeq \delta_\beta(u)$. By definition, we have $u=\alpha\beta$ for some $\alpha\in [1,c]$, and $\delta_\beta(u)=(\alpha,\beta)$. By \cite[Lemma 1.38]{springercat}, $s\preccurlyeq \delta_\beta(u)$ is then equivalent to $a\preccurlyeq \alpha$. We write $ax=\alpha$ and we obtain $ab=u=\alpha\beta=ax\beta$. By cancellativity, we get $b=x\beta$, $\beta\preccurlyeq^\Lsh b$ and $\beta\preccurlyeq b$ since every element of $[1,c]$ is balanced. Thus we have $\Div(\delta_\beta)\subset \Ss_\beta$.
\item If $s\in \Div^\Lsh(\delta_\beta)$, then there is some object $\beta\preccurlyeq^\Lsh v$ such that $s\preccurlyeq^\Lsh \delta_\beta(v')$. If we write $v'=\alpha\beta$, \cite[Lemma 1.38]{springercat} then gives that $a\preccurlyeq^\Lsh\alpha$. Since elements of $[1,c]$ are balanced, we write $\alpha=ya$ and we obtain $\beta y^{c^\eta}a^{c^\eta}=\beta \alpha^{c^\eta}=v=b a^{c^\eta}$. By cancellativity, we get $\beta y^{c^\eta}=b$ and $\beta\preccurlyeq b$. Thus we have $\Div^{\Lsh}(\delta_\beta)\subset \Ss_\beta$.
\item If $s\in \Ss_\beta$, then we have $\beta\preccurlyeq b$, and we can write $\beta y=b=x\beta$. We have $u=ab=ax\beta$. The triple $(a,x,\beta)$ in $\mathrm{Rel}$ induces the relation $(a,b)(x,\beta a^{c^\eta})=(ax,\beta)=\delta_\beta(u)$, thus $s\in \Div(\delta_\beta)$. Likewise, the tuple $(y^{c^{-\eta}},a,\beta)$ induces the relation $(y^{c^{-\eta}},a\beta)(a,b)=(y^{c^{-\eta}}a,\beta)=\delta_{\beta}(y^{c^{-\eta}} a\beta)$, thus $s\in \Div^\Lsh(\delta_\beta)$. We thus have $\Ss_\beta\subset \Div(\delta_\beta)$ and $\Ss_\beta\subset \Div^{\Lsh}(\delta_\beta)$.
\end{enumerate}
It remains to show that $\Div(\delta_\beta)=\Cc_\beta\cap \Ss$. Since $\delta_\beta(u)$ is always a simple element, the inclusion $\Div(\delta_\beta)\subset \Ss\cap\Cc_\beta$ is immediate. For the converse, let $s,t\in \Ss$ be composable and such that $st\in \Cc_\beta$. We write $s=(a,b)$ and $t=(d,e)$. By \cite[Corollary 1.39]{springercat}, the composition $st$ is simple if and only if $dx=b$ for some $x$, in which case the triple $(a,d,x)\in \mathrm{Rel}$ expresses the relation $st=(ad,x)$. We show that $st=(ad,x)\in \Ss_\beta$. 

The proof of \cite[Proposition 3.3]{springercat} gives that, in $[1,c]$, we have
\[b=ab\wedge ba^{c^\eta},~e=de\wedge ed^{c^\eta},~x=ab\wedge ed^{c^\eta}\]
and $b\wedge e \preccurlyeq x$. Conversely, $(a,d,x)\in \mathrm{Rel}$ gives that $dx=b$ and $xa^{c^\eta}=e$, thus $x\preccurlyeq b\wedge e$ and $x=b\wedge e$. Since $s,t\in \Div(\delta_\beta)$, the first part of the proof gives $\beta\preccurlyeq b,e$. We then have $\beta\preccurlyeq x=b\wedge e$ and $st\in \Ss_\beta$. An induction then proves that $\Ss\cap \Cc_\beta\subset \Div(\delta_\beta)$. 
\end{proof}

Let $\beta\in [1,c]$ as above, and let $s\in \Ss_\beta(u,v)$. Since $s\in \Div(\delta_\beta)$, there is some $\bbar{s}$ such that $s\bbar{s}=\delta_\beta(u)$. We have $\bbar{s}\in \Div^\Lsh(\delta_\beta)$ by definition, and $\bbar{s}\in \Div(\delta_\beta)$ by the above proposition. Thus there is some $\varphi_\beta(s)\in \Ss_\beta$ such that $\bbar{s}\varphi_\beta(s)=\delta_\beta(v)$. We deduce that $s\delta_\beta(v)=s\bbar{s}\varphi_\beta(s)=\delta_\beta(u)\varphi_\beta(s)$. The correspondence $s\mapsto \varphi_\beta(s)$ extends to an automorphism of $\Cc_\beta$ (and $\Gg_\beta$), and $\delta_\beta$ induces a natural transformation from the identify functor of $\Cc$ to $\varphi_\beta$.

\begin{prop}\label{lem:parabolic_subgroupoid_strongly_compatible}
Let $\Gg_\beta$ be a standard parabolic subgroupoid of $\Gg$.
\begin{enumerate}[(1)]
\item For $f\in \Cc$, the following assertions are equivalent
\begin{enumerate}[(i)]
\item $f\in \Cc_\beta$.
\item All the terms in the left-weighted factorization of $f$ lie in $\Cc_\beta$.
\end{enumerate}
Under these conditions, the left-weighted factorizations of $f$ in $\Cc$ and in $\Cc_\beta$ are equal.
\item For $f\in \Gg$, the following assertions are equivalent
\begin{enumerate}[(i)]
\item $f\in \Gg_\beta$
\item If $f=a^{-1}b$ is the reduced left-fraction decomposition of $f$ in $\Gg$, then $a,b\in \Gg_\beta$.
\end{enumerate}
Under these conditions, the reduced left-fraction decomposition of $f$ in $\Gg$ and in $\Gg_\beta$ are equal.
\end{enumerate}
Furthermore, the triple $(\Cc_\beta,\varphi_\beta,\delta_\beta)$ is a homogeneous Garside category, and the enveloping groupoid of $\Cc_\beta$ is $\Gg_\beta$.
\end{prop}
\begin{proof}
By Proposition \ref{prop:standard_parabolic_subgroupoids_in_springer} and \cite[Proposition VII.1.33]{ddgkm}, the category $\Cc_\beta$ is a parabolic subcategory of $\Cc$ in the sense of \cite[Proposition VII.1.30]{ddgkm}. In particular, $\Cc_\beta$ is closed by factors in $\Cc$, and \cite[Corollary VII.2.22]{ddgkm} gives that $\Cc_\beta$ is compatible with $\Ss$ in the sense of \cite[Definition VII.2.4]{ddgkm}, which is precisely $(1)$.

Likewise, proving $(2)$ amounts to prove that $\Cc_\beta$ is strongly compatible with $\Ss$ in the sense of \cite[Definition VII.2.23]{ddgkm}. Since $\Cc_\beta$ is closed under factor, it is closed under quotient in the sense of \cite[Definition I.1.6]{ddgkm}. Moreover, $\Ss_\beta=\Div^\Lsh(\delta_\beta)$ is closed under left-lcm by definition. We obtain by induction that $\Cc_\beta$ is closed under left-lcm. The result is then a direct application of \cite[Proposition VII.2.27]{ddgkm}.

Lastly, by \cite[Remark VII.1.34]{ddgkm}, and since $\delta_\beta$ is target injective, we have that $\delta_\beta$ is a Garside map for $\Cc_\beta$ in the sense of \cite[Definition V.2.19]{ddgkm}. Since $\Cc_\beta$ is a subcategory of $\Cc$, it is also homogeneous and cancellative, a length functor is given simply by restricting the one of $\Cc$. The triple $(\Cc_\beta,\varphi_\beta,\delta_\beta)$ is then a homogeneous Garside category in the sense of \cite{besgar}. Furthermore, $(2)$ gives that $\Gg_\beta$ is the enveloping groupoid of $\Cc_\beta$.
\end{proof}

\begin{rem}\label{rem:stable_under_factor}The fact that $\Cc_\beta$ is a parabolic subcategory of $\Cc$ in the sense of \cite[Proposition VII.1.30]{ddgkm} also proves that it is stable under prefix and suffix: if $f\preccurlyeq g$ or $f\preccurlyeq^\Lsh g$, then $g\in \Cc_\beta$ implies $f\in \Cc_\beta$. In particular, the atoms of $\Cc_\beta$ are exactly the atoms of $\Cc$ which lie in $\Cc_\beta$.
\end{rem}

At this point, our concept of standard parabolic subgroupoids may seem arbitrary. However, we will show later on that it coincides with the relative local braid groupoid of the last section. For now, we can easily show that standard parabolic subgroupoids are stable under intersection.

\begin{lem}\label{lem:gcd_delta_beta}
Let $\beta,\beta'\in [1,c]$. If $u\in \Ob(\Cc)$ is such that $\beta,\beta'\preccurlyeq u$, then $\beta\vee\beta'\preccurlyeq u$ and $\delta_{\beta\vee \beta'}(u)=\delta_\beta(u)\wedge \delta_{\beta'}(u)$.
\end{lem}
\begin{proof}
Let $u:=\alpha\beta=\alpha'\beta'$. By definition, we have $u=x(\beta\vee \beta')$, with $x\preccurlyeq \alpha\wedge \alpha'$ (note that, since every element of $[1,c]$ is balanced, the operation $\wedge$ and $\wedge^\Lsh$ coincide on $[1,c]$, as well as $\vee$ and $\vee^\Lsh$). Conversely, we write $\alpha=(\alpha\wedge \alpha')y$ and $\alpha'=(\alpha\wedge \alpha')y'$. We have $u=(\alpha\wedge \alpha')y\beta=(\alpha\wedge \alpha')y'\beta'$. By cancellativity, $y\beta=y'\beta'$ is a common multiple of $\beta$ and $\beta'$. We then have
$u=(\alpha\wedge \alpha')z(\beta\vee\beta')$ for some $z$, and $(\alpha\wedge \alpha')\preccurlyeq x$. By \cite[Lemma 1.38]{springercat}, we then have $\delta_\beta(u)\wedge\delta_{\beta'}(u)=(\alpha,\beta)\wedge(\alpha',\beta')=(\alpha\wedge \alpha',\beta\vee \beta')=\delta_{\beta\vee \beta'}(u)$ as claimed.
\end{proof}

\begin{prop}\label{prop:standard_parabolic_stable_intersection_divided}
Let $\Gg_\beta$ and $\Gg_{\beta'}$ be two standard parabolic subgroupoids of $\Gg$. If the intersection $\Gg_\beta\cap \Gg_{\beta'}$ is nonempty, then it is equal to $\Gg_{\beta\vee\beta'}$ and we also have $\Cc_{\beta}\cap\Cc_{\beta'}=\Cc_{\beta\vee\beta'}$.
\end{prop}
\begin{proof}
First, since $\Gg_\beta\cap \Gg_{\beta'}$ is nonempty, there exists an object $u\in \Ob(\Gg_\beta)\cap \Ob(\Gg_{\beta'})$. By definition, this means that $\beta,\beta'\preccurlyeq u$. By Lemma \ref{lem:gcd_delta_beta}, we have that $\beta\vee\beta'\preccurlyeq u$, thus $\Gg_{\beta\vee\beta'}$ is a well-defined standard parabolic subgroupoid of $\Gg$. Furthermore, for all $u\in \Ob(\Gg_{\beta}\cap\Gg_{\beta'})$, Lemma \ref{lem:gcd_delta_beta} also gives that $\delta_{\beta\vee\beta'}(u)=\delta_\beta(u)\cap\delta_{\beta'}(u)\in \Cc_\beta\cap\Cc_{\beta'}$, thus $\Cc_{\beta\vee\beta'}\subset \Cc_{\beta}\cap\Cc_{\beta'}$.

Conversely, let $f\in \Cc_\beta\cap \Cc_{\beta'}$ have greedy normal form $s_1\cdots s_r$. By Proposition \ref{lem:parabolic_subgroupoid_strongly_compatible}, we have $s_i\in \Cc_\beta\cap \Cc_{\beta'}$ for all $i\in \intv{1,r}$. If $u_i$ is the source of $s_i$ for $i\in \intv{1,r}$, we get
\[\forall i\in \intv{1,r},s_i\preccurlyeq \delta_\beta(u_i)\wedge \delta_{\beta'}(u_i)=\delta_{\beta\vee\beta'}(u_i).\]
In particular, all the $s_i$ lie in $\Cc_{\beta\vee\beta'}$, as well as $f$, we get $\Cc_\beta\cap\Cc_{\beta'}\subset \Cc_{\beta\vee\beta'}$ as claimed.

Now, let $f\in \Gg$ with reduced left-fraction decomposition $f=a^{-1}b$.  By Proposition \ref{lem:parabolic_subgroupoid_strongly_compatible}, we have
\[f\in \Gg_\beta\cap \Gg_\beta'\Leftrightarrow a,b\in \Cc_\beta\cap \Cc_{\beta'}\Leftrightarrow a,b\in \Cc_{\beta\vee\beta'}\Leftrightarrow f=a^{-1}b\in \Gg_{\beta\vee\beta'}.\]
Hence the result.
\end{proof}

Since standard parabolic subgroupoids and standard parabolic categories are stable under intersection (provided that said intersection is nonempty), we can give the following definition.

\begin{definition}\label{def:standard_parabolic_closure}Let $f\in \Gg$. The \nit{standard categorical parabolic closure} $\SCPC(f)$ of $f$ is defined as the intersection of the standard parabolic subgroupoids of $\Gg$ containing $f$. Let $x\in \Gg(u,u)$. The \nit{standard parabolic closure} $\SPC(x)$ of $x$ is given by $\SCPC(x)(u,u)$.
\end{definition}

\begin{rem}
This notion is not the same as the notion of parabolic closure defined in \cite[Corollary VII.1.36]{ddgkm}. In their sense, if $s\in \Ss(u,v)$ is an atom, then $\{1_u,s,1_v\}$ is a parabolic subcategory containing $s$, but it is not a Garside category in the sense of \cite[Definition 2.4]{besgar}.
\end{rem}

In practice, one can easily describe standard categorical parabolic closures. Let $f\in \Cc$ with left-weighted factorization $f=(a_1,b_1)\cdots (a_r,b_r)$. By Proposition \ref{lem:parabolic_subgroupoid_strongly_compatible}, a standard parabolic groupoid $\Gg_\beta$ contains $f$ if and only if $\beta\preccurlyeq b_i$ for all $i\in \intv{1,r}$. Thus $\SCPC(f)=\Gg_{\beta}$ where $\beta=\bigwedge_{i=1}^rb_i$.

 Likewise, if $f\in \Gg$ is written in reduced left-fraction decomposition $f=a^{-1}b$, then a standard parabolic groupoid $\Gg_\beta$ contains $f$ if and only if $a,b\in \Gg_\beta$. Thus, if $\SCPC(a)=\Gg_{\beta_1}$ and $\SCPC(b)=\Gg_{\beta_2}$, then $\SCPC(f)=\Gg_\beta$ where $\beta=\beta_1\wedge \beta_2$.

\begin{prop}\label{prop:ribbons_in_divided}
Let $\Gg_\beta$ be a standard parabolic subgroupoid of $\Gg$, and let $s\preccurlyeq \beta$ in $[1,c]$. For every object $u\in \Ob(\Gg_\beta)$, we have $s\preccurlyeq u$. We define a functor $\psi:\Gg_\beta\to \Gg$ as follows.
\begin{enumerate}[$\bullet$]
\item For $u\in \Ob(\Gg_\beta)$, $\psi(u)$ is defined as the target of $\sigma_u:=(s,s^{-1}u)\in \Ss(u,-)$.
\item For $f\in \Gg_\beta(u,v)$, we define $\psi(f):=\sigma_u^{-1}g\sigma_v\in \Gg(\psi(u),\psi(v))$.
\end{enumerate}
Furthermore, $\psi$ induces an isomorphism of categories between $\Gg_\beta$ and $\Gg_{\beta'}$ with $\beta':=s^{-1}\beta s^{c^\eta}$. This isomorphism restricts to an isomorphism between $\Cc_\beta$ and $\Cc_{\beta'}$. We say that $\psi$ is a \nit{ribbon} between $\Gg_\beta$ and $\Gg_{\beta'}$.
\end{prop}
\begin{proof}
First, an object $u\in \Ob(\Gg)$ is in $\Ob(\Gg_\beta)$ if and only if $\beta\preccurlyeq u$. We then have $s\preccurlyeq \beta\preccurlyeq u$. Thus $\sigma_u$ is well-defined for all object $u$ of $\Gg_\beta$. Furthermore, if $u,v$ are two distinct objects of $\Gg_\beta$, then $\sigma_u$ and $\sigma_v$ have different targets. In particular, $\psi$ is injective on objects. From this we easily deduce that $\psi$ is an isomorphism between $\Gg_\beta$ and $\psi(\Gg_\beta)$. We then have to show that $\psi(\Gg_\beta)=\Gg_{\beta'}$.

We show that $\psi$ restricts to a functor from $\Cc_\beta$ to $\Cc_{\beta'}$. Let $u\in \Ob(\Cc_\beta)$, with $\delta_\beta(u):=(\alpha,\beta)$. By definition, the target of $\sigma_u$ is 
\[\psi(u)=s^{-1}us^{c^{\eta}}=s^{-1}\alpha\beta s^{c^{\eta}}=s^{-1}\alpha s \beta'=\alpha^s\beta'\in \Ob(\Cc_{\beta'}).\]
Let now $(a,b)\in \Ss_\beta(u,v)$ for some $u,v\in \Ob(\Cc_\beta)$. We show that $\psi((a,b))\in \Ss_{\beta'}(\psi(u),\psi(v))$. Since $s\preccurlyeq \beta\preccurlyeq b$, $as=sa^s\in [1,c]$, and we have
\begin{align*}
\psi(a,b)=\sigma_u^{-1}(a,b)\sigma_v&=(s,s^{-1}u)^{-1}(a,b)(s,s^{-1}v)\\
&=(s,s^{-1}u)^{-1}(as,s^{-1}b)\\
&=(s,s^{-1}u)^{-1}(sa^s,s^{-1}b)\\
&=(a^s,s^{-1}bs^{c^{\eta}}).
\end{align*}

This is a simple element, which lies in $\Cc_{\beta'}$ since $\beta'\preccurlyeq s^{-1}bs^{c^{\eta}}$. We also obtain that $\psi$ preserves positivity, thus it induces a functor from $\Cc_\beta$ to $\Cc_{\beta'}$ as claimed. We deduce that $\psi(\Gg_\beta)\subset \Gg_{\beta'}$. Conversely, we show that $\Ss_{\beta'}\subset \psi(\Gg_\beta)$. Let $(a,b) \in \Ss_{\beta'}$. Since $\beta'\preccurlyeq^{\Lsh} b$, and since $ba^{c^{\eta}}\in [1,c]$, we have that $sa\in [1,c]$. The element $({}^sa,sbs^{c^{\eta}})\in \Ss_\beta$ is a preimage of $(a,b)$ by $\psi$. Since $\Ss_{\beta'}$ generates $\Gg_{\beta'}$, this shows that $\psi$ is an isomorphism from $\Gg_\beta$ to $\Gg_{\beta'}$. Lastly, as $\psi$ induces a bijection between $\Ss_\beta$ and $\Ss_{\beta'}$, it restricts to an isomorphism between $\Cc_\beta$ and $\Cc_{\beta'}$.
\end{proof}

\subsubsection{Correspondence between standard parabolic subgroupoids and relative local braid groupoids}\label{sec:correspondence_parabolic_localbraid}
We are now equipped to show that the isomorphism of groupoids $\Gg\simeq \Bb$ sends standard parabolic subgroupoids to relative local braid groupoids and vice versa.

Let $s:=(a,b)\in \Ss$. We consider $x_s$ the standard image of $s$ in $(V/W)^{\mu_d}$ \cite[Definition 11.20]{beskpi1}. 
\begin{enumerate}[\quad $\bullet$]
\item If $a=1$ or if $b=1$, then $\LLL(x_s)$ contains exactly one point in each sector $P_i$ for $i\in \intv{1,p}$, and $\clbl(x_s)=\cc(x_s)$. If $a=1$, then we consider $\gamma_s:[0,1]\to V/W$ to be the constant path equal to $x_s$. If $b=1$, then we consider $\gamma_s:[0,1]\to V/W$ to be the unique lift inside $V/W$ of the path $[0,1]\to E_n$ which sends $t$ to $(1-t)\LLL(x_s)$.
\item If $a\neq 1$ and $b\neq 1$, then $\LLL(x_s)$ contains exactly two points in each sector $P_i$ for $i\in \intv{1,p}$, and $\clbl(x_s)=(a,b,a^{c^{\eta}},\ldots,b^{c^{(p-1)\eta}})$. We then consider the path $\gamma:[0,1]\to E_n$ which consists in sliding the first point in each sector $P_i$ (which corresponds to $a^{c^{(i-1)\eta}}$ in the cyclic label) towards $0$ for $i\in \intv{1,p}$, as in Lemma \ref{lem:shrinkinglemma}. 
\end{enumerate}
In each of these cases, we define $z_s$ as the endpoint of the path $\gamma_s$. It is always a standard point. 

\begin{lem}
Let $s:=(a,b)\in \Ss$. 
\begin{enumerate}[\quad $\bullet$]
\item If $a=1$, then $z_s=x_s$. 
\item If $b=1$, then $z_s=0$. 
\item If $a\neq 1\neq b$, then $z_s$ is such that $\olbl(z_s)=(b,b^{c^\eta},\ldots,b^{c^{(p-1)\eta}})$.
\end{enumerate}
Furthermore, if $s':=(a',b')\in \Ss$ is such that $b=b'$, then $z_s=z_{s'}$. 
\end{lem}
\begin{proof}
The first point is immediate. The second point follows from \cite[Lemma 5.6]{beskpi1}. The third point follows from Lemma \ref{lem:shrinkinglemma}. Now, if $s':=(a',b')\in \Ss$ is such that $b=b'$ then we either have
\begin{enumerate}[\quad $\bullet$]
\item $a=1$, in which case $\ell_\Rr(b)=\ell_\Rr(b')=\frac{n}{p}$, and $a'=1$. Thus $z_s=z_{(1,b)}=z_{(1,b')}=z_{s'}$. 
\item $b=1$, in which case $b=1$, and $z_s=0=z_{s'}$.
\item $a\neq 1$ and $b\neq 1$, in which case $z_s$ and $z_{s'}$ share the same outer label and $\LLL(z_s)=\LLL(z_{s'})$ by construction. Thus, $z_s=z_{s'}$ by Corollary \ref{cor:lll_olbl_caractérisent_égalité}. 
\end{enumerate}
\end{proof}

Using this Lemma, we write $z_b$ instead of $z_s$ for $s=(a,b)\in \Ss$ from now on. The following proposition shows in particular that, in our study of parabolic subgroups up to conjugacy, we can restrict our attention to points of the form $z_b$ for $s=(a,b)\in \Ss$.

\begin{prop}\label{prop:replace_by_z_s}
Let $x\in (V/W)^{\mu_d}$, and let $b$ be the first term of $\cc(x)$. The point $z_b$ lies on the same stratum of the discriminant stratification as $x$. 
\end{prop}
\begin{proof}By Corollary \ref{cor:replace_by_standard_point}, it is sufficient to consider the case where $x$ is a standard point. First, if $x\in (X/W)^{\mu_d}$, then by construction, $z_{b}$ is just the standard image of $(1,b)$, which also lies in $(X/W)^{\mu_d}$, that is, on the same stratum as $x$.

Now, if $x\notin (X/W)^{\mu_d}$, then let $U$ be a confining neighborhood of $x$ in $(V/W)^{\mu_d}$. Since $\Uu^{\mu_d}$ is dense in $(V/W)^{\mu_d}$, we can consider some $y\in U\cap \Uu^{\mu_d}$. Consider the path $\gamma_y:[0,1]\to (V/W)^{\mu_d}$ defined in the proof of Proposition \ref{prop:conn_compo=cyc_cont}. The first two terms of the cyclic label of $\gamma_y(1)$ are of the form $s:=(\alpha,\beta)\in \Ss$. By Lemma \ref{lem:terms_of_outer_points_constant_in_confining_neighborhood}, we have $\beta=b$. By construction, the point $z_\beta=z_b$ is such that $\olbl(z_b)=(b,b^{c^\eta},\ldots,b^{c^{(p-1)\eta}})=\olbl(x)$. We can consider an obvious path $\theta:[0,1]\to E_n$ from $\LLL(x)$ to $\LLL(z_b)$, whose unique lift in $(V/W)^{\mu_d}$ is homotopically trivial, and gives that $x$ and $z_b$ lie on the same stratum of the discriminant hypersurface by Proposition \ref{prop:path_preserves_multiplicity_lll}. \end{proof}

As stated earlier, we can now show that two points of $(V/W)^{\mu_d}$ sharing the same cyclic content belong to the same stratum of the discriminant stratification.
\begin{cor}
Let $x,x'\in (V/W)^{\mu_d}$. If $\cc(x)=\cc(x')$, then $x$ and $x'$ belong to the same stratum of the discriminant stratification.
\end{cor}

Using Proposition \ref{prop:replace_by_z_s}, we show that parabolic subgroups are realized up to conjugacy as fundamental groups of the intersection of $(X/W)^{\mu_d}$ with confining neighborhoods of points of the form $z_b$ for $s=(a,b)\in \Ss$.  

\begin{cor}\label{cor:parabolics_are_confining2}
Let $B_0\subset \Bb(u,u)$ be a parabolic subgroup. There is a point of the form $z_b$ for some $s=(a,b)\in \Ss$, and a confining neighborhood $U$ of $z_b$ in $(V/W)^{\mu_d}$, such that $B_0$ is conjugate in $\Bb$ to the image of $\Bb_U(v,v)$ in $\Bb(v,v)$.
\end{cor}
\begin{proof}
Let $\eta$ be a normal ray such that $B_0$ is the image in $\Bb(u,u)$ of $\pi_1^{\loc}((X/W)^{\mu_d},\eta)$. Up to replacing $B_0$ by a conjugate subgroup, we can replace $\eta(1)$ by another point lying on the same stratum of the discriminant stratification. By Proposition \ref{prop:replace_by_z_s}, we can then assume that $\eta(1)$ has the form $z_b$ for some $s=(a,b)\in \Ss$.

Let $U$ be a neighborhood of $z_b$ in $(V/W)^{\mu_d}$, suitable for defining $\pi_1^{\loc}((X/W)^{\mu_d},\eta)$. Up to taking a subneighborhood of $U$, we can assume that $U$ is a confining neighborhood of $z_b$. Up to replacing $\eta$ and $z_b$ by $e^{i\epsi}\eta$ and $e^{i\epsi}z_b$ for $\epsi>0$ small enough, we can assume that $\eta(t)\in U\cap \Uu^{\mu_d}$ for some $t<1$ suitable for defining $\pi_1^{\loc}((X/W)^{\mu_d},\eta)$. If we denote by $v$ the connected component of $\eta(t)$ in $U\cap \Uu^{\mu_d}$, we get that $\pi_1^{\loc}((X/W)^{\mu_d},\eta)$ is identified with $\Bb_U(v,v)$. The restriction of the normal ray $\eta$ to $[0,t]$ gives a morphism $f$ from $u$ to $v$ in $\Bb$. The conjugation by $f$ in $\Bb$ sends the image of $\Bb_U(v,v)$ in $\Bb(v,v)$ to $B_0$ as claimed.
\end{proof}

At this stage, we still need to show that the functor $\Bb_U\to \Bb$ associated to some confining neighborhood $U$ in $(V/W)^{\mu_d}$ always sends $\Bb_U(u,u)$ to a parabolic subgroup of $\Bb(u,u)$. 

Let $s:=(a,b)\in \Ss$, that we fix from now on. We plan to show that, for any confining neighborhood of $z_b$ in $(V/W)^{\mu_d}$, the isomorphism $\Gg\to \Bb$ restricts to an isomorphism between $\Gg_b$ and $\Bb_U$.

Let $\sigma=(d,e)\in \Ss_b$. By definition of $\Ss_b$, we can write $e=e'b$ in $M(c)$ and consider the standard image $x_0$ of the tuple $(d,e',b)$ in $(V/W)^{\mu_d}$. Let $\gamma$ be the path in $E_n$ starting from $\LLL(x_0)$ and consisting in shrinking the points corresponding to $d$ and to $e'$ in the cyclic label towards $0$ as in Lemma \ref{lem:shrinkinglemma}, and let $\ttilde{\gamma}$ be the unique lift of $\gamma$ in $V/W$ starting from $x_0$. By Lemma \ref{lem:caractérisation_v/wmu_d}, Lemma \ref{lem:shrinkinglemma} and Proposition \ref{prop:trivialization_lll}, $\ttilde{\gamma}$ is a path in $(V/W)^{\mu_d}$ with $\ttilde{\gamma}(1)=z_b$.

For $t\in [0,1]$, let $H_\sigma(t,-)$ be the unique lift in $V/W$ of the path in $E_n$ starting from $\LLL(\gamma(t))$ and consisting in rotating the points of the cyclic support of $\gamma(t)$ corresponding to $d$ by an angle of $3\pi/4p$. We obtain a continuous map $H_\sigma:[0,1]\times [0,1]\to (V/W)^{\mu_d}$. By the Hurwitz rule, the path $H_\sigma(t,-)$ represents $\sigma$ in $\Bb$ for all $t<1$. 

Here are two examples with $p=3$, of paths of the form $H_\sigma(0,-)$ and $H_\sigma(1/2,-)$

\begin{center}
\begin{tikzpicture}[scale=1]
\draw [dashed] (-2,0)--(-2,2);
\draw [dashed] (-2,0)--(-2+2*0.866,-1);
\draw [dashed] (-2,0)--(-2-2*0.866,-1);

\node[draw,circle,inner sep=1pt,fill] (P1) at (-2-1.5*0.5,1.5*0.8660) {};
\node[draw,circle,inner sep=1pt,fill] (P2) at (-2-1.5*0.8660,1.5*0.5) {};
\node[draw,circle,inner sep=1pt,fill] (P3) at (-2-1.5*1,-1.5*0) {};
\node[draw,circle,inner sep=1pt,fill] (P5) at (-2-1.5*0.5,-1.5*0.8660) {};
\node[draw,circle,inner sep=1pt,fill] (P6) at (-2+1.5*0,1.5*-1) {};
\node[draw,circle,inner sep=1pt,fill] (P7) at (-2+1.5*0.5,-1.5*0.8660) {};

\node[draw,circle,inner sep=1pt,fill] (P9) at (-2+1.5*1,1.5*0) {};
\node[draw,circle,inner sep=1pt,fill] (P10) at (-2+1.5*0.8660,+1.5*0.5) {};
\node[draw,circle,inner sep=1pt,fill] (P11) at (-2+1.5*0.5,1.5*0.8660) {};

\draw[->,>=latex,rounded corners=1pt] (P3) to[bend right] (-2-1.5*0.7071,1.5*-0.7071);
\draw[->,>=latex,rounded corners=1pt] (P7) to[bend right] (-2+1.5*0.9659,-1.5*0.2588);
\draw[->,>=latex,rounded corners=1pt] (P11) to[bend right] (-2-1.5*0.2588,1.5*0.9659);

\draw [dashed] (3,0)--(3,2);
\draw [dashed] (3,0)--(3+2*0.866,-1);
\draw [dashed] (3,0)--(3-2*0.866,-1);

\node[draw,circle,inner sep=1pt,fill] (P1) at (3-1.5*0.5,1.5*0.8660) {};
\node[draw,circle,inner sep=1pt,fill] (P2) at (3-0.75*0.8660,0.75*0.5) {};
\node[draw,circle,inner sep=1pt,fill] (P3) at (3-0.75*1,-0.75*0) {};
\node[draw,circle,inner sep=1pt,fill] (P5) at (3-1.5*0.5,-1.5*0.8660) {};
\node[draw,circle,inner sep=1pt,fill] (P6) at (3+0.75*0,0.75*-1) {};
\node[draw,circle,inner sep=1pt,fill] (P7) at (3+0.75*0.5,-0.75*0.8660) {};

\node[draw,circle,inner sep=1pt,fill] (P9) at (3+1.5*1,1.5*0) {};
\node[draw,circle,inner sep=1pt,fill] (P10) at (3+0.75*0.8660,+0.75*0.5) {};
\node[draw,circle,inner sep=1pt,fill] (P11) at (3+0.75*0.5,0.75*0.8660) {};

\draw[->,>=latex,rounded corners=1pt] (P3) to[bend right] (3-0.75*0.7071,0.75*-0.7071);
\draw[->,>=latex,rounded corners=1pt] (P7) to[bend right] (3+0.75*0.9659,-0.75*0.2588);
\draw[->,>=latex,rounded corners=1pt] (P11) to[bend right] (3-0.75*0.2588,0.75*0.9659);
\end{tikzpicture}
\end{center}

\begin{prop-def}\label{propdef:f_sigma}
Let $U$ be a confining neighborhood of $z_b$ in $(V/W)^{\mu_d}$. For all $\sigma\in \Ss_b$, the path $H_\sigma(t,-)$ lies in $U$ for $t<1$ big enough, we denote by $f_\sigma$ the morphism it induces in $\Bb_U$. \\ The map $\Ss_b\to \Bb_U$, $\sigma\mapsto f_\sigma$ extends to a groupoid morphism $\psi:\Gg_b\to \Bb_U$, which fits in a commutative square of functors
\[\xymatrix{\Gg_b\ar@{^(->}[r] \ar[d]_-{\psi} & \Gg \ar[d]^-{\wr} \\ \Bb_U \ar[r] &\Bb}\]
\end{prop-def}
\begin{proof}
We keep the notation from above. Let $\theta\in [0,1]$. By Lemma \ref{lem:shrinkinglemma}, we have $H_\sigma(1,\theta)=z_s$. Since $U$ is a neighborhood of $z_s$, and since $H_\sigma$ is continuous, there is a neighborhood $K_\theta$ of $(1,\theta)$ in $[0,1]\times [0,1]$ that is sent inside of $U$ by $H$. The union $\bigcup K_\theta$ is an open cover of the compact set $\{(1,\theta)~|~\theta\in I\}$. By extracting a finite cover, we see that there is some $t_0\in [0,1[$ such that $H_\sigma(t,-)$ is a path in $U$ for $t>t_0$.

Furthermore, as all the paths $H_\sigma(t,-)$ are homotopic for $t<1$, we see that the morphism $f_\sigma$ does not depend on the choice of $t<1$ big enough.

Since $(\Cc_b,\varphi_b,\delta_b)$ is a Garside category, \cite[Proposition VI.1.11]{ddgkm} gives that $\Cc_b$ and the groupoid $\Gg_b$ are generated by the elements of $\Ss_b$, with the relations induced by the triples of the form $(x,y,z)$ with $b\preccurlyeq z$. Thus, replacing the circular tunnels of \cite[Definition 11.25]{beskpi1} by paths of the form $H_\sigma(t,-)$ for $t<1$ big enough, we can imitate the proof of \cite[Lemma 11.26]{beskpi1} to obtain that the map $\sigma\mapsto f_\sigma$ is compatible with the defining relations of $\Cc_b$ and $\Gg_b$. Thus $\sigma\mapsto f_\sigma$ induces a functor $\psi:\Gg_b\to \Bb_U$. 

Lastly, since the path $H_\sigma(t,-)$ represents $\sigma$ in $\Bb$ for all $t<1$, the image of $f_\sigma$ under the natural functor $\Bb_U\to \Bb$ is $\sigma$, thus the above square is commutative as claimed.
\end{proof}

From this we finally deduce that standard parabolic subgroupoids and relative local braid groupoids are identified by the isomorphism $\Gg\to \Bb$.

\begin{theo}\label{theo:isomorphism_parabolic_subgroupoids}
Let $s=(a,b)\in \Ss$, and let $U$ be a confining neighborhood of $z_b$ in $(V/W)^{\mu_d}$. The functor $\psi:\Gg_b\to \Bb_U$ of Proposition-Definition \ref{propdef:f_sigma} is an isomorphism of groupoids.
\end{theo}
\begin{proof}
The commutative square of functors of Proposition-Definition \ref{propdef:f_sigma} proves in particular that $\psi$ is faithful. To prove that $\psi$ is full, we use the same generic position argument as in the proof of \cite[Theorem 11.28]{beskpi1}: every morphism $f$ in $\Bb_U$ can be represented by a path $\gamma$ such that at any given $t\in [0,1]$, at most one point in $\LLL(\gamma(t))$ lies on the half-line $i\R_{\geqslant 0}$. This expresses $f$ as a composition of paths homotopic to some $f_\sigma$ with $\sigma\in \Ss_b$.

Now, by Proposition \ref{prop:conn_compo=cyc_cont}, we know that $\psi$ is bijective on objects, which terminates the proof.
\end{proof}

In particular, we obtain that the groupoid $\Bb_U$ does not depend on the confining neighborhood $U$, but only on the point $z_b$. Lastly, we show the main theorem of this section.

\begin{theo}\label{theo:parabolics_up_to_conj}
A subgroup $B_0\subset \Bb(u,u)=\Gg(u,u)$ is parabolic if and only if there is a standard parabolic subgroupoid $\Gg_\beta$ of $\Gg$ and some $f\in \Gg(u,v)$ such that $(B_0)^f=\Gg_\beta(v,v)$.
\end{theo}
\begin{proof}
First, let $\Gg_\beta$ be a standard parabolic subgroupoid of $\Gg$. Let also $f\in \Gg(u,v)$ be such that $v\in \Ob(\Gg_\beta)$. By definition we have $v=\alpha\beta$, and we can set $s:=(\alpha,\beta)$. Let $x_s$ be the standard image of $s$ in $(X/W)^{\mu_d}$. The natural path $\eta$ from $x_s$ to $z_\beta$ consisting in shrinking points towards $0$ is a normal ray. The image of $\pi_1^{\loc}((X/W)^{\mu_d},\eta)$ in $\pi_1((X/W)^{\mu_d},x_s)=\Bb(v,v)$ is a parabolic subgroup of the latter by definition.

Let $U$ be a small enough confining neighborhood of $z_\beta$ in $(V/W)^{\mu_d}$. By definition of a local fundamental group, we have identifications
\[\pi_1^{\loc}((X/W)^{\mu_d},\eta)\simeq \pi_1(U\cap (X/W)^{\mu_d},\eta)=\Bb_U(v,v).\]
Thus, the image of $\Bb_U(v,v)$ in $\Gg(v,v)$ is a parabolic subgroup of $\Gg(v,v)$. By Theorem \ref{theo:isomorphism_parabolic_subgroupoids}, this image is $\Gg_\beta(v,v)$, which is then a parabolic subgroup of $\Gg(v,v)$.

Let now $\delta$ be a path from $u$ to $v$ in $(V/W)^{\mu_d}$ representing $f$ in $\Gg$. The isomorphism between $\Gg(u,u)$ and $\Gg(v,v)$ induced by $\delta$ is a mere change of basepoint, and as such, it preserves the set of parabolic subgroups. Thus $B_0:=(\Gg_\beta(v,v))^f$ is a parabolic subgroup of $\Gg(u,u)=\Bb(u,u)$.

Conversely, let $B_0\subset \Bb(u,u)$ be a parabolic subgroup. By Corollary \ref{cor:parabolics_are_confining2}, there is some $\alpha\in \Bb(u,v)$ such that $(B_0)^{\alpha}$ is the image of $\Bb_U(v,v)$ in $\Bb(v,v)$, where $U$ is a confining neighborhood of some $z_\beta$ in $(V/W)^{\mu_d}$. We have by Theorem \ref{theo:isomorphism_parabolic_subgroupoids} that the image of $\Bb_U(v,v)$ is $\Gg_b(v,v)$, which terminates the proof.
\end{proof}

\section{The particular case of $B(G_{31})$}

In this Section, we restrict our attention to the complex braid group $B(G_{31})$. Our aim is to prove the main theorems of \cite{paratresses} for this group, by using its regular embedding inside of the well-generated group $B(G_{37})$. We can then use the associated Springer groupoid and the Garside description of parabolic subgroups obtained in the last Section to adapt the Garside-theoretic arguments of \cite[Section 4,5,6]{paratresses} to this context.

\subsection{Diagrams for $B(G_{31})$}\label{sec:diagrams_b31}
Let $W:=G_{37}$ be the complexified Coxeter group of type $E_8$. The integer $d=4$ is regular for this group. Using the notation of the last Section, we have $h=30$, $d=4$, $q=15$ and $p=2$. The integer $\eta$ of \cite[Lemma 1.28]{springercat} is $-7$ in this context. But, as the $15$-th power of a Coxeter element in $G_{37}$ is central, we can replace $\eta$ with $8$. By \cite[Theorem 13.3]{beskpi1}, we have $|O|=88$ and $|\Ss|=2691$ and $|\mathrm{Rel}|=16359$ in this case. The data we use for computation is that of \cite[Section 4.2 and Section 4.3]{springercat}. 

We denote by $\Bb_{31}$ the Springer groupoid attached to this data. We denote by $\Cc_{31}\subset \Bb_{31}$ the associated Springer category.

Let $u$ be on object of $\Bb_{31}$. Following \cite[Definition 3.6]{springercat}, we associate to $u$ the (finite) set of \nit{atomic loops} in $\Cc_{31}(u,u)$. By \cite[Theorem 3.31]{springercat}, an element $\sigma$ of $\Bb_{31}(u,u)$ is a braided reflection if and only if there is an object $v$ of $\Bb_{31}$ such that $\sigma$ is conjugate in $\Bb_{31}$ to an atomic loop in $\Bb_{31}(v,v)$. The atomic loops in $\Cc_{31}(u,u)$ always form a generating set of $\Bb_{31}(u,u)$ by \cite[Proposition 4.7]{springercat}.

By \cite[Theorem 4.2 and Section 4.3.1]{springercat}, there is an object $u_0$ in $\Bb_{31}$ which has five atomic loops $s,t,u,v$ and $w$. These atomic loops generate $\Bb_{31}(u_0,u_0)\simeq B(G_{31})$ with the following presentation 
\begin{equation}\label{eq:presentation_b31}B(G_{31})\simeq \Bb_{31}(u_0,u_0)=\left\langle s,t,u,v,w\left| \begin{array}{l} st=ts,~vt=tv,~wv=vw,\\ suw=uws=wsu,\\svs=vsv,~vuv=uvu,~utu=tut,~twt=wtw\end{array}\right. \right\rangle.\end{equation}
This presentation is summarized in the following diagram \cite[Table 3]{bmr}.
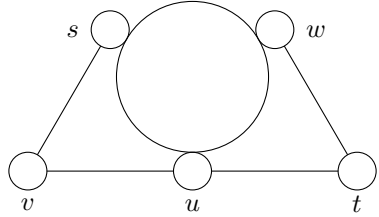
\begin{center}
\begin{tikzpicture}[scale=0.25]
\draw (0,0) circle (4);

\draw (30:5) circle (1);
\node[right] at (24:6) {$w$};
\node[left] at (156:6) {$s$};
\node[below] at (270:6) {$u$};
\node[below] at (-8.66,-5.96) {$v$};
\node[below] at (8.66,-5.8) {$t$};

\draw (150:5) circle (1);
\draw (270:5) circle (1);
\draw (-8.66,-5) circle (1);
\draw (8.66,-5) circle (1);
\draw (-7.66,-5)--(-1,-5);
\draw (1,-5)--(7.66,-5);
\draw (-8.16,-4.134)--(-4.83,1.634);
\draw (8.16,-4.134)--(4.83,1.634);
\end{tikzpicture}\end{center}

Let us denote by $\bbar{s},\bbar{t},\bbar{u},\bbar{v},\bbar{w}$ the images of $s,t,u,v,w$ in the reflection group $G_{31}$. The elements $\bbar{s},\bbar{t},\bbar{u},\bbar{v},\bbar{w}$ generate $G_{31}$, and a presentation is obtained from Presentation (\ref{eq:presentation_b31}) by adding the relations $\bbar{s}^2=\bbar{t}^2=\bbar{u}^2=\bbar{v}^2=\bbar{w}^2=1$. It is known that every parabolic subgroup of $G_{31}$ is, up to conjugacy, generated by some subset $S$ of $\{\bbar{s},\bbar{t},\bbar{u},\bbar{v},\bbar{w}\}$. And a presentation is given by taking the relations in Presentation (\ref{eq:presentation_b31}) which only involve the elements in $S$ (plus the relations $r^2=1$ for $r\in S$). The following theorem states that the situation is the same at the level of the braid group $B(G_{31})$.

\begin{theo}\label{theo:lattice_of_parabolics_of_B31}
The lattice of parabolic subgroups of $\Bb_{31}(u_0,u_0)$ up to conjugacy is given by
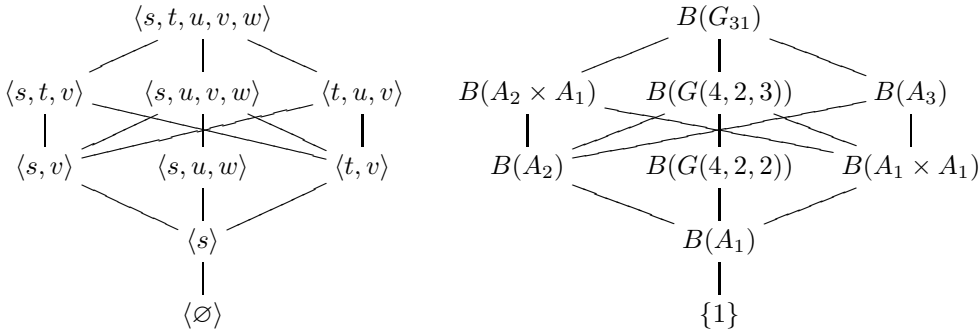
\[\xymatrixcolsep{1pc}\xymatrixrowsep{1pc}\xymatrix{&\langle s,t,u,v,w \rangle & \\
\langle s,t,v\rangle \ar@{-}[ur] & \langle s,u,v,w\rangle \ar@{-}[u]& \langle t,u,v\rangle\ar@{-}[ul]\\
\langle s,v\rangle \ar@{-}[u] \ar@{-}[ur]\ar@{-}[urr]&\langle s,u,w\rangle \ar@{-}[u]&\langle t,v\rangle\ar@{-}[ull] \ar@{-}[ul] \ar@{-}[u]\\
& \langle s\rangle\ar@{-}[ul] \ar@{-}[u] \ar@{-}[ur]& \\
& \langle \varnothing\rangle \ar@{-}[u]&}\hspace{0.5cm}
\xymatrix{&B(G_{31}) & \\
B(A_2\times A_1) \ar@{-}[ur] & B(G(4,2,3)) \ar@{-}[u]& B(A_3)\ar@{-}[ul]\\
B(A_2) \ar@{-}[u] \ar@{-}[ur]\ar@{-}[urr]&B(G(4,2,2)) \ar@{-}[u]& B(A_1\times A_1)\ar@{-}[ull] \ar@{-}[ul] \ar@{-}[u]\\
& B(A_1)\ar@{-}[ul] \ar@{-}[u] \ar@{-}[ur]& \\
& \{1\}\ar@{-}[u]&}
\]
The lattice on the right, where $A_n$ denotes the complex reflection group $G(1,1,n+1)$, describes the isomorphism type of the parabolic subgroups given on the left. For each such parabolic subgroup $\langle S\rangle\subset \Bb_{31}(u_0,u_0)$, a presentation is given by taking the relations in Presentation (\ref{eq:presentation_b31}) which only involve elements of $S$.
\end{theo}
\begin{proof}
First, we show that the considered subgroups are indeed parabolic. Let $S$ be one of the following families of atomic loops in $\Bb_{31}(u_0,u_0)$
\[\{s,t,v\},\{s,u,v,w\},\{t,u,v\}.\]
We denote by $\Gg_b$ the intersection of all the standard parabolic subgroupoids of $\Bb_{31}$ which contain $S$. We check by direct computations that $S$ is exactly the set of atomic loops of $\Bb_{31}(u_0,u_0)$ which lie in $\Gg_b(u_0,u_0)$. We then use the method and algorithms detailed in \cite[Section 4.1]{springercat} on the presentation of $\Gg_b$ to show that $S$ generates $\Gg_b(u_0,u_0)$ with the required presentation. As in \cite[Lemma 4.6]{springercat}, we choose a Schreier transversal $T$ \cite[Definition A.10]{springercat} for $\Gg_b$ rooted in $u_0$ such that the elements of $S$ appears as generators of the presentation of $\Gg_b(u_0,u_0)$ induced by $T$. We then prove, using \cite[Algorithm 4.3, 4.4 and 4.5]{springercat}, that $S$ generates $\Gg_b(u_0,u_0)$, and that all the relators of the presentation of $\Gg_b(u_0,u_0)$ induced by $T$ are deduced from the relators of the required presentation.

The three groups $\langle s,t,v\rangle$, $\langle s,u,v,w\rangle$ and $\langle t,u,v\rangle$ are thus parabolic subgroups of $\Bb_{31}(u_0,u_0)\simeq B(G_{31})$ (and a presentation is given by taking the relations in Presentation (\ref{eq:presentation_b31}) which only involve elements of $S$).

The subgroup $\langle s,t,v\rangle$ is a braid group of type $B(A_2\times A_1)$. It contains the subgroups $\langle s,t\rangle, \langle s,v\rangle, \langle s\rangle$ and $\langle \varnothing \rangle$ as parabolic subgroups by \cite[Proposition 3.1]{paratresses} (we deduce in particular these groups have the desired presentation). These groups are then also parabolic subgroups of $\Bb_{31}(u_0,u_0)\simeq B(G_{31})$ by \cite[Proposition 2.5]{paratresses}. Likewise, $\langle s,u,v,w\rangle$ is a parabolic subgroup of type $B(G(4,2,3))$, which contains $\langle s,u,w\rangle$ as a parabolic subgroup. The group $\langle s,u,w\rangle$ is then a parabolic subgroup of $\Bb_{31}(u_0,u_0)$ again by \cite[Proposition 2.5]{paratresses}, and it admits a presentation given by taking the relations in Presentation (\ref{eq:presentation_b31}) which only involve elements of $S$.

The groups we consider are thus all parabolic subgroups of $\Bb_{31}(u_0,u_0)\simeq B(G_{31})$. Furthermore, the images of these groups in $G_{31}$ form a complete set of representatives for the parabolic subgroups of $G_{31}$ under conjugacy (see \cite[Table C.12]{orlikterao}). By \cite[Proposition 2.6]{paratresses}, the groups we consider then form a complete set of representatives for the parabolic subgroups of $B(G_{31})$ under conjugacy. Lastly, the lattice of parabolic subgroups of $B(G_{31})$ up to conjugacy is the same as the one of $G_{31}$ by \cite[Proposition 2.5]{paratresses}, which terminates the proof.

\end{proof}

We now plan to use this theorem to show that a standard parabolic subgroup in $\Bb_{31}$ is always generated by the atomic loops it contains. We begin by showing that ribbons preserve atomic loops.

\begin{lem}\label{rem:atoms_and_atomic_loops_in_parabolic}
Let $\Cc_\beta$ be a standard parabolic subcategory of $\Bb_{31}$, and let $s\preccurlyeq \beta$ in $[1,c]$. The ribbon $\Cc_\beta\to \Cc_{\beta'}$ induced by $s$ sends the atomic loops in $\Cc_\beta$ to atomic loops in $\Cc_{\beta'}$. Furthermore, every atomic loop in $\Cc_\beta$ is the image under $\psi$ of an atomic loop of $\Cc_{\beta}$.
\end{lem}
\begin{proof}
By \cite[Proposition 3.3, Lemma 3.5, Lemma 4.5]{springercat}, the atomic loops in $\Cc_{31}$ are exactly the composition of two atoms in $\Cc_{31}$ which happen to be endomorphisms. Since the atoms of $\Cc_{31}$ are exactly its elements of length $1$ \cite[Lemma 3.2]{springercat}, atomic loops are exactly the endomorphisms of length $2$ in $\Cc_{31}$. Since the length functor on a standard parabolic subcategory is the restriction of the length functor of $\Cc$, the atomic loops in $\Cc_\beta$ (resp. $\Cc_{\beta'}$) are exactly its endomorphisms of length $2$. 

In the proof of Proposition \ref{prop:ribbons_in_divided}, we saw that $\psi(a,b)=(a^s,s^{-1}bs^{c^\eta})$. Since $\ell_\Rr$ is invariant under conjugacy, we have 
\[\ell(\psi(a,b))=\ell_\Rr(a^s)=\ell_\Rr(a)=\ell(a,b).\]
Thus, ribbons preserve the length function. In particular they preserve endomorphism of length 2, that is, atomic loops.
\end{proof}

We then show that, up to isomorphism induced by ribbons, one can always assume that a standard parabolic subcategory contains the object $u_0$. The following lemma is shown through direct computations.

\begin{lem}\label{lem:ribbons_to_u_0}
Let $\Cc_\beta$ be a standard parabolic subcategory of $\Bb_{31}$. There is a finite sequence of ribbons which sends $\Cc_\beta$ to a standard parabolic subcategory of $\Bb_{31}$ containing $u_0$ as an object.
\end{lem}
\begin{proof}We show that there are two sequences $\beta=\beta_1,\ldots,\beta_n$ and $s_1,\ldots,s_{n-1}$ in $[1,c]$ such that
\begin{enumerate}
\item $s_i\preccurlyeq \beta_i$ and $\beta_{i+1}=s_i^{-1}\beta s_i^{c^8}$ for all $i\in \intv{1,n-1}$,
\item $\beta_n\preccurlyeq u_0$.
\end{enumerate} 
For $i\in \intv{1,n-1}$, the element $s_i$ provides a ribbon $\psi_i:\Cc_{\beta_i}\to \Cc_{\beta_{i+1}}$. The sequence $\psi_{n-1}\circ\cdots\circ\psi_1$ gives the desired result.

The proof that such sequences exist is computational. We start with the set $B:=\{\beta\in [1,c]~|~\beta\preccurlyeq u_0\}$. Then, for every $\beta\in B$, and every $s\preccurlyeq \beta$, we add $s^{-1}\beta s^{c^8}$ to $B$. We iterate this process until no new element can be reached. The fact that every divisor of any $u\in \Ob(\Bb_{31})$ lies in $B$ shows the claim.
\end{proof}

A useful consequence of this lemma is that standard parabolic subgroupoids are always connected. This can also be seen as a corollary of Theorem \ref{theo:isomorphism_parabolic_subgroupoids}.

\begin{lem}\label{lem:standard_parabolic_connected}
Every standard parabolic subgroupoid of $\Bb_{31}$ is connected. That is, if $\Gg_\beta$ is a standard parabolic subgroupoid of $\Bb_{31}$, and $u,v\in \Ob(\Gg_\beta)$, then $\Gg_\beta(u,v)$ is nonempty.
\end{lem}
\begin{proof}
By Lemma \ref{lem:ribbons_to_u_0}, every standard parabolic subgroupoid is isomorphic (through a sequence of ribbons) to a standard parabolic subgroupoid containing $u_0$. There are $20$ such subgroupoids. Direct computations show that these 20 subgroupoids are all connected.
\end{proof}

We are now equipped to show that a standard parabolic subgroup in $\Bb_{31}$ is always generated by the atomic loops it contains.

\begin{prop}\label{prop:atomic loops generate}
Let $\Gg_\beta$ be a standard parabolic subgroupoid of $\Bb_{31}$, and let $u\in \Ob(\Gg_\beta)$. The group $\Gg_{\beta}(u,u)$ is generated by the atomic loops of $\Bb_{31}(u,u)$ which it contains.
\end{prop}
\begin{proof}
Let $S$ be the set of atomic loops of $\Bb_{31}(u,u)$ which lie in $\Cc_\beta(u,u)$. By Lemma \ref{lem:ribbons_to_u_0}, there is an isomorphism of categories $\psi:\Cc_\beta\to \Cc_{\beta'}$ where $\Cc_{\beta'}$ is a standard parabolic subcategory with $u_0\in \Ob(\Cc_{\beta'})$. By Lemma \ref{rem:atoms_and_atomic_loops_in_parabolic}, $\psi$ sends $S$ to the set of atomic loops of $\psi(u)$ which lie in $\Cc_{\beta'}(\psi(u),\psi(u))$. Thus we can replace $\beta$ by $\beta'$ and $u$ by $\psi(u)$ to assume that $u_0\in \Ob(\Cc_\beta)$.

Now, by Lemma \ref{lem:standard_parabolic_connected}, we can consider a morphism $f:u\to u_0$. It is then sufficient to show that $S^f$ generates $\Gg_{\beta}(u_0,u_0)$. Since Theorem \ref{theo:lattice_of_parabolics_of_B31} already gives a generating set $R$ for $\Gg_\beta(u_0,u_0)$, it is sufficient to show that $S^f$ generates $R$, which is done by a case-by-case analysis: every element of $R$ is conjugate by an element of $\langle S^f \rangle$ to an element of $S^f$.
\end{proof}

\subsection{Swaps, recurrent elements and positive conjugators}

Most of the usual solutions to the conjugacy problem in Garside groups have been adapted to the category context in \cite[Chapter VIII]{ddgkm}. In \cite{paratresses} is introduced a new simple procedure to treat theoretical aspects of conjugacy and centralizer in Garside groups. This procedure will be crucial for defining parabolic closure in $\Bb_{31}$. We restate the main results of \cite[Section 4.4]{paratresses} in the general case of a Garside category. We fix a Garside category $(\Cc,\phi,\Delta)$ (in the sense of \cite{besgar}), and we denote by $\Gg$ its enveloping groupoid.

\begin{definition}\label{def:swap}
The \nit{left-swap} function (or just the \nit{swap} function) is the map $\sw$ which sends an endomorphism $x=f^{-1}g$ in $\Gg$ (written as a reduced left-fraction) to $\sw(x)=gf^{-1}$.
\end{definition}
Notice that saying that $x=f^{-1}g$ is an endomorphism amounts to saying that $f$ and $g$ have same target and same source. Thus $gf^{-1}$ is defined and is an endomorphism.

\begin{definition}\label{def:recurrent} Given an endomorphism $x$ in $\Gg$, we say that $x$ is \nit{recurrent for swap} (or just \nit{recurrent}) if $\sw^m(x)=x$ for some $m>0$. The set of recurrent elements conjugate to $x$ is denoted by $\RE(x)$.\end{definition}

\begin{lem}\label{deltrans}
Let $x\in \Gg$ be an endomorphism. If $y\in \mathsf{R}(x)$, then $\phi(y)\in \mathsf{R}(x)$.
\end{lem}
\begin{proof}
Since $\phi$ preserves the lattice structures on $\Gg$, it follows that it preserves the reduced left-fraction decomposition. Hence, $\sw(\phi(z))=\phi(\sw(z))$. That is, applying a left swap commutes with $\phi$.

Given $y\in \mathsf{R}(x)$, there is some $m>0$ such that $\sw^m(y)=y$. By the argument in the previous paragraph, we have $\sw^m(\phi(y))=\phi(\sw^m(y))=\phi(y)$. Therefore $\phi(y)\in \mathsf{R}(x)$.
\end{proof}

These notions are immediate adaptations of the case of Garside groups as considered in \cite[Section 4.4]{paratresses}. By using the same arguments as in \cite[Section 4.4]{paratresses}, we get the following results.

\begin{prop}\label{prop:properties_of_swap}Let $x\in \Gg$ be an endomorphism. 
\begin{enumerate}[(1)]
\item There are integers $0\leqslant m<n$ such that $\sw^m(x)=\sw^n(x)$. The elements in the set $\{\sw^m(x),\ldots,\sw^{n-1}(x)\}$ are recurrent, and this set will be called a \nit{circuit for swap}.
\item If $x$ is conjugate to a positive endomorphism, then $\RE(x)$ is equal to the set of positive conjugates of $x$.
\item If $x$ is conjugate to a negative endomorphism, then $\RE(x)$ is equal to the set of negative conjugates of $x$.
\end{enumerate}
\end{prop}

By definition, $\RE(x)$ coincides with the set of circuits for swap in the conjugacy class of $x$.

\begin{rem}\cite[Remark 4.16]{paratresses}\label{rem:rapidité_swap} If $x$ is conjugate to a positive endomorphism, and $\sw^i(x)=a_i^{-1}b_i$ is the reduced left-fraction decomposition of $\sw^i(x)$ for each $i\geqslant 0$, then $c=a_{m-1}\ldots a_1a_0$ is the minimal morphism in $\Cc$, with respect to $\preccurlyeq^\Lsh$, such that $cxc^{-1}$ is positive (provided $m$ is the smallest nonnegative integer such that $\sw^m(x)$ is positive). In other words, iterated swaps conjugate $x$ to a positive endomorphism in the fastest possible way (conjugating by elements of $\Cc$ on the left). 
\end{rem}

The swap function obviously acts on the elements of $\RE(x)$. The notion of transport for swap gives a way in which the swap function also acts on conjugating morphisms between two elements of $\RE(x)$.

\begin{prop}\cite[Proposition 4.19]{paratresses}\label{prop:transport_for_swap} Let $u,v\in \Ob(\Gg)$, and let $y\in \Gg(u,u),z\in \Gg(v,v)$ which are conjugate by $\alpha\in \Cc(u,v)$, that is, $y^\alpha=z$. Consider the reduced left-fraction decompositions $y=f^{-1}g$ and $z=h^{-1}k$. The \nit{transport of $\alpha$ at $y$} is defined as as $\alpha^{(1)}=f\alpha h^{-1}$, and the following conditions hold:
\begin{enumerate}
\item $\alpha^{(1)}=f\alpha h^{-1}=g\alpha k^{-1}=f\alpha\wedge g\alpha$.
\item $\sw(y)^{(\alpha^{(1)})}=\sw(z)$.
\end{enumerate}
\end{prop}

The transport satisfies very useful properties, in particular it is periodic on $\RE(x)$.

\begin{lem}\cite[Lemma 4.21]{paratresses}\label{lem:swap_agit_sur_R(x)} Let $u\in \Ob(\Gg)$ and let $x\in \Gg(u,u)$. Given two morphisms $y,z\in \mathsf{R}(x)$ and some positive $\alpha\in \Cc$ such that $y^\alpha=z$, then there exists $k>0$ such that $\sw^k(y)=y$ and $\alpha^{(k)}=\alpha$.
\end{lem}

By following the same arguments as in \cite[Section 4.5]{paratresses}, we then obtain the following important property, satisfied by all sets which are used in general to solve the conjugacy problem in Garside categories:

\begin{prop}\cite[Proposition 4.22]{paratresses}\label{prop:gcd_in_R(x)} Let $u\in \Ob(\Gg)$ and let $x\in \Gg(u,u)$. Given $y\in \mathsf{R}(x)$, if $\alpha,\beta\in \Gg(u,-)$ are such that $y^\alpha,y^\beta\in \mathsf{R}(x)$, then $y^{\alpha\wedge \beta}\in \mathsf{R}(x)$.
\end{prop}

Given an endomorphism $x\in \Gg$ which is conjugate to a positive endomorphism, we denote by $C^+(x)$ the set of positive conjugates of $x$. One can compute the set $C^+(x)$ in two steps. First, one finds one element $y\in C^+(x)$ by applying iterated swaps to $x$ (we have $C^+(x)=\RE(x)$ by Proposition \ref{prop:properties_of_swap}). Then, starting with $y$, one computes the directed graph $\Gamma_{C^+(x)}$, defined as follows:
\begin{enumerate}[\quad $\bullet$]
\item The vertices of $\Gamma_{C^+(x)}$ correspond to the elements of $C^+(x)$.
\item There is an arrow labeled $g$ with source $u$ and target $v$ if and only if $g\in \Cc$ is non-trivial, $u^g=v$ and $u^h\notin C^+(x)$ whenever $1\prec h\prec g$.
\end{enumerate}
The arrows in $\Gamma_{C^+(x)}$ are called \nit{minimal positive conjugators} (of $x$). In general a morphism $f\in \Cc(u,-)$ such that $x^f\in \Cc$ will be called a \nit{positive conjugator} (of $x$).

The following Lemma ensures in particular that the graph $\Gamma_{C^+(x)}$ is connected.
\begin{lem}\label{lem:decomposition_positive_conjugator}
Let $x\in \Cc(u,u)$, and $f\in \Cc(u,v)$. If $f$ is a nontrivial positive conjugator of $x$, then there is a path in $\Gamma_{C^+(x)}$ from $x$ to $x^f$ whose composition is $f$.
\end{lem}
\begin{proof}
If $f$ is a minimal positive conjugator, the claim clearly holds. Otherwise, since $f$ is not minimal, it can be decomposed as $f=a_1b_1$, where $a_1$ and $b_1$ are positive with $1\prec a_1\prec f$ and $x^{a_1}\in C^+(x)$. If $a_1$ is not a minimal positive conjugator, we keep going and find $a_2$ such that $1\prec a_2\prec a_1 \prec f$ and $x^{a_2}\in C^+(x)$. Since we cannot have an infinite descending chain of positive elements in a Garside groupoid, this process must stop, and we will obtain a minimal positive conjugator $f_1$ for $x$, with $f=f_1\beta_1$ and $\beta_1\in \Cc$. If $\beta_1$ is not trivial, we apply the same reasoning to $\beta_1$ and find a minimal positive conjugator $f_2$ for $x^{f_1}$, such that $f=f_1f_2\beta_2$ and $\beta_2\in \Cc$. This process must also terminate, since a positive element cannot be decomposed as an arbitrarily large product of positive elements. So some $\beta_r$ will be trivial, and $f=f_1\cdots f_r$ will be a product of minimal positive conjugators, starting at $x$. This shows the claim.
\end{proof}

From now on, we fix $x\in \Cc(u,u)$ an endomorphism. Since the graph $\Gamma_{C^+(x)}$ is connected, one can compute the whole graph, starting from $x$, provided that one knows how to compute the minimal positive conjugators for any given element. This is explained in \cite{paratresses} and in \cite{centralizergarside}. The arguments there adapt directly to the case of a groupoid.

Take an atom $a\in \Cc(u,-)$, that is, a morphism in $\Cc$ with no nontrivial prefix in $\Cc$. The set
\[M_a(x)=\{\alpha\in \Cc(u,-)~|~ a\preccurlyeq \alpha\text{~and~}x^\alpha \in C^+(x)\}=\{\alpha\in \Cc(u,-)~|~a\preccurlyeq \alpha \text{~and~}x^\alpha\in \Cc\}\]
is nonempty ($\Delta(u)\in M_a(x)$), and is closed under $\wedge$ by Proposition \ref{prop:gcd_in_R(x)} since $\RE(x)=C^+(x)$ by Proposition \ref{prop:properties_of_swap}. Hence, as $M_a(x)$ is formed by morphisms in $\Cc$, and every $\preccurlyeq$-chain in $\Cc$ must have a minimal element, it follows that $M_a(x)$ has a unique $\preccurlyeq$-minimal element, that we denote by $\rho_a(x)$.

Every minimal positive conjugator of $x$ must have an atom as a prefix, so it must be equal to $\rho_a(x)$ for some atom $a$. Actually, the set of minimal positive conjugators of $x$ is precisely the set of $\preccurlyeq$-minimal elements in the set $\mathcal M(x)=\{\rho_a(x)~|~ a\in \mathcal A(u,-)\}$, where $\mathcal A$ is the set of atoms of $\Cc$. Since the set $\mathcal M(x)$ is finite, we can compute the arrows starting at $x$ if we are able to compute $\rho_a(x)$ for every atom $a$.

It is known that right-lcms in a Garside category correspond to the classical notion of pushout. If $\alpha\in \Cc(u,v),\beta\in \Cc(u,v')$, then the pushout of $\alpha$ and $\beta$ yields a diagram
\[\xymatrix{ u\ar[r]^-\beta \ar[d]_-{\alpha} &v' \ar[d]^-{\beta'}\\ v\ar[r]_-{\alpha'} & w}\]
where $\alpha\alpha'=\beta\beta'=\alpha\vee\beta$. Concatenating several pushout diagrams produces a new pushout diagram, if one only reads the product of the arrows in the perimeter of the diagram. The element $\alpha'$ is usually denoted by $\alpha\sous \beta$, and it is called the \nit{right-complement of $\alpha$ in $\beta$}. That is, $\alpha\sous \beta=\alpha^{-1}(\alpha\vee\beta)$ and, similarly, $\beta\sous \alpha=\beta^{-1}(\beta\vee\alpha)$. We can then express every pushout diagram as follows (omitting the objects)
\[\xymatrix{\ar[r]^{\beta} \ar[d]_-{\alpha} & \ar[d]^-{\beta\sous\alpha} \\ \ar[r]_-{\alpha\sous\beta} &}\]

Let now $x\in \Cc(u,u)$ be an endomorphism, and let $a\in \Cc(u,-)$ be an atom. As in \cite{paratresses}, we consider the sequence $c_1\preccurlyeq c_2\preccurlyeq \ldots\preccurlyeq \rho_a(x)$ of \nit{converging prefixes} of $\rho_a(x)$, defined by $c_0=a$ and $c_{i+1}=c_i\vee x\backslash c_i$. This sequences stabilizes, and $\rho_a(x)=c_m$ for the smallest $m$ such that $c_m=c_{m+1}$. Notice that all $c_i$ are simple, since $\rho_a(x)\preccurlyeq \Delta(u)$.

If we write $x=u_1\cdots u_r$ as a product of atoms, the computation of $c_{i+1}$ starting from $c_i$ is given by the following concatenation of pushout diagrams.
\[\xymatrixrowsep{3pc}\xymatrixcolsep{3pc}\xymatrix{ \ar[r]^-{u_1} \ar[d]_-{c_j=s_{j,0}} & \ar[r]^-{u_2} \ar[d]_-{s_{j,1}} & \ar@{..>}[r] \ar[d]_-{s_{j,2}} & \ar[r]^-{u_r} \ar[d]^-{s_{j,r-1}} & \ar[d]^-{s_{j,r}=x\backslash c_j} \ar[r]^-{c_j} & \ar[d]\\
\ar[r]_-{v_1} & \ar[r]_-{v_2} & \ar@{..>}[r] & \ar[r]_-{v_r} & \ar[r] & }\]

For every $j\in \intv{0,m-1}$ and every $i\in \intv{0,r}$, the elements $s_{j,i}=(u_1\cdots u_i)\backslash c_j$ will be called the \nit{pre-minimal conjugators} for $a$ and $x$. Notice that the pre-minimal conjugators are not necessarily prefixes of $\rho_a(x)$, but the converging prefixes $c_0,\ldots,c_m$ are.

\subsection{Existence of parabolic closure for $B(G_{31})$}
In this section we show that every element $x$ in $B(G_{31})$ is contained in a parabolic subgroup $\PC(x)$ of $B(G_{31})$ which is minimal with respect to inclusion among all parabolic subgroups of $B(G_{31})$ which contain $x$. Our method is reminiscent of \cite[Section 5.3]{paratresses}.

\begin{lem}\label{lem:sundial}
Let $\Gg_\beta$ be a standard parabolic subgroupoid of $\Bb_{31}$. Let also $\Aa_{\beta}$ (resp. $\Aa^\beta$) be the set of atoms of $\Cc_\beta$ (resp. the set of atoms of $\Cc_{31}\setminus \Cc_\beta$ whose source lies in $\Ob(\Gg_{\beta})$). Let $A_0:=\{(a,b)\in \Aa^\beta~|~a\preccurlyeq \beta\}$.

There is a sequence $A_0\subset \ldots \subset A_r=\Aa^\beta$ such that, for every $i>0$, every atom $s=(a,b):u\to v\in A_i$, and every atom $\sigma\in \Aa_\beta$ with source $u$, we have either
\begin{enumerate}
\item $\sigma\backslash s$ has the form $(a,b')$, where $ab'$ is the target of $\sigma$,
\item $\sigma\backslash s$ is divisible by some $t\in A_{i-1}$.
\end{enumerate}
Furthermore, if $(a,b)\in A_i$, then all atoms of the form $(a,b')$ in $\Aa^\beta$ also lie in $A_i$.
\end{lem}

\begin{rem}\hfill
\begin{enumerate}[(1)]
\item Let $\sigma:=(x,y)$ and $s=(a,b)$ be atoms sharing the same source. Condition $(1)$ in the above lemma is equivalent to stating that $xa$ is simple in $M(c)$. Indeed, if condition $(1)$ is satisfied, we have $\sigma\vee s=(xa,a^{-1}y)$ and $xa$ is simple. Conversely, if $xa$ is simple in $M(c)$, then $xa=x\vee a$ by \cite[Lemma 5.23]{paratresses}. The lcm $s\vee \sigma$ then has the form $(xa,z)$ by \cite[Lemma 1.38]{springercat}, thus $\sigma\backslash s=(a,zx^{c^8})$ and condition $(1)$ holds.
\item By definition, an atom $s$ lies in $A_0$ if and only if some composition of the form $\delta_\beta(u)s$ is simple. By Proposition \ref{prop:ribbons_in_divided}, this is also equivalent to the fact that $s$ induces a ribbon between $\Gg_\beta$ and some other standard parabolic subgroupoid.
\end{enumerate}
\end{rem}

\begin{proof}
For every standard parabolic subgroupoid $\Gg_\beta$, one applies by computer the following pseudocode.
\begin{enumerate}
\item Set $A:=A_0$.
\item Set $\Nn:=\varnothing$.
\item Repeat
\begin{enumerate}[(a)]
\item $A\leftarrow A\cup \Nn$ \hfill \texttt{\# $A$ will be $A_i$ for $i=0,1,\ldots$}
\item For $s:=(a,b)\in \Aa^\beta\setminus A$, with target $u$ set $D_s:=\{\sigma\in \Aa_\beta(u,v)~|~ \sigma\backslash s\neq (a,a^{-1}v)\}$. 
\item Set $\Nn:=\{s\in \Aa^\beta\setminus A~|~\forall \sigma\in D_s, \sigma\backslash s$ is divisible by some $t\in A\}$.
\item $\Nn=\{\sigma=(a,b)\in \Nn~|~ \{s':=(a,b')\in \Aa_\beta\}\subset \Nn\}$ \newline \hfill \texttt{\# Ensures that all atoms with the same first term lie in the same $A_i$}
\end{enumerate}
\item Until $\Nn=\varnothing$.
\item If $A=A^\beta$ then return {\tt true} else return {\tt false}.
\end{enumerate}
The fact that this algorithm terminates in all cases with the statement {\tt true} proves the claim. Since ribbons preserve atoms and lcms, it is actually sufficient to apply this algorithm on a set of representatives for standard parabolic subgroupoids.
\end{proof}

\begin{prop}\label{prop:conj_pos_min}
Let $x\in \Cc_{31}(u,u)$ and $s\in \Ss(u,-)$ an atom such that $\rho_s(x)$ is a minimal positive conjugator. Let also $\Gg_\beta$ be the standard categorical parabolic closure of $x$. We have either
\begin{enumerate}
\item $\rho_s(x)=s$ lies in the set $A_0$ associated to $\Gg_\beta$ in Lemma \ref{lem:sundial},
\item $\rho_s(x)$ lies in $\Cc_\beta$.
\end{enumerate}
\end{prop}

\begin{proof}
We can assume that $\Gg_\beta\neq \Bb_{31}$, otherwise (2) always holds and the result is trivial. Let $A_0\subset A_1\subset\ldots \subset A_r$ be the sequence associated to $\Gg_\beta$ by Lemma \ref{lem:sundial}. Note that $s$ always lies in one of the sets $\Aa^\beta, \Aa_\beta$ described in Lemma \ref{lem:sundial}.  

Let $s:=(a,b)\in A_0$. By Proposition \ref{prop:ribbons_in_divided}, $s$ induces a ribbon $\psi$ between $\Gg_\beta$ and $\Gg_{\beta'}$, where $\beta'=a^{-1}\beta a^{c^8}$. This ribbon restrict to an isomorphism of monoids between $\Cc_\beta(u,u)$ and $\Cc_{\beta'}(\psi(u),\psi(u))$, equal to the restriction to $\Cc_\beta(u,u)$ of conjugation by $s$ in $\Bb_{31}$. We obtain in particular that $x^s=\psi(x)\in \Cc_{\beta'}(\psi(u),\psi(u))$ lies in $\Cc_{31}$. Since $s$ is an atom, and since $x^s\in \Cc_{31}$, we have $\rho_s(x)=s$ and we are in case (1).

Now, suppose that $s\in A_i\setminus A_0$ for some $i>0$. We claim that $\rho_s(x)$ is not minimal. We show this claim by induction on $i$.

Let $s:=(a,b)\in A_1$, and write $x=x_1\cdots x_k$ as a product of atoms, with $x_i=(\alpha_i,\beta_i)$. Since $s\notin A_0$, there is at least one $i\in \intv{1,k}$ such that $a\not\preccurlyeq \beta_i$. Otherwise we have $a\preccurlyeq \bigwedge_{i=1}^k \beta_i=\beta$. Let $j\in \intv{1,k}$ be the first index such that $a\not\preccurlyeq \beta_j$. We compute the pre-minimal conjugators, starting with $s_{0,0}=s$. We write $s_{0,i}:=(a_{i},b_{i})$. If $a\preccurlyeq a_{t}$ for some $t<j$, then $a$ divides the first term of $x_t\backslash s_{0,t}=s_{0,t+1}$. This shows that $a\preccurlyeq a_{j-1}$. Let $s'$ be the atom of the form $(a,b')$ with same source as $x_j$. We have $s'\in A_1$ by definition of $A_1$, and $s'\preccurlyeq s_{0,j-1}$ since $a\preccurlyeq a_{j-1}$ by \cite[Lemma 1.38]{springercat}. Since $a\not\preccurlyeq \beta_j$, $\alpha_j a$ is not simple. Thus $x_j\backslash s'$ admits some $d\in A_0$ as a divisor. Since $d\in A_0$, all subsequent preminimal conjugator admit $d$ as a prefix. We thus have $\rho_d(x)=d\preccurlyeq \rho_s(x)$, which is not minimal.

Now, suppose that $i>1$ and that the claim holds for smaller values of $i$. For $s\in A_i\setminus A_0$, we get using the above argument that some pre-minimal conjugator for $s$ will admit an atom $s'\in A_{i-1}$ as a prefix. Again by Lemma \ref{lem:sundial}, all subsequent pre-minimal conjugators will admit an atom from $A_{i-1}$ as a prefix. Hence, some converging prefix for $s$ and $x$ admits some atom $s'\in A_{i-1}$ as a prefix. It follows that $s'\preccurlyeq \rho_s(x)$, and $\rho_{s'}(x)\preccurlyeq \rho_s(x)$. If $s'\in A_0$, then $\rho_{s'}(x)=s'$ does not admit $s$ as a prefix, thus $\rho_{s'}(x)\neq \rho_s(x)$ and the latter is not minimal. If $s'\notin A_0$, then $\rho_{s'}(x)$ is not minimal by induction hypothesis, and neither is $\rho_s(x)$.

By Lemma \ref{lem:sundial}, this shows that, for an atom $s\notin \Cc_\beta$, we have that $\rho_s(x)$ is minimal if and only if $s\in A_0$, with $\rho_s(x)=s$ in this case. 

Now, if $s\in \Cc_\beta$, the construction of the converging prefixes of $\rho_s(x)$ shows that $\rho_s(x)\in \Cc_\beta$ since standard parabolic subcategories are stable under lcm and under factor. This terminates the proof.
\end{proof}

\begin{prop}\label{prop:conj_pos_min_preserve_supp}
Let $x\in \Cc_{31}(u,u)$ and $y\in \Cc_{31}(v,v)$. Let also $\rho$ be a minimal positive conjugator of $x$, we have $\SPC(x)^\rho=\SPC(x^\rho)$.
\end{prop}
\begin{proof}
For $x\in \Cc_{31}(u,u)$, write $\Delta_x=\delta_\beta(u)$, where $\Gg_{\beta}(u,u)=\SPC(x)$. By definition of the length functor on $\Cc_{31}$, we have that $\ell(\Delta_x)=\ell(u\beta^{-1},\beta)=4-\ell_\Rr(\beta)$ depends only on $\beta$.

By Proposition \ref{prop:conj_pos_min}, we can assume that $\rho$ either is an atom of the form $(a,b)$ with $a\preccurlyeq \beta$ in $[1,c]$, or that $\rho$ lies in $\Cc_{\beta}$. In the first case, we have $\rho=\sigma_u$ is part of a ribbon $\psi:\Cc_{\beta}\to \Cc_{\beta'}$. If we write $x:=s_1\cdots s_r$, with $s_i=(a_i,b_i)$ for $i\in \intv{1,r}$. We have
\[y=x^\rho=\psi(x)=(a_1^a,a^{-1}b_1a^{c^8})\cdots (a_r^a,a^{-1}b_ra^{c^8}).\]
Thus $\SPC(y)=\Gg_{\beta'}(v,v)$ by Definition (see the discussion after Definition \ref{def:standard_parabolic_closure}). Moreover, as $\psi$ is an isomorphism of categories, we have
\[\SPC(x)^\rho=\psi(\Gg_{\beta}(u,u))=\Gg_{\beta'}(v,v)=\SPC(y)\]
as we wanted.

Now, we claim that, for all minimal positive conjugator $\rho$ of any endomorphism $x$, we have $\ell(\Delta_{x^\rho})\leqslant \ell(\Delta_x)$. Indeed, if $\rho\notin \Cc_\beta$, we have showed that $\ell(\Delta_{x^\rho})=\ell(\psi(\Delta_x))=\ell(\Delta_x)$. And, if $\rho\in \Cc_\beta$, we have $x^\rho\in \Cc_\beta$ and $\SPC(x^\rho)\subset \Cc_\beta(v,v)$. Whence $\Delta_{x^\rho}\preccurlyeq \delta_\beta(v)$ and $\ell(\Delta_{x^\rho})\leqslant \ell(\delta_\beta(v))=\ell(\Delta_x)$.

Now, assume that $\rho:u\to v$ is a minimal positive conjugator of $x$ which lies in $\Cc_\beta$. There is a big enough integer $k$ such that $(\delta_\beta(v))^k:v\to v$ is central in $\Gg_\beta(v,v)$. The element $f:=\rho^{-1}(\delta_\beta(v))^k$ is then a positive element which conjugates $x^\rho$ to $x$. Now, $f$ can be written as a composition of minimal positive conjugators by Lemma \ref{lem:decomposition_positive_conjugator}. We showed that the conjugation by a minimal positive conjugator cannot increase the length of the corresponding $\Delta$. This implies that $\ell(\delta_\beta(v))=\ell(\delta_{\beta'}(v))$. It follows that $\delta_{\beta'}(v)=\delta_\beta(v)$, thus $\beta=\beta'$ and, since $\rho\in \Cc_\beta$,
\[\SPC(x)^\rho=\Gg_\beta(u,u)^\rho=\Gg_\beta(v,v)=\Gg_{\beta'}(v,v)=\SPC(x^\rho).\]
\end{proof}

\begin{theo}\label{theo:strong_support_preservingness}
Let $x\in \Bb_{31}(u,u)$ and $y\in \Bb_{31}(v,v)$ be recurrent endomorphisms. Let also $\alpha\in \Bb_{31}(u,v)$ be such that $x^\alpha=y$. One has
\[\SPC(x)^{\alpha}=\SPC(y).\] 
\end{theo}
\begin{proof}
First, if $\alpha$ is not positive, we can replace it by a positive endomorphism $\Delta^k(u)\alpha$, where $\Delta^k(u)$ lies in the center of $\Bb_{31}(u,u)$. Then
\[(\SPC(x))^\alpha=\left((\SPC(x))^{\Delta^{-k}(u)}\right)^{\Delta^k(u)\alpha}=(\SPC(x))^{\Delta^k(u)\alpha},\]
and $x^{\Delta^k(u)\alpha}=x^\alpha=y$. Hence, replacing $\alpha$ with $\Delta^k(u)\alpha$ if necessary, we can assume that $\alpha$ is positive. The proof then proceeds in two steps. 

\underline{Step 1: Assume that $x,y\in \Cc_{31}$}. By Lemma \ref{lem:decomposition_positive_conjugator}, $\alpha$ can be decomposed as a composition of minimal positive conjugators $\alpha_1\cdots \alpha_r$. Since we already showed in Proposition \ref{prop:conj_pos_min_preserve_supp} that standard parabolic closure are preserved under conjugation by minimal positive conjugators, it follows that it will be preserved by $\alpha$
\[(\SPC(x))^\alpha=(\SPC(x))^{\alpha_1\cdots\alpha_r}=(\SPC(x^{\alpha_1}))^{\alpha_2\cdots\alpha_r}=\cdots=\SPC(x^{\alpha_1\cdots\alpha_r})=\SPC(y).\]

\underline{Step 2: Assume that $x,y$ are recurrent}. If $x$ is conjugate to a positive endomorphism, then $y,z\in \RE(x)=C^+(x)$ by Proposition \ref{prop:properties_of_swap}. Hence the result holds  by step 2. If $x$ is conjugate to a negative endomorphism, then $y^{-1},z^{-1}\in \RE(x)^{-1}=C^+(x^{-1})$ by Proposition \ref{prop:properties_of_swap}. Hence $\SPC(y^{-1})^\alpha=\SPC(z^{-1})$. Now, notice that $\SPC(t)=\SPC(t^{-1})$ for every endomorphism $t$ in $\Bb_{31}$, hence $\SPC(y)^{\alpha}=\SPC(z)$, again by step 2.

Finally, suppose that $x$ is conjugate to neither a positive nor a negative element. For every $i\geqslant 0$, let $\sw^i(y)=a_i^{-1}b_i$ (resp. $\sw^i(z)=c_i^{-1}d_i$) be the reduced left-fraction decomposition of $\sw^i(y)$ (resp. of $\sw^i(z)$). We have $y^\alpha=z$. As $\alpha$ is positive, we have
\[(a_0\alpha)^{-1}(b_0\alpha)=\alpha^{-1}y\alpha=z=c_0^{-1}d_0.\]
Consider the transport
\[\alpha^{(1)}=a_0\alpha\wedge b_0\alpha=a_0\alpha c_0^{-1}=b_0\alpha d_0^{-1}.\]
Notice that we have the commutative diagrams of conjugations
\[\xymatrix{z\ar[r]^-{c_0^{-1}} & \sw(y)  & & & z \ar[r]^-{d_0^{-1}}& \sw(z)\\ y\ar[u]^-{\alpha} \ar[r]_-{a_0^{-1}} & \sw(y) \ar[u]_-{\alpha^{(1)}} & & &   y \ar[u]^-{\alpha} \ar[r]_-{b_0^{-1}}& \sw(y) \ar[u]_-{\alpha^{(1)}}}\]
Now recall that both $y$ and $z$ are recurrent, so there exists some $k>0$ such that $\sw^k(y)=y$, $\sw^k(z)=z$ and $\alpha^{(k)}=\alpha$ by Lemma \ref{lem:swap_agit_sur_R(x)}. We obtain the following commutative diagram of conjugations:
\[\xymatrix@C=12mm@R=12mm{
z  \ar[r]^{c_0^{-1}}  & \sw(z) \ar[r]^{c_1^{-1}} & \sw^2(z) \ar@{.>}[r] & \sw^{k-1}(z) \ar[r]^{c_{k-1}^{-1}} & z
\\
y \ar[u]^{\alpha} \ar[r]_{a_0^{-1}} & \sw(y) \ar[u]_{\alpha^{(1)}} \ar[r]_{a_1^{-1}} & \sw^2(y) \ar[u]_{\alpha^{(2)}} \ar@{.>}[r] & \sw^{k-1}(y) \ar[u]_{\alpha^{(k-1)}} \ar[r]_{a_{k-1}^{-1}} & y \ar[u]_{\alpha^{(k)}=\alpha}
}
\]
Simplifying the diagram, we have
\[\xymatrix@C=12mm@R=12mm{
z  \ar[rr]^{(c_{k-1}\cdots c_0)^{-1}}  & & z
\\
y \ar[u]^{\alpha} \ar[rr]_{(a_{k-1}\cdots a_0)^{-1}} & & y \ar[u]_{\alpha}
}\]
Let us denote $g_1=a_{k-1}\cdots a_0$ and $h_1=c_{k-1}\cdots c_0$. Both elements are positive, and we have $\alpha=g_1\alpha h_1^{-1}$. 

Now notice that we also had $\alpha^{(1)}=b_0\alpha d_0^{-1}$. Repeating the above arguments, if we define $g_2=b_{k-1}\cdots b_0$ and $h_2=d_{k-1}\cdots d_0$, we have $\alpha=\alpha^{(k)}=g_2 \alpha h_2^{-1}$. Therefore $\alpha= g_1g_2 \alpha h_2^{-1}h_1^{-1}$.
That is
$$
    (g_1g_2)^\alpha = h_1h_2.
$$
Since $g_1g_2$ and $h_1h_2$ are positive, step 2 gives that $\SPC(g_1g_2)^\alpha=\SPC(h_1h_2)$. The proof will then finish by showing that $\SPC(g_1g_2)=\SPC(y)$ and that $\SPC(h_1h_2)=\SPC(z)$.
Recall that $g_1g_2=a_{k-1}\cdots a_0b_{k-1}\cdots b_0$ where all factors in this expression are positive elements. Hence $a_0,b_0\in \SCPC(g_1g_2)$ and $\SPC(y)\subset\SPC(g_1g_2)$. On the other hand, since $y\in \SPC(y)$, all elements $\{\sw^i(y)\}_{i\geqslant 0}$ belong to $\SCPC(y)$. Hence, all positive elements $a_{k-1},\ldots,a_0,b_{k-1},\ldots,b_0$ belong in $\SCPC(y)$. Therefore $\SPC(g_1g_2)\subset \SPC(y)$, and hence $\SPC(g_1g_2)=\SPC(y)$.

The same argument shows that $\SPC(h_1h_2)=\SPC(z)$, and this finally implies that $\SPC(y)^\alpha=\SPC(z)$, as we wanted to show.
\end{proof}

\begin{theo}\label{theo:spc_is_pc_for_recurrent}
Let $x\in \Bb_{31}(u,u)$ be recurrent for swap. Then $\SPC(x)$ is the smallest parabolic subgroup of $\Bb_{31}(u,u)$ which contains $x$.
\end{theo}
\begin{proof}
Let $\Gg_\beta$ be the smallest standard parabolic groupoid which contains $x$. By definition, we have $\SPC(x)=\Gg_\beta(u,u)$. This group is a parabolic subgroup of $\Bb_{31}(v,v)$ by Theorem \ref{theo:parabolics_up_to_conj}. Conversely, let $B_0\subset \Bb_{31}(u,u)$ be a parabolic subgroup containing $x$. We can assume that $B_0\neq \Bb_{31}(u,u)$, otherwise it is clear that $\SPC(x)\subset B_0$. Since $B_0$ is parabolic, Theorem \ref{theo:parabolics_up_to_conj} gives that there exists a morphism $f\in \Bb_{31}(u,v)$ such that $(B_0)^f=\Gg_{\beta'}(v,v)$ is a standard parabolic subgroup of $\Bb_{31}(v,v)$. If follows that $x^f\in \Gg_{\beta'}(v,v)$. By applying iterated swaps to $x^f$, we obtain that $(x^f)^g$ is recurrent for some $g\in \Gg_{\beta'}(v,w)$. Since $x^f\in \Gg_{\beta'}$, it follows that $x^{fg}\in \Gg_{\beta'}(w,w)$. By definition of the standard parabolic closure, we have $\SPC(x^{fg})\subset \Gg_{\beta'}(w,w)$. Now, since both $x$ and $x^{fg}$ are recurrent, Theorem \ref{theo:strong_support_preservingness} gives
\[(\Gg_\beta(u,u))^{fg}=\SPC(x)^{fg}=\SPC(x^{fg})\subset \Gg_{\beta'}(w,w).\]
Therefore
\[\SPC(x)=\Gg_{\beta}(u,u)\subset (\Gg_{\beta'}(w,w))^{g^{-1}f^{-1}}=(\Gg_{\beta'}(v,v))^{f^{-1}}=B_0.\]
This shows the minimality of $\Gg_{\beta}(u,u)$, hence $\SPC(x)$ is the smallest parabolic subgroup of $\Bb_{31}(u,u)$ which contains $x$.
\end{proof}

\begin{theo}\label{theo:parabolic_closure}
For every element $x\in B(G_{31})$, there exists a unique minimal parabolic subgroup $\PC(x)$ of $B(G_{31})$ containing it, and we have $\PC(x^m)=\PC(x)$ for every $m\neq 0$.
\end{theo}
\begin{proof} Fix $u_0\in \Ob(\Bb_{31})$ such that $B(G_{31})=\Bb_{31}(u_0,u_0)$. There is a morphism $f\in \Bb_{31}(u_0,v)$ such that $y:=x^f\in \Bb_{31}(v,v)$ is recurrent swap (simply apply iterated swaps to $x$). By Theorem \ref{theo:spc_is_pc_for_recurrent}, we know that $\SPC(y)$ is the smallest parabolic subgroup of $\Bb_{31}(v,v)$ which contains $y$. Let $B_0:=(\SPC(y))^{f^{-1}}$.

Let $B_1\subset \Bb_{31}(u_0,u_0)$ be a parabolic subgroup containing $x$. Since conjugating by $f$ corresponds to a change of basepoint in $(X/W)^{\mu_d}$, we know that $(B_1)^f$ is a parabolic subgroup of $\Bb_{31}(v,v)$. Since $y\in (B_1)^f$, we have $\SPC(y)\subset (B_1)^f$. Thus $(\SPC(y))^{f^{-1}}\subset B_1$, and $(\SPC(y))^{f^{-1}}$ is the smallest parabolic subgroup of $\Bb_{31}(u_0,u_0)$ containing $x$.

Let now $m\neq 0$. Up to conjugacy, we can assume that $x$ is recurrent. By definition, $\sw(x^m)$ is obtained by conjugating $x^m$ by $1_{u_0} \wedge x^m$. Since both $1_{u_0}$ and $x^m$ conjugate $x$ to itself in $\Bb_{31}$, Proposition \ref{prop:gcd_in_R(x)} proves that $1_{u_0}\wedge x^m$ conjugate $x$ to a recurrent element. Hence we can conjugate the pair $(x,x^m)$ be $x^m\wedge 1_{u_0}$, to apply $\sw$ to the second coordinate while keeping the first coordinate inside $\RE(x)$. Iterating this we can assume that both $x$ and $x^m$ are recurrent elements.

If both $x$ and $x^m$ are recurrent, then we just need to show that $x\in \SPC(x^m)$ to show that $\PC(x)=\PC(x^m)$. Indeed we have $x^m\in \PC(x)$ by definition, thus we only have to show that $x\in \PC(x^m)=\SPC(x^m)$.

Now, if either $x$ or its inverse lie in $\Cc_{31}$, it is immediate that $\SPC(x)=\SPC(x^m)$. If this is not the case, we consider a Garside map $\Delta^n$ such that both the left-numerator and the left-denominator of $x$ are simple. We can then write $x=s_0^{-1}t_0$ and $\sw^i(x)=t_{i-1}s_{i-1}^{-1}=s_i^{-1}t_i$. One can show that an immediate induction that $x^m=y_0^{-1}y_1^{-1}\cdots y_{m-1}^{-1}x_{m-1}\cdots x_0$. One can then follow the arguments in \cite[Theorem 8.2]{paraartin} and \cite[Theorem 2.9]{rigidity} to conclude that this expression is a reduced left-fraction. In particular, the standard categorical parabolic closure of $x^m$ contains $y_0$ and $x_0$. We then have $x\in \SPC(x^m)$ as we wanted to show.
\end{proof}

This is \cite[Theorem 1.1]{paratresses} in the case of $B(G_{31})$. With the above notation, the group $\PC(x)$ is called the \nit{parabolic closure} of $x$ in $B(G_{31})$. Note that the above proof also gives us the following proposition.

\begin{prop}\label{prop:parabolic_closure_conjugate}
Let $x\in \Bb_{31}(u,u)$ and $f\in \Bb_{31}(u,v)$. We have $\PC(x)^f=\PC(x^f)$.
\end{prop}

\subsection{Intersection of parabolic subgroups}
We show in this section that the intersection of parabolic subgroups of $B(G_{31})$ is again a parabolic subgroup. Our argument is an adaptation of that of \cite[Section 6]{paratresses} to our context. In particular, we also begin by defining, for every parabolic subgroup $B_0$ of $B(G_{31})$, a special element $z_{B_0}\in Z(B_0)$, which encodes the conjugacy class of $B_0$ in $\Bb_{31}$. In the case where $B_0$ is irreducible, we will see that $z_{B_0}$ is a generator of $Z(B_0)$ as claimed in the introduction. We start by defining $z_{B_0}$ when $B_0$ is a standard parabolic subgroup of some $\Bb_{31}(u,u)$. This definition is a straightforward generalization of \cite[Definition 6.1]{paratresses} to our context. 

\begin{definition}
Let $u\in \Ob(\Bb_{31})$, and let $B_0:=\Gg_\beta(u,u)$ be a standard parabolic subgroup of $\Bb_{31}(u,u)$ We define $z_{B_0}=(\delta_\beta(u))^e$, where $e$ is the smallest positive integer such that $(\delta_\beta(u))^e$ lies in the center of $B_0$.
\end{definition}

A difference with the case of a monoid is that the integer $e$ in the above definition depends a priori on the object $u\in \Ob(\Gg_\beta)$. The following lemma ensures that $e$ actually depends only on the standard parabolic subgroupoid $\Gg_\beta$.

\begin{lem}\label{lem:centre_parabolic}
Let $\Gg_\beta$ be a standard parabolic subcategory, and let $k\neq 0$ be an integer. If $(\delta_\beta(u))^k\in Z(\Gg_\beta(u,u))$ for some $u\in \Ob(\Gg_\beta)$, then $(\delta_\beta(v))^k\in Z(\Gg_\beta(v,v))$ for all $v\in \Ob(\Gg_\beta)$.
\end{lem}
\begin{proof}
Since ribbons are isomorphisms which preserve the Garside map, it is sufficient by Lemma \ref{lem:ribbons_to_u_0} to prove the result for standard parabolic subcategories which contain $u_0$. Given such a subcategory $\Cc_\beta$, we directly determine on a computer, for an object $u\in \Ob(\Cc_\beta)$, the smallest integer $k_u$ such that $(\delta_\beta(u))^{k_u}\in Z(\Gg_\beta(u,u))$. This is possible thanks to Proposition \ref{prop:atomic loops generate}: we just have to check that $k_u$ is the smallest positive integer such that $(\delta_\beta(u))^{k_u}\in Z(\Gg_\beta(u,u))$ commutes with the atomic loops lying in $\Gg_\beta(u,u)$.

 We then simply observe that, for all $u\in \Ob(\Cc_\beta)$, the Garside automorphism $\pphi_\beta^{k_u}$ is trivial. This proves that, for all $v\in \Ob(\Cc_\beta)$, $(\delta_\beta(v))^{k_u}$ is central in $\Gg_\beta(v,v)$.
\end{proof}

\begin{lem}\label{lem:conj_pos_min_preserves_zb0}
Let $B_0:=\Gg_\beta(u,u)$ be a standard parabolic subgroup of $\Bb_{31}(u,u)$.  If $\rho:u\to v$ is a minimal positive conjugator of $z_{B_0}$, then $(z_{B_0})^\rho=z_{(B_0)^\rho}$.
\end{lem}
\begin{proof}
By Proposition \ref{prop:conj_pos_min}, we either have $\rho\in \Cc_\beta$ or $\rho=(a,b)$ is an atom of $\Cc_{31}$ such that $a\preccurlyeq \beta$. In the first case, we have $(B_0)^\rho=\Gg_\beta(v,v)$ and the result follows from Lemma \ref{lem:centre_parabolic}. In the second case, there is a ribbon $\psi:\Cc_\beta\to \Cc_{\beta'}$ with $\beta'=a^{-1}\beta a^{c^8}$. Since $\psi$ is an isomorphism of categories which sends $\delta_{\beta}$ to $\delta_{\beta'}$, the smallest power of $\delta_{\beta'}(v)$ which is central in $\Gg_{\beta'}(v,v)$ is equal to $\psi(z_{B_0})$, whence the result.
\end{proof}

Let now $B_0\subset B(G_{31})\simeq \Bb_{31}(u_0,u_0)$ be arbitrary. We want to conjugate $B_0$ in $\Bb_{31}$ to a standard parabolic subgroup $B_1$, and to define $z_{B_0}$ as the conjugate of $z_{B_1}$ in $B_0$. In order for such a definition to be valid, we need to show that it does not depend on the choice of $B_1$. For this we need the following proposition.

\begin{prop}\label{prop:zb0_are_conjugate}
Let $B_1:=\Gg_{\beta_1}(u_1,u_1)$ and $B_2:=\Gg_{\beta_2}(u_2,u_2)$ be two standard parabolic subgroups in $\Bb_{31}$. Let also $f\in \Bb_{31}(u_1,u_2)$ be such that $(B_1)^f=B_2$. We have $(z_{B_1})^f=z_{B_2}$.
\end{prop}
\begin{proof}
Let $f:u_1\to u_2$ be such that $(B_1)^f=B_2$. Up to left-multiplying $f$ by a big enough power of $\Delta(u_1)$ which is central in $\Bb_{31}(u_1,u_1)$, we can assume that $f\in \Cc_{31}(u_1,u_2)$.

Let $z:=(z_{B_1})^f\in B_2\subset \Gg_{\beta_2}$. Since $z$ is conjugate to the positive morphism $z_{B_1}$, we can apply iterated swaps to $z$ so that we get a conjugate $z^g\in \Cc_{31}$ of $z$. By construction of iterated swaps, we have $g\in \Cc_{\beta_2}$ and $z^g\in \Cc_{\beta_2}(v,v)$ for some $v\in \Ob(\Gg_{\beta_2})$.

The morphism $fg:u_1\to v$ is a positive conjugator of $z_{B_1}$. By Lemma \ref{lem:decomposition_positive_conjugator}, we have a sequence of morphisms $z_{B_1}=z_1,\ldots,z_m=z^g$, such that $fg=\rho_1\cdots \rho_{m-1}$, with $\rho_i$ a minimal positive conjugator of $z_i$ and $z_i^{\rho_i}=z_{i+1}$ for $i\in \intv{1,m-1}$.

By an immediate induction using Lemma \ref{lem:conj_pos_min_preserves_zb0}, we have $z^g=(z_{B_1})^{fg}=z_{(B_1)^{fg}}$. Since $(B_1)^{fg}=(B_2)^g$, it remains to show that $z_{B_2}^g=z_{(B_2)^g}$. Since $g$ is a morphism in $\Cc_{\beta_2}$, we have $(z_{B_2})^g=z_{(B_2)^g}$ by Lemma \ref{lem:centre_parabolic}.
\end{proof}

\begin{definition}
Let $B_0$ be a parabolic subgroup of $\Bb_{31}(u,u)$. Let $f\in \Bb_{31}(u,v)$ be such that $(B_0)^f$ is standard, say $(B_0)^f=\Gg_\beta(v,v)$. Then we define $z_{B_0}=(z_{\Gg_\beta(v,v)})^{f^{-1}}$.
\end{definition}

\begin{prop}
Under the above conditions, the element $z_{B_0}$ is well-defined. 
\end{prop}
\begin{proof}
Suppose that $f_1:u\to v_1$ and $f_2:u\to v_2$ are such that $(B_0)^{f_1}=\Gg_{\beta_1}(v_1,v_1)=:B_1$ and $(B_0)^{f_2}=\Gg_{\beta_2}(v_2,v_2)=:B_2$. We need to show that $(z_{B_1})^{f_1^{-1}}=(z_{B_2})^{f_2^{-1}}$. But this follows from Proposition \ref{prop:zb0_are_conjugate} since we have $(B_1)^{f_1^{-1}f_2}=B_2$, and this implies that $(z_{B_1})^{f_1^{-1}f_2}=z_{B_2}$.
\end{proof}

\begin{prop}\label{prop:B_0_parabolic_closure_z_B_0}
Let $B_0$ be a parabolic subgroup of $\Bb_{31}(u,u)$. Then $B_0=\PC(z_{B_0})$.
\end{prop}
\begin{proof}
If $B_0=\Gg_{\beta}(u,u)$ is standard, then $z_{B_0}=(\delta_\beta(u))^k$ for some positive integer $k$. We then have $\PC(z_{B_0})=\SPC(z_{B_0})=\Gg_\beta(u,u)=B_0$ by definition of standard parabolic closure. If $B_0$ is not standard, then let $f$ be such that $B_1:=(B_0)^f$ is standard. We have $B_1=\PC(z_{B_1})$ and $z_{B_0}=(z_{B_1})^{f^{-1}}$ by definition of $z_{B_0}$. By Proposition \ref{prop:parabolic_closure_conjugate}, we have $\PC(z_{B_0})=\PC(z_{B_0}^f)^{f^{-1}}=(B_1)^{f^{-1}}=B_0$.
\end{proof}

\begin{prop}\label{prop:z_B_0_caractérise_conjugaison}
Let $B_1\subset \Bb_{31}(u_1,u_1)$ and $B_2\subset \Bb_{31}(u_2,u_2)$ be two parabolic subgroups. For every $f\in \Bb_{31}(u_1,u_2)$, one has $(B_1)^f=B_2$ if and only if $(z_{B_1})^f=z_{B_2}$.
\end{prop}
\begin{proof}
Suppose that $(B_1)^f=B_2$. Let $g\in \Bb_{31}(u_1,v)$ be such that $(B_1)^g=\Gg_{\beta}(v,v)=:B_0$ is a standard parabolic subgroup. The morphism $f^{-1}g$ is such that $(B_2)^{f^{-1}g}=\Gg_{\beta}(v,v)$ and, by definition, we have $z_{B_1}=(z_{B_0})^{g^{-1}}$ and $z_{B_2}=(z_{B_0})^{g^{-1}f}$. Therefore $(z_{B_1})^f=z_{B_2}$. Conversely, suppose that $(z_{B_1})^f=z_{B_2}$. Then, by Proposition \ref{prop:parabolic_closure_conjugate} and Proposition \ref{prop:B_0_parabolic_closure_z_B_0}, $(B_1)^f=\PC((z_{B_1})^f)=\PC(z_{B_2})=B_2$.
\end{proof}

Using this proposition, we can show that, when $B_0$ is irreducible, $z_{B_0}$ is indeed a generator of the center of $Z(B_0)$.

\begin{prop}Let $u\in \Ob(\Bb_{31})$. If $B_0\subset \Bb_{31}(u,u)$ is an irreducible parabolic subgroup. Then $z_{B_0}$ is a generator of $Z(B_0)$.
\end{prop}
\begin{proof}
Since $\Bb_{31}$ is a connected groupoid, we can consider a morphism $f\in \Bb_{31}(u,u_0)$. The group $B_1:=(B_0)^f$ is a parabolic subgroup of $B(G_{31})=\Bb_{31}(u_0,u_0)$, and we have $z_{B_1}=(z_{B_0})^f$ by Proposition \ref{prop:z_B_0_caractérise_conjugaison}. By Theorem \ref{theo:lattice_of_parabolics_of_B31}, there is some $g\in \Bb_{31}(u_0,u_0)$ so that $B_2:=(B_1)^g=\langle R\rangle$ is generated by one of the sets $R$ of braided reflections given in Theorem \ref{theo:lattice_of_parabolics_of_B31}. We also have $z_{B_2}=(z_{B_1})^g$. It remains to show that $z_{B_2}$ is indeed a generator of the center of $B_2$. For all sets $R$ given in Theorem \ref{theo:lattice_of_parabolics_of_B31}, an expression of a generator of $\langle R\rangle$ as a word in $R$ is known, and we check by direct computations that this element is indeed equal to $z_{\langle R\rangle}$.
\end{proof}

Another consequence of Proposition \ref{prop:z_B_0_caractérise_conjugaison} is that the standardicity of a parabolic subgroup in $\Bb_{31}$ can be checked by looking only at $z_{B_0}$.

\begin{prop}\label{prop:standard_iff_z_B_0_positive}
Let $B_0\subset \Bb_{31}(u,u)$ be a parabolic subgroup. Then $B_0$ is standard if and only if $z_{B_0}$ is positive.
\end{prop}
\begin{proof}
If $B_0=\Gg_\beta(u,u)$ is standard, then $z_{B_0}$ is a positive power of $\delta_\beta(u)$ by definition, so it is positive. Conversely, suppose that $z_{B_0}$ is positive. It is in particular recurrent, and $\PC(z_{B_0})=\SPC(z_{B_0})$ by Theorem \ref{theo:spc_is_pc_for_recurrent}. On the other hand, Proposition \ref{prop:B_0_parabolic_closure_z_B_0}, we have $B_0=\PC(z_{B_0})$. Therefore $B_0=\SPC(z_{B_0})$ is standard.
\end{proof}

Let us then prove that the intersection of parabolic subgroups is a parabolic subgroup. The proof proceeds by, considering two parabolic subgroups $B_1,B_2\subset \Bb_{31}(u,u)$, constructing an element $\alpha\in B_1\cap B_2$ such that $B_1\cap B_2=\PC(\alpha)$. The element $\alpha$ should be chosen to be maximal in some sense. The following definition can be thought of as a kind of rank of an endomorphism in $\Bb_{31}$.

\begin{definition}
For every endomorphism $\gamma$ of $\Bb_{31}$ we define an integer $r(\gamma)$ as follows. Conjugate $\gamma$ to a recurrent endomorphism $\gamma'$. Let $\Gg_\beta$ be the standard categorical parabolic closure of $\gamma'$. Then let $r(\gamma)$ be the length in $\Bb_{31}$ of $\delta_\beta(u)$, where $u$ is the source of $\gamma$.
\end{definition}

By direct examination, we have that, if $\Gg_\beta(u,u)$ and $\Gg_{\beta'}(v,v)$ are two conjugate standard parabolic subgroups of $\Bb_{31}$, then the lengths of $\delta_\beta(u)$ and of $\delta_{\beta'}(v)$ are equal. Thus the integer $r(\gamma)$ does not depend on the choice of a recurrent conjugate $\gamma'$ of $\gamma$.

We are now equipped to prove \cite[Theorem 1.2]{paratresses} in the case of $B(G_{31})$. Again our proof relies on similar arguments to those of \cite[Theorem 6.11]{paratresses}, adapted to our context.

\begin{theo}\label{theo:intersection_parabolic_subgroups}
Every intersection of a family of parabolic subgroups of $B(G_{31})$ is a parabolic subgroup of $B(G_{31})$.
\end{theo}
\begin{proof}
The argument in \cite[Section 6.3]{paratresses} shows that it is enough to show that the intersection of two parabolic subgroups $B(G_{31})$ is still parabolic.

We choose once and for all $u_0\in \Ob(\Bb_{31})$ so that $B(G_{31})=\Bb_{31}(u_0,u_0)$. First, we can assume that $B_1\cap B_2\neq \{1_u\}$, otherwise the result is trivial. Then, we consider a nontrivial element $\alpha\in B_1\cap B_2$ such that $r(\alpha)$ is maximal, and we will show that $B_1\cap B_2=\PC(\alpha)$. 

Up to conjugation of $\alpha,B_1$ and $B_2$ by the same element, we can assume that $\alpha$ is a recurrent endomorphism. Hence, if $u$ is the source of $\alpha$, and $B_0:=\Gg_{\beta}(u,u)=\PC(\alpha)$, we need to show that $B_0=B_1\cap B_2$.

Since $\alpha$ belongs to $B_1$ and to $B_2$, it follows that $B_0=\PC(\alpha)\subset B_1\cap B_2$. So we have one inclusion. Notice that this also yields $z_{B_0}\in B_1\cap B_2$.

Let us then show the inclusion $B_1\cap B_2\subset B_0$. Let $w\in B_1\cap B_2$. In order to show that $w\in B_0$, we will consider the sequence of element $\gamma_m=w(z_{B_0})^m\in B_1\cap B_2$, for $m>0$. We will first show, using a similar argument as in \cite{paratresses}, that $\gamma_m$ is conjugate to a positive endomorphism, for $m$ big enough.

Consider the reduced left-fraction decomposition $w=a^{-1}b$, where $a$ is a composition of $r$ simple morphisms and $b$ is a composition of $s$ simple morphisms. Moreover, if we write $z_{B_0}=\delta_\beta(u)^k$, we have that $z_{B_0}$ is a product of $k$ simple morphisms. We have that $\gamma_m=a^{-1}b(z_{B_0})^m$, so the denominator $D_L(\gamma_m)$ is the composition of at most $r$ simple morphisms, and the numerator $N_L(\gamma_m)$ is the product of at most $s+km$ simple morphisms. Notice that these numbers could decrease, but not increase, if one applies iterated swaps to $\gamma_m$.

We can now conjugate $\gamma_m\in B_1\cap B_2$ to a recurrent element $\ttilde{\gamma}_m$ by iterated swaps. Let $\Gg_{\beta_m}$ be the standard categorical parabolic closure of $\ttilde{\gamma}_m$ and denote $n:=r(\alpha)$. By definition of $\alpha$, if $u_m$ us the source of $\ttilde{\gamma}_m$, we have $r(\gamma_m)=\ell(\delta_{\beta_m}(u_m))\leqslant r(\alpha)=n$.

As we pointed out, if we consider the reduced left-fraction decomposition $\ttilde{\gamma}_m=x_m^{-1}y_m$, then $x_m$ is the product of at most $r$ simple elements, and $y_m$ is the product of at most $s+km$ simple elements. Since $\ttilde{\gamma}_m$ is recurrent, we have $x_m,y_m\in \SCPC(\ttilde{\gamma}_m)$. It follows that the simple factors in the normal forms of $x_m$ and $y_m$ have length at most $r(\gamma_m)\leqslant n$. Let $N_m$ be the number of simple factors in the normal form of $y_m$, whose length is smaller than $n$. Since $y_m$ has at most $s+km$ simple factors in its normal form, we have that $\ell(y_m)\leqslant (s+km)n-N_m$. Thus, we have $\ell(\ttilde{\gamma}_m)\leqslant (s+km)n-N_m$. Since $\ttilde{\gamma}_m$ is a conjugate of $\gamma_m$ in $\Bb_{31}$, we have $\ell(\ttilde{\gamma}_m)=\ell(\gamma_m)=\ell(w)+knm$. We then have $\ell(w)+knm\leqslant (s+km)n-N_m$, so $N_m\leqslant sn-\ell(w)$, where the right-hand side of the inequality does not depend on $m$. On the other hand, since $\ell(\gamma_m)=\ell(w)+nkm$ goes to infinity as $m$ goes to infinity, we have that $\ell(y_m)$ also goes to infinity as $m$ goes to infinity. In particular the number of simple factors in the decomposition of $y_m$ is unbounded. This all implies, in particular, that there exists $M>0$ such that, for every $m>M$, some factor in the left normal form of $y_m$ has length $n$. Since the simple elements in $\Cc_{\beta_m}$ have length at most $n$, it follows that some factor in the left normal form of $y_m$ has the form $\delta_{\beta_m}(u)$, and that $\ell(\delta_{\beta_m}(u))=n$. But, as $x_m\in\Cc_{\beta_m}$, we have $x_m=1$ since otherwise there would be cancellation between $x_m$ and $y_m$. Hence, for $m$ big enough $(m>M)$, we have $\ttilde{\gamma}_m=y_m\in \Cc_{31}$.

Notice that, for every $m>M$, the left normal form of $\ttilde{\gamma}_m$ is $\delta_{\beta_m}(u_m)^{m-N_m}s_1\cdots s_{N_m}$. Let us denote $R_m=s_1\cdots s_{N_m}$ the non-Delta part of the left normal form of $\ttilde{\gamma}_m$. Since $N_m$ is bounded above by a number independent on $m$, it follows that the sequence $\{R_m\}_{m\geqslant 1}$ can take a finitely many possible values.

Also, for every $m>M$, let $c_m$ be the $\preccurlyeq^\Lsh$-minimal morphism in $\Cc_{31}$ such that $c_m\gamma_m c_m^{-1}$ is positive. We know from Remark \ref{rem:rapidité_swap} that $c_m$ is precisely the conjugating element for iterated swaps, hence $c_m\gamma_m c_m^{-1}=\ttilde{\gamma}_m$. By \cite[Proposition VIII.2.5]{ddgkm}, there exists a conjugating element from $\gamma_m$ to a positive element, whose length is bounded above by $\inf(\gamma_m)n$. Since $\inf(\gamma_m)=N_m$ is bounded by a integer not depending on $m$, it follows that the length of the positive morphism $c_m$ is bounded above by an integer not depending on $m$. That is, the sequence $\{c_m\}_{m\geqslant 1}$ can take a finite number of possible values.

Let $e>0$ be such that, for all standard parabolic subgroupoid $\Gg_b$ of $\Bb_{31}$, we have $(\varphi_\beta)^e=1_{\Gg_\beta}$ (this exists since there is a finite number of standard parabolic subgroupoids). Since the sequences $\{R_{em}\}_{m\geqslant 1}$ and $\{c_{em}\}_{m\geqslant 1}$ can take finitely many possible values, and there are only a finite number of standard parabolic subgroupoids, there exists integers $m_1,m_2$, with $M<m_1<m_2$, such that $c_{em_1}=c_{em_2}$, $R_{em_1}=R_{em_2}$ and $\Gg_{\beta_{em_1}}=\Gg_{\beta_{em_2}}$. Let us denote $c:=c_{em_1}=c_{em_2}$, $R=R_{em_1}=R_{em_2}$, $\Gg_\beta=\Gg_{\beta_{em_1}}=\Gg_{\beta_{em_2}}$ and $t=m_2-m_1$. We have $\ttilde{\gamma}_{em_1}=c\gamma_{em_1}c^{-1}$ and:
\[\ttilde{\gamma}_{em_2}=c\gamma_{em_2}c^{-1}=c\beta_{em_1}(z_{B_0})^{et}c^{-1}=(c\gamma_{em_1}c^{-1})(c(z_{B_0})^{et}c^{-1})=\ttilde{\gamma}_{em_1}(c(z_{B_0})^{et}c^{-1}).\]
On the other hand, let us denote $B_3:=\Gg_\beta(u_{em_1},u_{em_1})$. By definition of $e$, we have $\delta_\beta(u_{em_1})=(z_{B_3})^k$ for some $k>0$. Also, since $R=R_{em_1}=R_{em_2}$, it follows that $N:=N_{em_1}=N_{em_2}$. Hence:
\begin{align*}
\ttilde{\gamma}_{em_2}=\delta_{\beta}(u_{em_1})^{em_1-N}R&=\delta_{\beta}(u_{em_1})^{em_2-em_1}\delta_{\beta}(u_{em_1})^{em_1-N}R\\
&=\delta_{\beta}(u_{em_1})^{et}\ttilde{\gamma}_{em_1}\\
&=(z_{B_3})^{kt}\ttilde{\gamma}_{em_1}=\ttilde{\gamma}_{em_1}(z_{B_3})^{kt}.
\end{align*}
Therefore, $c(z_{B_0})^{et}c^{-1}=(z_{B_3})^{kt}$. That is, $((z_{B_0})^{et})^c=(z_{B_3})^{kt}$.

Now, notice that $\PC((z_{B_3})^{et})=B_3$ and $\PC((z_{B_0})^{et})=B_0$. Hence, by Proposition \ref{prop:parabolic_closure_conjugate}, $(B_3)^c=B_0$. On the other hand, we know that $\PC(\ttilde{\gamma}_{em_1})=\Gg_{\beta}(u_{em_1},u_{em_1})$. Hence, again by Proposition \ref{prop:parabolic_closure_conjugate}, $\PC(\gamma_{em_1})=\PC((\ttilde{\gamma}_{em_1})^c)=(B_3)^c=B_0$. This implies that $\gamma_{em_1}\in B_0$. Hence, as $w=\beta_{em_1}z_{B_0}^{-em_1}$, we finally obtain that $w\in B_0$ as we wanted to show.
\end{proof}

\subsection{Characterization of adjacency in the curve graph}
As in \cite{paratresses}, we define the \nit{curve graph} $\Gamma$ for $B(G_{31})$ as the graph whose vertices are irreducible parabolic subgroups of $B(G_{31})$, and where two such subgroups $B_1$ and $B_2$ are adjacent if $B_1\neq B_2$ and either $B_1\subset B_2$, $B_2\subset B_1$ or $B_1\cap B_2=[B_1,B_2]=\{1\}$. As in \cite[Section 6.2]{paratresses}, the adjacency in $\Gamma$ is characterized very easily in terms of the elements $z_{B_1}$ and $z_{B_2}$.

\begin{theo}\label{theo:characterization_of_adjacency}
Two irreducible parabolic subgroups $B_1,B_2\subset B(G_{31})$ are adjacent in $\Gamma$ if and only if $z_{B_1}$ and $z_{B_2}$ are distinct and commute.
\end{theo}
\begin{proof}
Suppose that $B_1$ and $B_2$ are adjacent. We have $B_1\neq B_2$ and $z_{B_1}\neq z_{B_2}$ by Lemma \ref{prop:B_0_parabolic_closure_z_B_0}. If $B_1\subset B_2$, since $z_{B_2}$ is central in $B_2$ it follows that $[z_{B_1},z_{B_2}]=1$. If $B_2\subset B_1$ the argument is the same. Finally, if $B_1\cap B_2=[B_1,B_2]=1$, every element of $B_1$ commutes with every element of $B_2$, hence $z_{B_1}$ commutes with $z_{B_2}$ also in this case. 

Conversely, suppose that $z_{B_1}$ and $z_{B_2}$ are distinct and that they commute. Since $z_{B_i}$ is well-defined and depends only on $B_i$ for $i=1,2$, we have $B_1\neq B_2$. Then, as in the proof of \cite[Proposition 6.12]{paratresses}, we simultaneously conjugate $B_1$ and $B_2$ to standard parabolic subgroups in $\Bb_{31}$.

Since $B_1$ is a parabolic subgroup, there is a morphism $f:u\to v$ in $\Bb_{31}$ such that $(B_1)^f$ is a standard parabolic subgroup by Theorem \ref{theo:parabolics_up_to_conj}. We can replace $B_1$ and $B_2$ by $(B_1)^f$ and $(B_2)^f$ to assume that $B_1$ is standard, and that $z_{B_1}$ is positive by Proposition \ref{prop:standard_iff_z_B_0_positive}.

Now consider the left-fraction decomposition $a^{-1}b$ of $z_{B_2}$. By definition, we have $a\wedge b=1$ and, left-multiplying by $a^{-1}$, we get $1\wedge a^{-1}b=1\wedge z_{B_2}=a^{-1}$.

Since $z_{B_1}$ is positive and commutes with $z_{B_2}$, we have that $z_{B_1}=(z_{B_1})^{z_{B_2}}$ is positive, hence recurrent. By Proposition \ref{prop:gcd_in_R(x)}, $(z_{B_1})^{1\wedge z_{B_2}}=(z_{B_1})^{a^{-1}}$ is a recurrent conjugate of $z_{B_1}$. Since $z_{B_1}$ is positive, all its recurrent conjugates are positive by Proposition \ref{prop:properties_of_swap}. It follows that, if we simultaneously conjugate $z_{B_1}$ and $z_{B_2}$ by $a^{-1}$, we replace $z_{B_1}$ by a positive conjugate, and we replace $z_{B_2}$ by $\sw(z_{B_2})$.

The obtained conjugates of $z_{B_1}$ and $z_{B_2}$ satisfy the same initial hypotheses; they commute and the first one is positive. Hence we can iterate the process described above, replacing $z_{B_1}$ by a positive conjugate, $\sw(z_{B_2})$ by $\sw^2(z_{B_2})$ and so on. Since $z_{B_2}$ is conjugate to a positive element by a series of iterated swaps, we finally obtain simultaneous conjugates of $z_{B_1}$ and $z_{B_2}$ which are both positive (and commute).

If we denote by $f$ the common conjugator of $z_{B_1}$ and $z_{B_2}$ that we consider, we have, for $i=1,2$, $B_i':=(B_i)^f=\PC(z_{B_i})^f=\PC((z_{B_i})^f)$. Since $(z_{B_i})^f$ is positive by construction, it follows from Proposition \ref{prop:standard_iff_z_B_0_positive} that $B_1'$ and $B_2'$ are both standard parabolic subgroups.

We can then assume, up to a simultaneous conjugation, that $B_1=\Gg_{\beta_1}(u,u)$ and $B_2=\Gg_{\beta_2}(u,u)$ are standard parabolic subgroups. Let $R_1$ (resp. $R_2$) be the set of atomic loops of $\Bb_{31}(u,u)$ which are contained in $B_1$ (resp. in $B_2$). We have $\langle R_1\rangle=B_1$ and $\langle R_2\rangle=B_2$ by Proposition \ref{prop:atomic loops generate}. By listing the finite number of standard parabolic subgroups of $\Bb_{31}(u,u)$, we show through direct computations that $z_{B_1}$ and $z_{B_2}$ commute if and only if we have $R_1\subset R_2$, or $R_2\subset R_1$, or $xy=yx$ for all $(x,y)\in R_1\times R_2$. This terminates the proof.
\end{proof}

\printbibliography
\end{document}